\documentclass[12pt, a4paper,leqno]{amsart}
\usepackage{amsmath,amsthm,amscd,amssymb,amsfonts, amsbsy}
\usepackage{latexsym}
\usepackage{txfonts}
\usepackage{exscale}

\usepackage[colorlinks,citecolor=red,pagebackref,hypertexnames=false]{hyperref}

\calclayout
\allowdisplaybreaks

\numberwithin{equation}{section}

\usepackage{pgf}
\usepackage{color}

\newcommand{\nn}{\nonumber}

\theoremstyle{plain}
\newtheorem{theorem}[equation]{Theorem}
\newtheorem{lemma}[equation]{Lemma}

\newtheorem{proposition}[equation]{Proposition}

\newtheorem{lem}[equation]{Lemma}
\newtheorem{pro}[equation]{Proposition}
\newtheorem{thm}[equation]{Theorem}

\theoremstyle{definition}

\newtheorem{definition}[equation]{Definition}

\theoremstyle{remark}
\newtheorem{remark}[equation]{Remark}
\newtheorem{rem}[equation]{Remark}

\newcommand{\ms}{\medskip}

\newcommand{\R}{\mathbb{R}}

\newcommand{\bD}{\mathbb{D}}
\newcommand{\bZ}{\mathbb{Z}}

\newcommand{\bP}{\mathbb{P}}

\renewcommand{\H}{\mathcal H}
\newcommand{\Int}{\mathrm{int}\,}

\renewcommand{\d}{\partial}
\newcommand{\dist}{\,\mathrm{dist}\,}
\newcommand{\sm}{\setminus}
\newcommand{\supp}{\mathrm{supp}}
\newcommand{\diam}{\mathrm{diam}}

\newcommand{\wt}{\widetilde}
\newcommand{\wh}{\widehat}
\newcommand{\ol}{\overline}
\newcommand{\ub}{\underbar}

\newcommand{\cF}{\Theta}
\newcommand{\cP}{{\mathcal P}}
\newcommand{\cQ}{{\mathcal Q}}
\newcommand{\cJ}{{\mathcal J}}
\newcommand{\cR}{{\mathcal R}}
\newcommand{\1}{{\mathbf 1}}

\newcommand{\RR}{{\mathbb{R}}}
\newcommand{\NN}{{\mathbb{N}}}

\newcommand{\CC}{{\mathbb{C}}}
\newcommand{\eps}{\varepsilon}

\newcommand{\bp}{\noindent {\it Proof}.\,\,}
\newcommand{\ep}{\hfill$\Box$ \vskip 0.08in}

\newcommand{\dint}{\int\!\!\!\int}

\newcommand{\dd}{\mathbb{D}}

\newcommand{\C}{\mathcal{C}}
\newcommand{\po}{\partial\Omega}

\newcommand{\A}{\mathcal{A}}
\newcommand{\F}{\mathcal{F}}

\newcommand{\M}{\mathcal{M}}
\newcommand{\W}{\mathcal{W}}

\newcommand{\Ak}{\mathfrak{A}}
\newcommand{\Ab}{\mathbb{A}}
\newcommand{\Jac}{{\rm{Jac}}}

\newcommand{\kF}{\mathcal{F}}
\newcommand{\km}{\mathfrak{m}}
\newcommand{\ka}{\mathfrak{a}}

\renewcommand{\P}{\mathcal{P}}

\renewcommand{\emptyset}{\o}

\def\div{\mathop{\operatorname{div}}\nolimits}

\begin{document}

\title[Harmonic measure on low-dimensional UR sets]{Harmonic measure is absolutely continuous with respect to the Hausdorff measure on all low-dimensional uniformly rectifiable sets}
\author{G. David and S. Mayboroda.}
\newcommand{\Addresses}{{
  \bigskip
  \vskip 0.08in \noindent --------------------------------------
\vskip 0.10in

  \footnotesize

  G.~David, \textsc{
  Universit\'e Paris-Saclay, CNRS, Laboratoire de mathématiques d'Orsay, 91405, Orsay, France} 
\par\nopagebreak
  \textit{E-mail address}: \texttt{Guy.David@universite-paris-saclay.fr }
 \medskip

S.~Mayboroda, \textsc{School of Mathematics, University of Minnesota, 206 Church St SE, Minneapolis, MN 55455 USA}\par\nopagebreak
  \textit{E-mail address}: \texttt{svitlana@math.umn.edu}
}}
\date{}

\thanks{G. David was partially supported by the European Community H2020 grant GHAIA 777822, and the Simons
Foundation grant 601941, GD. S. Mayboroda was partly supported by the NSF RAISE-TAQS grant DMS-1839077 and the Simons foundation grant 563916, SM. }
\maketitle

\ms\noindent{\bf Abstract.} It was recently shown that the harmonic measure is absolutely continuous with respect to the Hausdorff measure on a domain with an $n-1$ dimensional uniformly rectifiable boundary, in the presence of now well understood additional topological constraints. The topological restrictions, while mild, are necessary, as the counterexamples of C. Bishop and P. Jones show, and no analogues of these results have been available for higher co-dimensional sets. 

In the present paper we show that for any $d<n-1$ and for any domain with a $d$-dimensional uniformly rectifiable boundary the elliptic measure of an appropriate degenerate elliptic operator is absolutely continuous with respect to the Hausdorff measure of the boundary. There are no 
topological or dimensional restrictions contrary to the aforementioned results. 

\ms\noindent{\bf R\'esum\'e en Fran\c cais.}
On sait que la mesure harmonique associ\'ee \`a  un domaine de $\R^n$ dont a fronti\`ere est
uniform\'ement rectifiable de dimension $n-1$ est absolument continue par rapport \`a la mesure de surface, sous des conditions topologiques r\'ecemment bien comprises. 
Ces conditions, bien que faibles, sont n\'ecessaires, comme l'ont montr\'e des contre exemples de 
C. Bishop and P. Jones. On ne disposait pas jusqu'ici de r\'esultats analogues lorsque la fronti\`ere est 
de codimension plus grande.

On d\'emontre dans cet article que lorsque la fronti\`ere est uniform\'ement rectifiable 
de dimension $d < n-1$, la mesure elliptique associ\'ee \`a des op\'erateurs elliptiques
d\'eg\'en\'er\'es appropri\'es est absolument continue par rapport \`a la mesure de Hausdorff,
sans avoir besoin de condition topologique suppl\'ementaire.

\ms\noindent{\bf Key words/Mots cl\'es.}
Harmonic measure, elliptic measure, uniform rectifiability, domains with large co-dimensional boundaries/
mesure harmonique, mesure elliptique, uniforme rectifiabilit\'e, domaines \`a fronti\`ere de 
grande co-dimension.

\ms\noindent

\tableofcontents

\section{Introduction}
\label{Sintro} 

Spectacular achievements of the past 20 years at the interface of harmonic analysis, geometric measure theory, and PDEs have finally identified the necessary and sufficient conditions for the absolute continuity of harmonic measure with respect to the Hausdorff measure of an $(n-1)$-dimensional set. In some very informal terms, the problem is as follows. 
The harmonic measure 
of a subset $E$ of the boundary of a domain $\Omega$, $\omega^X(E)$, is the probability that a Brownian traveler, starting at $X\in \Omega$, would exit through the set $E$ rather than its complement. The celebrated 1924 Wiener criterion has identified all boundary points where the harmonic functions are continuous and hence, the harmonic measure is classically well-defined. However, the quantitative information, that is, the question whether the resulting probability is reasonably related to the Hausdorff measure of the set $E$, in other words, whether the Brownian travelers see the portions of the boundary in accordance with their size, turned out to be much more delicate. In PDE terms it is analogous to the question whether the Dirichlet boundary value problem is well-posed with the $L^p$ (rather than continuous) data, with the appropriate dependence of solutions on the $L^p$ size of the data on the boundary \cite{MZ}. 

It is quite remarkable that the key geometric notion in this context was identified already in 1916, 
when F. and M. Riesz  proved that the harmonic measure is absolutely continuous with respect to the Lebesgue measure in a simply connected planar domain bounded by a rectifiable curve \cite{Rfm}. Rectifiability is the property that the set can be covered by a countable collection of Lipschitz graphs, 
modulo a subset of measure zero. 
Extending this result to higher dimensions took more than a century and a development of harmonic analysis, singular integrals, 
and corona decomposition techniques on uniformly rectifiable sets.
We do not aim to provide a detailed overview in this introduction, but let us mention that 
the key milestones were perhaps 
Dahlberg's treatment of Lipschitz domains in \cite{Da}, then 
results on 
2-sided and 1-sided NTA domains with uniformly rectifiable boundaries in \cite{JK} and \cite{HM}, and then, finally, the 
discovery of necessary and sufficient geometric conditions  that
were recently identified in \cite{AHMMT}. One of the main 
problems was that the uniform rectifiability is not sufficient \cite{BJ}, for, in addition, the domain has to exhibit some 
quantitative connectedness, and the exact, sharp form of this connectedness condition seemed elusive for many years. It is interesting to point out that 
the converse result for 
this seemingly PDE question was ultimately resolved thanks to the big advancements in 
the  understanding of  
singular integrals and other harmonic analysis objects on the uniformly rectifiable sets, for instance, the resolution of the David-Semmes conjecture regarding the $L^2$ boundedness of Riesz transforms in \cite{NTV}.

However, all of these results have been restricted to the sets with $n-1$ dimensional boundaries, and with very few exceptions, 
up to recently virtually no theorems treated the higher co-dimensional case
such as the complement of a curve in $\RR^3$.
In \cite{DFM2} the authors of this paper, together with J. Feneuil, have launched a 
program devoted to the degenerate elliptic operator $L = -\div A \nabla$ 
whose matrix of coefficients $A$ has eigenvalues 
roughly proportional to
$\dist (\cdot, \po)^{n-d-1}$, a power of the distance to the boundary. 
Here the additional weight takes the dimension of the boundary
into account, and something like this is 
needed because the usual harmonic functions do not ``see" the lower dimensional sets.
The authors proved the existence of the associated elliptic measure $\omega_L$, together with
fundamental properties such as the H\"older continuity of solutions, the maximum principle,
the doubling property for $\omega_L$, the comparison principle, and estimates for the Green function. See \cite{DFM1, DFM2}, and more recently \cite{DFM-mixed} for an even more 
general setting of domains with boundaries of mixed dimension that will be used here. 
They also proved in \cite{DFM3} a first absolute continuity result for $\omega_L$ as 
in the theorem below, but restricted to the special case of Lipschitz graphs with small constants.

The present paper culminates this line of research by establishing the $A_\infty$ property (quantitative absolute continuity) of the elliptic measure with respect to the surface measure in the most general setting of uniformly rectifiable domains.

\begin{theorem}\label{main-intro} Let $E$ be a $d$-dimensional  uniformly rectifiable set in $\RR^n$, $d\leq n-2$,  and $\mu$ be a uniformly rectifiable measure on $E$. Let $\omega$ be the harmonic measure associated to the operator 
$L=-\div D_{\mu, \alpha}^{-(n-d-1)} \nabla$ in $\RR^n\setminus E$,  
with the distance function
\begin{equation} \label{eq9.40-intro}
D_{\mu, \alpha}(X) 
= \Big\{ \int_E |X-y|^{-d-\alpha} d\mu(y) \Big\}^{-1/\alpha}, \quad \alpha>0.
\end{equation}
Then $\omega$ is 
absolutely continuous with respect to $\mu$,
and its density is weight in the Muckenhoupt class $A_\infty(\mu)$.
\end{theorem}

We chose to state the theorem for unbounded (uniformly rectifiable) sets, but the case of bounded sets would work as well, with minor modifications.

All the 
notions in the theorem 
are customary quantifiable analogues of the properties we just discussed (rectifiability and absolute continuity), but for completeness let us nonetheless carefully recall the definitions.

We say that the closed set $E \subset \R^n$ is Ahlfors regular of dimension $d$ 
when there is an accompanying Ahlfors regular measure $\mu$ on $E$,
which means that
\begin{equation} \label{a1}
C_0^{-1} r^d \leq \mu(E\cap B(x,r)) \leq C_0 r^d
\ \text{ for $x\in E$ and } r> 0.
\end{equation}
We know that $\mu$ is then equivalent to $\H^d_{\vert E}$, 
the restriction to $E$ of the Hausdorff measure,
but we prefer to keep the flexibility of choosing a different $\mu$. The constant $C_0$ will be sometimes referred to as the AR constant of $\mu$.

The notion of uniform rectifiability  was officially introduced by David and Semmes in \cite{DS}.  
One of the many equivalent definitions (see Chapter I.1 of \cite{DS}) is the following.

\begin{definition}\label{def-ur} A $d$-dimensional Ahlfors regular measure $\mu$ 
is uniformly rectifiable if there exist $\theta, C_1 > 0$ so that, for each $x\in \supp (\mu)$ and $r > 0$, there is a Lipschitz mapping $g$ from the $d$-dimensional ball 
$B_d(0, r) \subset \RR^d$ 
to $\RR^n$ such that $g$ has Lipschitz norm less than or equal to $C_1$ and
$$\mu (B(x, r)\cap g(B_d(0, r))) \geq \theta r^d. $$
That is,  $\supp(\mu)$ contains ``big pieces of Lipschitz images of $\RR^d$". An Ahllfors regular set $E\in \RR^n$ is called uniformly rectifiable if $\H^d|_E$ is uniformly rectifiable. We refer to $C_1$ and $\theta$ as the UR constants of $\mu$ (or simply the UR constants of $E$). This is a quantified analogue of the classical concept of rectifiability of a set, as discussed above. 
 \end{definition}
 
Finally, the $A_\infty$ absolute continuity of the 
elliptic (harmonic) measure  with respect to the Hausdorff measure is defined as follows.

\begin{definition}\label{Ainfty} Let $E$ be a $d$-dimensional Ahlfors regular set in $\RR^n$,
$\mu$ be an Ahlfors regular measure on $E$, 
and set $\Omega=\RR^n\setminus E$.
Define an elliptic operator $L$ as above, and denote by $\omega^X$ the corresponding 
elliptic measure with pole at $X$. We say that the elliptic measure $\omega$ is 
$A_\infty$-absolutely continuous
with respect to $\mu$ if for every choice of $\tau_0 \in (0,1)$ and $\epsilon \in (0,1)$, 
there exists $\delta \in (0,1)$, such that for each choice of $x\in E$, $r>0$, a Borel set $F\subset B(x,r)\cap E$, and a corkscrew point $X = A_{x,r}$ 
(i.e., chosen as in \eqref{Corkscrew} below),  
\begin{equation} \label{AiTh}
\frac{\omega^X(F)}{\omega^X(B(x,r)\cap E)} < \delta \, \Rightarrow \, 
\frac{\mu(F)}{\mu(B(x,r)\cap E)} < \epsilon.
\end{equation}
\end{definition}

The notion is tiny bit refined compared to the notion of a (single) $A_\infty$ weight, because
$\omega$ is actually a family of probability measures parameterized by $X \in \Omega$,
and our definition accounts for this in a standard way, and as in the classical case this $A_\infty$
property self-improves into a seemingly stronger condition that says that the 
harmonic measure and $\mu$ are virtually a power of each other, in the sense of Definition~\ref{Ainfty-bis} below. We refer to \cite{Jou, Rub} for general information on 
Muckenhoupt weights, and to 
Section~\ref{Sdef} for the special case of harmonic measure. 

In the conclusion of Theorem~\ref{main-intro}, we also obtain, naturally, 
that $\delta$ depends on $\tau_0, \eps$ and the AR and UR constants of the set $E$ only.

\ms 
Let us discuss Theorem~\ref{main-intro} in more detail. 
Even though the motivation has come from by now ``classical" work in domains with $(n-1)$ dimensional boundaries, the result itself and its proof are actually different, perhaps surprisingly stronger, than their classical counterparts.  

In contrast with the classical case of co-dimension $1$, there is no need for an additional 
topological connectedness conditions here.  The lower dimensional nature of $E$ takes 
care of the topology, 
our boundary is so small that $\Omega = \R^n \sm E$ is sufficiently connected, 
and contrary to Bishop-Jones counterexamples, in this setting we prove that the uniform rectifiability by itself is sufficient for the absolute continuity of harmonic measure. 

Even more intriguing is the situation with the converse. 
Analogous results for traditional co-dimension one boundaries suggest that Theorem~\ref{main-intro} is of the nature of the best possible, i.e., rectifiability should be necessary and sufficient for the absolute continuity of harmonic measure with respect to the Hausdorff measure of the boundary. In our setting, however, there is a surprising ``magical" counterexample. 
When $d+\alpha = n-2$ (which forces $d < n-2$), it turns out that the distance function 
in \eqref{eq9.40-intro} is, in fact, the Green function with the pole at infinity, and
$\omega$ is automatically $A_\infty$ with respect to the Hausdorff measure 
{\it on any Ahlfors regular set}, even when $d$ is not an integer. See \cite{DEM}. 
To the best of our knowledge, this is essentially a unique case where
one can explicitly derive the Green function for an arbitrary Ahlfors regular set. At this point  we tend to believe that it is a miraculous algebraic cancellation and for other values of $\alpha$ the $A^\infty$ property of harmonic measure implies uniform rectifiability, i.e., the condition in Theorem~\ref{main-intro} will be proved to be necessary and sufficient.

Returning to the discussion of the statement of Theorem~\ref{main-intro}, observe that
the higher co-dimensional setting requires a rather peculiar choice for the operators $L$ above. 
The fact that the coefficients are roughly  proportional to a certain power of distance to the boundary is almost a necessity, dictated
by the scaling considerations, Sobolev embeddings, etc. 
The usual Laplacian would not work, because the harmonic functions do not even see 
sets of dimension $d\leq n-2$.
However, what is perhaps more surprising, working with the conventional Euclidean distance would ruin our proof of absolute continuity for the elliptic 
elliptic measure,  even on a small Lipschitz graph -- see the discussion in \cite{DFM3}. 
The distance $D_{\mu,\alpha}$ of \eqref{eq9.40-intro} 
turns out to be a correct substitute, smoothing out appropriately at all scales. 
It may also have some special algebraic properties: as we mentioned above, for $\alpha = n-2-d$ 
it actually coincides with the Green function with the pole at infinity. It appears to be a powerful and perhaps still not completely understood version of the distance function used in geometric analysis -- see \cite{DEM}. 

These observations lead to a question:  what is the range of the operators for which one could establish an analogue of Theorem~\ref{main-intro}? 
Certainly not every degenerate elliptic operator for which $A(X)$ has roughly the size
of $\dist(X, E)^{-n+d+1}$ will give an absolutely continuous elliptic measure. Even in the classical
case of co-dimension $1$, there are many 
counterexamples (see \cite{MM, CFK}), and most absolute continuity
results concern operators with a special form, or which are small perturbations of the Laplacian.
Here $L= - \div D_{\mu, \alpha}^{-(n-d-1)} \nabla$ plays the role of the Laplacian,
and in this case too Theorem~\ref{main-intro} also holds for any 
elliptic operator which is a Carleson perturbation of $L$; see \cite{Bruno}.
That is, the set of ``good" operators is ultimately quite large, comparable to the classical 
scenario.

Let us now discuss the proof. A brutal attempt to adapt the ``classical" approach to lower dimensional boundaries collapses spectacularly. The general principle, that we still want to follow,
is to start from \cite{DFM3}, which provides the desired absolute continuity result when
$E$ is a small Lipschitz graph, and then use approximations and stopping time arguments to
extend this result to sawtooth domains and then to general 
uniformly rectifiable sets (in the classical case, we would need some additional connectivity).
The machinery that allows this is very beautiful, and also quite intricate, but let us highlight at least two big problems which will make the situation different.

The main obvious one is that the principal
engine of the construction is a comparison 
between domains. That is, we like to consider $\Omega$, and hide parts of the boundary that we do not control outside of sawtooth domains, and for the intersection of $\Omega$ with the sawtooth domain, use the maximum principle to relate the two corresponding elliptic measures.
Here we manage to construct better approximating domains whose boundaries 
coincide with $E$ in some places, but how are we going to hide the rest of the boundary and
compare the harmonic measures on the intersection?

Another unpleasant issue is that in the classical case, we typically deal with the Laplacian, which
has a local definition, unlike our operator $L$ whose coefficients depend on $E$, and even tend
to infinity along the boundary. In other words, all the domains, including the approximating 
ones, come with their own operator $L$, which is not local either, and which typically carries a ``memory" of the original domain, and of course we will need to be more careful about which operator we take when we make comparisons.

After many different attempts, we decided that the possibility to hide a bad piece of
$E$ behind a sawtooth boundary was too important to be avoided, and this forced us
to consider domains with boundaries of mixed dimensions, i.e., where some part of the 
boundary would be $(n-1)$-dimensional (and could be used to hide parts of $E$ that we
do not control), and other parts would be $d$-dimensional. Then we had to adapt the 
theory of degenerate elliptic operators to such domains, and in particular understand how to
relate the size of our coefficients to the mass of a doubling measure on the boundary.
For the boundary of the sawtooth regions, for instance, it helps that the coefficient of
$L$ is related to the distance to $E$, but is not necessarily singular along the sawtooth.
Fortunately the reader will not have to deal with the extension of the elliptic theory here, 
because most of it was taken care of in \cite{DFM-mixed}.

Even so, we also have to replace, rather than shield, ``bad" portions of the set to create a better, Reifenberg flat surface. This construction takes a good part 
of the present paper and heavily relies on the parametrization of Reifenberg flat domains in \cite{DT} and on the $A_\infty$ 
results for the elliptic measure on lower dimensional small Lipschitz graphs that we proved in \cite{DFM3}. Such a replacement is performed at all scales which eventually have to be glued via a certain extrapolation argument. Here again, we start from the Hofmann-Martell approach to extrapolation in \cite{HM} but quickly learn that our hypothesis is quite a bit weaker than theirs, requiring a careful re-working of the entire scheme. Finally, last but not the least, is the obstacle arising from the fact that our operator is not local, the coefficients are the powers of distance to the original boundary, and so changing the domain entails changing the operator. One has to prove that such a change, at the correct scales, is not detrimental to the entire enterprise. An experienced reader can take a look at Theorem~\ref{l9.39} and compare it to a (considerably easier) standard version in \cite{DJK}. 

Finally, before we turn to the body of the paper, let us mention that during the preparation of this manuscript an alternative approach to some of its results has been developed in \cite{F}. 

\section{Definitions and preliminaries related to the elliptic theory and properties of weights}
\label{Sdef}

We are given an Ahlfors regular set $E$ of dimension $d$ in $\R^n$ and an accompanying Ahlfors regular measure $\mu$ on $E$. We denote $\Omega=\RR^n\setminus E$.

Before we proceed, we need to say a few words about 
{\it{dyadic pseudocubes}.}
We shall assume that a net of dyadic pseudocubes has been chosen on $E$. We use the cubes
given by \cite{Ch}, except that we will find it more convenient to use scales that are powers of
$10$, because this way we will be closer to the notation of \cite{DT}. We systematically set
\begin{equation} \label{a2}
r_k = 10^{-k} \ \text{ for } k \in \bZ,
\end{equation}
and we suppose we chose a collection $\bD = 
{\dd(E)=} \cup_{k \in \bZ} \bD_k$
(our ``dyadic cubes''),  with the usual properties. We will in particular use the facts that
for each cube $Q\in \bD_k$, there is a center {$x_Q$} such that 
\begin{equation} \label{a3}
E \cap B(x_Q,C^{-1}r_k) \subset Q \subset B(x_Q, r_k),
\end{equation}
and that the different cubes $Q \in \bD_k$ are disjoint, that when $k \leq l$,
$Q \in \bD_k$, and $R \in \bD_l$, then either $R \subset Q$ or else
{$R\cap Q=\emptyset$}.

To get this, we merely assume that the cubes $Q$ are Borel sets,
but there is a ``small boundary condition'', that we shall not use in full before
the last Carleson estimates for the control on different distances, that implies in
particular that
\begin{equation} \label{a4}
\mu(\ol Q \sm Q) = 0 \ \text{ for } Q\in \bD.
\end{equation}
We shall often denote by $k(Q)$ the generation of $Q$ (i.e., the integer $k$ such that 
$Q\in \bD_k$), and set $l(Q) = r_{k(Q)} = 10^{-k(Q)}$. 
{For any $Q\in \dd(E)$ we let 
$$\dd(Q)=\dd_Q=\{Q'\in \dd(E):\, Q'\subseteq Q\}.$$}

\begin{definition}\label{cs} The corkscrew points for $\Omega$  are points $A_{x,r} \in \Omega$,
associated to $x\in E$ and $r > 0$, such that (for some constant $\tau_0 > 0$)
\begin{equation} \label{Corkscrew}
\tau_0 r \leq \dist(A_{x,r},E) \leq |A_{x,r}-x| \leq r.
\end{equation}
Corkscrew points exist for all $x, r$
whenever $E$ is any Ahlfors regular set of dimension $d < n-1$, and we can take $\tau_0$
to depend only on $n$, $d$, and $C_0$ from \eqref{a1}; 
see Lemma 11.46 in \cite{DFM2}.

It will also be convenient to use corkscrew points associated to the dyadic decomposition of $E$, and for any $Q\in \dd(E)$ we write $A_Q:=A_{x_Q, C^{-1} l(Q)}$ where the constant $C$ is from \eqref{a3}.
\end{definition}

Since the set $E$ satisfies \eqref{a1}, it enters the scope of the elliptic theory 
developed in \cite{DFM2}. Let us recall some of the main properties that will be needed. 

Let $L = - \div \Ak \nabla$ be a degenerate elliptic operator, in the sense that  $\Ak : \Omega \to \mathbb{M}_n(\R)$ satisfies 
\begin{eqnarray} \label{1.2.1} 
&& \dist(X,\Gamma)^{n-d-1} \Ak(X)\xi \cdot \zeta \leq C_1  |\xi| \,|\zeta|
\ \text{ for } X\in \Omega \text{ and } \xi, \zeta \in \R^n, \\[4pt]
\label{1.2.2}
&& \dist(X,\Gamma)^{n-d-1} \Ak(X)\xi \cdot \xi \geq C_1^{-1}  |\xi|^2
\ \text{ for } X\in \Omega \text{ and } \xi \in \R^n,
\end{eqnarray}
for some $C_1 \geq 1$. We say that $u$ is a weak solution of $Lu=0$, if 
$u \in  W^{1,2}_{loc}(\Omega)$ and
\begin{equation} \label{a2.7}
\int_\Omega \Ak \nabla u\cdot \nabla v = 0 \qquad \forall v\in C^\infty_0(\Omega).
\end{equation}
Here, $W^{1,2}_{loc}(\Omega)$
is the set of functions $u\in L^2_{loc}(\Omega)$ whose derivative 
(in the sense of distribution on $\Omega$) also lies in $L^2_{loc}(\Omega)$.

For each $X \in \Omega$, we can define a (unique)
probability measure $\omega^X=\omega^X_{\Omega, L}$ on $E$, with the property that 
for any bounded measurable function $f$ on $E$, the function $u_f$ defined by 
\begin{equation} \label{a2.6}
u_f(X) = \int_E f(y) d\omega^X(y)
\end{equation}
is a weak solution.
This is only stated in \cite{DFM2} when $f \in C^0_0(E)$ is continuous and 
compactly supported in $E$ (see Lemma 9.30 and (iii) of Lemma 9.23 there,
and also (8.1) and (8.14) for the definitions) 
and when $f$ is a characteristic function of Borel set (see Lemma 9.38 there); 
the general case would not be hard, but we do not need it anyway. It has been proved in \cite{DFM2} that $\omega$ is doubling on the 
Ahlfors regular set $E$, in the sense of Definition~\ref{DDoub}.

There is also a dense subclass on which we can say a little more. 
Denote by ${\mathcal M}(E)$ the set of measurable functions on $E$ and then 
define the Sobolev space
\begin{equation}\label{defH}
H=\dot H^{1/2}(E) :=\left\{g\in {\mathcal M}(E): \,
\int_E\int_E {|g(x)-g(y)|^2 \over |x-y|^{d+1}}d\mu(x) d\mu(y)<\infty\right\}.
\end{equation}
\noindent The class $H \cap C^0_0(E)$ is dense in $C^0_0(E)$
(see about 13 lines above (9.25) in \cite{DFM2} for the proof of density), 
and if $f \in H \cap C^0_0(E)$, the solution $u_f$ defined by \eqref{a2.6}
lies in the Sobolev space $W^{1,2}(\Omega,\dist(X,E)^{d+1-n}dX)$, which means that
\begin{equation} \label{a2.8}
\int_\Omega |\nabla u_f(X)|^2 \dist(X,E)^{d+1-n} dX <+\infty,
\end{equation}
and also 
\begin{equation} \label{a2.9}
\text{$u_f$ has a continuous extension to $\R^n$, which coincides with $f$ on $E$.} 
\end{equation}
See (i) of Lemma 9.23 in \cite{DFM2}, together with its proof eight lines above (9.25).

It should be stressed that since $\omega^X$ is a probability measure, $u_f$
is a nondecreasing function of $f \geq 0$. This is of course a manifestation of the maximum principle.

Recall now the definition of the absolute continuity of the elliptic measure in Definition~\ref{Ainfty}.
Here, $\delta$ depends on $\tau_0$ and $\eps$ as well possibly other parameters, and we shall refer to the latter as 
$A_\infty$ constants  of $\omega_L$. Typically, those include ellipticity constants of $L$ as well as some geometric characteristics of the set. We will try to be prudent below listing them carefully.

Let us temporarily narrow down to $\Omega_0=\RR^n\setminus \RR^d$,
where we shall write the generic point as $X = (x,t)$, with
$x\in \R^d$ and $t\in \R^{n-d}$.
The following theorem is the starting point of our absolute continuity results.

\begin{definition}\label{cm} We say that a matrix-valued (or scalar-valued) measurable function $F$ on $\Omega_0=\RR^n\setminus \RR^d$ satisfies the Carleson measure condition with constant $C$ and write $F\in CM(C)$ or simply $F\in CM$ if there is a constant $C \geq 0$
such that for $X = (x,0)\in \R^d$ and $R>0$,
\begin{equation} \label{d2}
\int_{\Omega_0 \cap B(X,R)} |F(y,t)|^2 \, \frac{dydt}{|t|^{n-d}} \leq C R^d.
\end{equation}
\end{definition}

\begin{thm}\cite{DFM3} \label{Itsf1} 
Let $A_0$ be a degenerate elliptic matrix satisfying \eqref{1.2.1} and \eqref{1.2.2} in 
$\Omega_0 = \R^n \setminus \R^d$, and set then $L_0 = - \div A_0 \nabla$.
Define the rescaled matrix $\A$ by $\A = |t|^{n-d-1} A_0$, so that now 
$\A$ satisfies the usual ellipticity bounds and
$L_0 = - \div |t|^{d+1-n}\A \nabla$, and assume that $\A$ has the following block structure: 
\begin{equation}\label{Ieqsf1}
\A(X)= \left( \begin{array}{cc}
\A^1(X) & \A^2(X) 
\\  \C^3(X) & b(X)I_{n-d}+  \C^4(X)  \end{array} \right), 
\end{equation}
where $\A^1(X)$ is a matrix in $M_{d\times d}$, $\A^2(X)$ is a matrix in 
$M_{d\times (n-d)}$, $b$ is a function on $\Omega_0$, $I_{n-d}$ is the identity 
matrix in $M_{(n-d)\times (n-d)}$, and in addition we can find constants $C \geq 0$ 
and $\lambda \geq 1$ such that
\begin{equation} \label{a2.32}
\lambda^{-1} \leq b \leq \lambda \ \text{ on } \Omega_0,
\end{equation}
\begin{equation} \label{a2.33}
|t|\nabla b \in CM(C),
\end{equation}
\begin{equation} \label{a2.34}
\C^3, \, \C^4 \in CM(C).
\end{equation}
Then the harmonic measure 
$\omega^X_{\Omega_0,L_0}$ 
is $A_\infty$-absolutely continuous with respect to the Lebesgue measure on $\R^d$ (with the Definition~\ref{Ainfty}).
\end{thm}

Definition~\ref{Ainfty} might seem confusing at first because traditionally 
the $A_\infty$ 
property of one weight with respect to another is formulated somewhat differently. However, 
when the weight is doubling, the $A_\infty$-absolute continuity as stated in Definition~\ref{Ainfty} is equivalent to the traditional $A_\infty$ 
property (see, e.g., \cite{CF}, Lemma~5, where it is proved for any two doubling measures). Let us list the main definitions and results to this effect, complemented by their 
dyadic counterparts, as they will be used throughout the paper for a variety of measures. 

\begin{definition}\label{DDoub} Let $E$ be a $d$-dimensional Ahlfors regular set in $\RR^n$ (or more generally, any metric space). A nontrivial measure $\omega$ on $E$ is doubling if for any $x\in E$, $r>0$, 
\begin{equation} \label{aa21}
0<\omega (B(x, 2r)\cap E) <C \omega (B(x, r)\cap E)<\infty,
\end{equation}
with a uniform constant $C$. 

When $E$ possesses a dyadic structure, e.g., 
when $E$ is a $d$-dimensional Ahlfors regular set in $\RR^n$
and we chose a collection of pseudocubes $\bD(E)$ as near \eqref{a2},
we say that $\omega$ is dyadically doubling on $Q_0\subset \dd(E)$ if for every $Q\in \dd(Q_0)$ and every dyadic ``child" of $Q$, $Q'$, there exists a uniform constant $C$ such that 
$$0<\omega (Q) <C \omega (Q')<\infty.$$ 
\end{definition} 

In the particular case when $\omega$ is 
in fact a family $\{ \omega^X \}$ of elliptic measures on $E$, we say that $\omega$ is doubling 
if for any surface ball $\Delta(x, r)=B(x,r)\cap E$, $x\in E$, $r>0$, the harmonic measure with a pole at $A_{x, r}$, $\omega^{A_{x, r}}$, is doubling on on $\Delta(x, r)$, and the constant $C$ in the doubling property is independent of $x, r$. Similar definitions apply in the dyadic case, with the pole at $A_Q$.
Equivalently, one could say that $\omega$ is doubling  when \eqref{aa21} holds for $\omega^X$ as long as $X$ is far enough from $B(x,r)$ (for instance, $X \in \Omega \sm B(x,4r)$), and similarly,  in the definition of dyadically doubling for the family 
$\omega = \{ \omega^X \}$, we would only ask for the doubling
condition when $X$ is far enough from $Q$.

In Lemma 11.102 of \cite{DFM2}, it is proved  that for any $d$-dimensional Ahlfors regular set $E$ in $\RR^n$, $d<n-1$, the harmonic measure of any elliptic operator is doubling and, hence, dyadically doubling on $E$.

Next we say a little more about local versions of the $A_\infty$ condition.
For a single measure $\omega$, we would use the following definition.

\begin{definition}\label{Ainfty-gen} 
Let $E$ be a $d$-dimensional Ahlfors regular set in $\RR^n$ and 
$\mu$ be an Ahlfors regular measure on $E$.
Given any surface ball $\Delta(x, r)=B(x,r)\cap E$, $x\in E$, $r>0$, 
we say that a doubling measure $\omega$ on $E$ is 
$A_\infty$-absolutely continuous 
with respect to 
$\mu$ on $\Delta(x, r)$ (denoted by $A^\infty(\Delta(x, r))$),  
if for every $\epsilon \in (0,1)$, there exists $\delta \in (0,1)$, such that for 
every surface ball $\Delta'=B'\cap E$, $B'\subseteq B(x,r)$ and every Borel set $F\subset \Delta'$
\begin{equation} \label{AiTh-gen}
\frac{\omega(F)}{\omega(\Delta')} < \delta \Rightarrow \frac{\mu(F)}{\mu(\Delta')} < \epsilon.
\end{equation}
Similarly, replacing surface balls by dyadic cubes, we say that a dyadically doubling 
measure $\omega$ is  (dyadically) $A_\infty$-absolutely continuous 
with respect to $\mu$ on $Q\in \dd(E)$ 
(denoted by $A^\infty_\dd(Q)$),   
if for every $\epsilon \in (0,1)$, there exists $\delta \in (0,1)$ such that for every 
$Q'\in \dd(Q)$ and every Borel set $F\subset Q'$, 
\begin{equation} \label{AiTh-gen-2}
\frac{\omega(F)}{\omega(Q')} < \delta \Rightarrow \frac{\mu(F)}{\mu(Q')} < \epsilon.
\end{equation}
\end{definition}

As the reader can guess from Definition~\ref{Ainfty}, in the particular case when 
$\omega = \{ \omega^X \}$ is (a family of) harmonic measures on $E$, 
we say that $\omega$ is $A_\infty$-absolutely continuous with respect to  
if for any surface ball $\Delta(x, r)=B(x,r)\cap E$, $x\in E$, $r>0$, 
the harmonic measure with a pole at $A_{x, r}$, $\omega^{A_{x, r}}$, is 
$A_\infty$-absolutely continuous with respect to $\mu$ on $\Delta(x, r)$, 
and the choice of $\delta$ depends on $\epsilon$ (and $\tau_0$ in the definition of the corkscrew point as well as AR and doubling constants) but not on $x, r$. In fact it would then be possible,
using estimates on $\omega$ (that the reader may find in \cite{DFM2}),
to deduce the same estimates for other poles $X\in \Omega \setminus B(x, 4r)$. All these definitions are equivalent, in the sense that
we get estimates for $\delta$ with different definitions that depend only
on those with the initial definition, our bounds for $E$ and $L$, and the
various corkscrew constants. Similar definitions apply in the dyadic case, with a pole at $A_Q$.

The $A_\infty$ condition is known to imply a stronger form of absolute continuity,
which we define now, starting with the case of a single measure.

\begin{definition}\label{Ainfty-bis} Let $E$ be a $d$-dimensional Ahlfors regular set in $\RR^n$ and $\mu$ be an Ahlfors regular measure on $E$, $\Omega=\RR^n\setminus E$. Given any surface ball $\Delta(x, r)=B(x,r)\cap E$, $x\in E$, $r>0$, we say that a Borel measure 
$\omega$ is strongly  absolutely continuous in $\Delta(x, r)$,
with respect to $\mu$,  if there are positive constants $C$ and $\theta$ such that for every surface ball $\Delta'=B'\cap E$, $B'\subseteq B(x,r)$ and every Borel set $F\subset \Delta'$, 
\begin{equation} \label{eqAinf}
\omega(F)\leq C\,\left(\frac{\mu(F)}{\mu(\Delta')}\right)^{\theta}\omega(\Delta').
\end{equation}

Similarly, replacing surface balls by dyadic cubes, we say that a dyadically 
doubling measure $\omega$ is 
strongly dyadically  absolutely continuous in $Q$,
with respect to $\mu$,
if there are positive constants $C$ and $\theta$ such that for every $Q'\in \dd(Q)$ and every Borel set $F\subset Q'$, 
\begin{equation} \label{eqAinf}
\omega(F)\leq C\,\left(\frac{\mu(F)}{\mu(Q')}\right)^{\theta}\omega(Q').
\end{equation}
\end{definition} 

Much as above, in the particular case when 
$\omega = \{ \omega^X \}$  is in fact a family of harmonic measures on $E$,
we say that $\omega$ is 
strongly absolutely continuous 
with respect to $\mu$ if for any surface ball $\Delta(x, r)=B(x,r)\cap E$, 
$x\in E$, $r>0$, the harmonic measure with a pole at $A_{x, r}$, $\omega^{A_{x, r}}$, 
is strongly  absolutely continuous
with respect to $\mu$, and 
the constants
$C, \theta$ depend on $\tau_0$ in the definition of the corkscrew point as well as AR and doubling constants  but not on $x, r$. Similar definitions apply in the dyadic case, with the pole at $A_Q$.

\begin{rem}\label{rAinfty}
It was proved in \cite{CF}, Lemma~5, 
that for any two doubling measures $\mu$ and $\omega$,
if $\omega$ is $A_\infty$-absolutely continuous with respect to $\mu$, then
$\omega$ is also strongly  absolutely continuous with respect to $\mu$,
and also $\mu$ is also strongly  absolutely continuous with respect to $\omega$.
Also see \cite{Jou} or \cite{Rub}.

The dyadic (and local) analogues of these facts were established, e.g., in \cite{HM}, Appendix B, Remark 2.10. That is, under the definitions above, a dyadically doubling measure $\omega$ is 
strongly dyadically  absolutely continuous  on $Q\in \dd(E)$ with respect to $\mu$ 
if and only if it is dyadically $A_\infty$ absolutely continuous  with respect to $\mu$. In fact, both are equivalent to an (apparently) weaker statement that there exist $0<\epsilon, \delta<1$ such that \eqref{AiTh-gen} (respectively, \eqref{AiTh-gen-2}) hold -- the latter property is referred to as comparability for doubling measures. Moreover, 
$A_\infty$ and its local, dyadic, and strong versions
are equivalence relationships, in the sense that for instance, if $\omega$ is strongly absolutely continuous with respect to $\mu$ with some constants $C>0$, $\theta>0$, 
then $\mu$ is strongly absolutely continuous 
with respect to $\omega$ with some other constants $C'>0$, $\theta'>0$; see \cite{CF}, Lemma~5 for the standard case and \cite{HM}, Lemma B.7 for the dyadic one.

\end{rem}

\section{Preliminary geometric considerations}\label{intro}

In Sections~\ref{intro}--\ref{sSum} we define a correct change of variables, adapted to
a stopping time region associated to a uniformly rectifiable set of integer dimension $d$ in $\R^n$.  To be more precise, for any stopping time region subject to some flatness and regularity constraints we construct a Reifenberg flat set $\Sigma$ which coincides with our initial set $E$ in the ``base" of the sawtooth and which has a nice parametrization, 
in fact coming from  a nice change of variables transforming $\RR^n\setminus \Sigma$ into $\RR^n \setminus\RR^d$.

The change of variables will be inspired by that in \cite{DT}, but unfortunately we need an array of properties which was not explicitly targeted in \cite{DT}. Indeed, we need to use it similarly to the change of variables of \cite{DFM3}, to ensure the absolute continuity of a certain elliptic measure on the underlying set. In both cases, the philosophy is to respect the orthogonal direction to the tangent plane to the boundary set. However, the details are quite different and we will have to devote a considerable effort to the proof that (a slightly modified) construction from \cite{DT} satisfies the desired properties. We try to take notations that are fairly close to those of \cite{DT}, which we shall cite abundantly. To start, let us describe a stopping time region.

We are given an Ahlfors regular set $E$ of dimension $d$ in $\R^n$. In our end-game applications $E$ will be uniformly rectifiable, but we do not need to assume this for the moment. The definition of the stopping time regions will take care of the regularity needed for the first few chapters.

In this section, we are given a \ub{stopping time region} $\cF$, with some definite constraints on
how it is built, and associate to it a few geometric objects. 
There will be a specific way to construct $\cF$ from its top cube $Q_0$, 
but let us keep some latitude, 
without making our life too complicated. So we start from a cube $Q_0$, and without loss
of generality we assume that
\begin{equation} \label{a5}
Q_0 \in \bD_0.\end{equation}
Then $\cF$ will be a subset of $\dd(Q_0)$, the set of subcubes of $Q_0$. 
For $Q \in \cF$ and $0 \leq k \leq k(Q)$, denote by $R_k(Q)$
the cube of $\bD_k$ that contains $Q$; thus $R_k(Q)$ is an ancestor of $Q$
and $R_k(Q) \subseteq Q_0$. 
We demand that $Q_0 \in \cF$ (otherwise, there in no construction to be done)
and that $\cF$ is hereditary, which means that
\begin{equation} \label{a7}
R_k(Q) \in \cF \ \text{ for $Q \in \cF$ and } 0 \leq k \leq k(Q) 
\end{equation}
(i.e., 
if $Q \in \cF$, then all its ancestors between $Q$ and $Q_0$ lie in $\cF$).  

For the remaining properties of $\cF$, we need to choose a large constant
$M \geq 1$, a very small constant $\varepsilon_1 > 0$, and another constant
$\delta_1 > 0$, in practice much larger than $\varepsilon_1$. 
Apparently our construction will not put any constraint on $\delta_1$, except for the
fact that some constants will become very large when $\delta_1$ is large.
We will take $M$ quite large, depending on other geometric constants of the construction,
and then $\varepsilon_1$ will need to be small enough, depending on $n$, $d$, 
the constant $C_0$ in \eqref{a1}, and $M$.
This includes a dependence on our choice of $\bD$ through the constant in \eqref{a3}, 
but we can choose $\bD$ once and for all, with a constant in \eqref{a3} that depends only on $n$, $d$, and $C_0$.

It will simplify our definition if we assume that for each $Q \in \cF$, a $d$-plane
$P(Q)$ has been chosen, with the following properties. First of all, $P(Q)$
is quite close to $E$ near $Q$. That is, if we define a normalized Hausdorff 
distance between sets $d_{x,r}(F,G)$ by
\begin{equation} \label{a8}
d_{x,r}(F,G) = r^{-1} \sup_{y\in F \cap B(x,r)} \dist(y,G)
+ r^{-1} \sup_{y\in G \cap B(x,r)} \dist(y,F),
\end{equation}
then we require that
\begin{equation} \label{a9}
d_{x_Q, M l(Q)}(E,P(Q))  \leq \varepsilon_1 \ \text{ for } Q \in \cF,
\end{equation}
where $l(Q) = 10^{-k(Q)}$ is the official sidelength  
of $Q$ and the center $x_Q$ is as in \eqref{a3}.
We also measure the average distance from points of $E$ near $Q$ to $P(Q)$, and 
encode them into numbers $\beta(Q)$ such that
\begin{equation} \label{a10}
\int_{E \cap B(x_Q, M l(Q))} \dist(y,P(Q))\, d\mu(y)\leq l(Q)^{d+1} \beta(Q),
\end{equation}
where $\mu$ is the measure on $E$ that we started with (but its precise choice does not 
matter here). These numbers are close to 
the $\beta$-numbers of P. Jones
associated to $E$ and computed with $L^1$-norms, but we reserve the right to make 
$\beta(Q)$ larger than the actual number $\beta_1(x_Q,M l(Q))$ 
and choose the $P(Q)$ differently.
Then we define a Jones function $J$ on $\cF$ by setting
\begin{equation} \label{a11}
J(Q) = \sum_{0 \leq k \leq k(Q)} \wt\beta(R_k(Q))^2,
\end{equation}
where as before $R_k(Q)$ is the ancestor of $Q$ which is of generation $k$,
and unfortunately we need need to replace $\beta(Q)$ with a slightly larger, more regular,
function of $Q$, namely
\begin{equation} \label{a11bis}
\wt\beta(R) = \sup\big\{ \beta(S)\, ; \, S \in \cF(k(R)) \text{ and }
\dist(R,S) \leq M l(R) \big\}, 
\end{equation}
where 
\begin{equation} \label{b4}
\cF(k) = \big\{ Q \in \cF \, ; \, l(Q) = r_k \big\}.
\end{equation}
Notice that we may count the same set $R$ twice in \eqref{a11}, if successive ancestors
of $Q$ happen to be given by the same subset of $E$. This is all right, and probably 
even more reasonable. Notice also that replacing $\beta(R)$ with $\wt\beta(R)$ will not cost
us much in practice; we will just need to control $E$ (and possibly $\mu$) on an even larger ball.
Finally observe that $J(Q) \geq J(R)$ when $Q \subset R \subset Q_0$ and 
$k(Q) \geq k(R) \geq 0$. We demand that
\begin{equation} \label{a12}
J(Q) \leq \delta_1 \ \text{ for } Q \in \cF.
\end{equation}

\ms
This completes the list of conditions that we put on $\cF$. 
We do not need to say yet how we produce $\cF$, but the algorithm that will
be used later is as follows. For each cube $Q$, we shall define a quantity $\alpha(Q)$,
for instance using the $\alpha$-numbers coming from \cite{To} 
and choose a plane $P(Q)$ that is nearly optimal in the definition of $\alpha(Q)$.
These numbers will be introduced in Section \ref{SDist}; for the moment we do not need to 
know what they are.

Then we will start from the top cube $Q_0$, and decide to remove a cube
$Q \in \bD(Q_0)$, as well as all its descendants, as soon as $\alpha(Q) > \varepsilon_0$
or 
$$J_\alpha(Q) := \sum_{0 \leq k \leq k(Q)} \alpha(R_k(Q))^2 \geq \delta_0.$$
It will turn out that the numbers $\alpha(Q)$ control the properties \eqref{a9} and \eqref{a12},
in the sense that if $\varepsilon_0$ is chosen small enough, then \eqref{a9} follows from
the fact that $\alpha(Q) < \varepsilon_0$, and similarly \eqref{a12} follows from the 
fact that $J_\alpha(Q) \leq \delta_0$.

\begin{remark} \label{r0} 
There are constraints on $M$ and $\varepsilon_1$. The first ones will come soon,
to verify the CCBP properties at the beginning of the next section, and then there
will be other ones in the last section. Since we want to keep some freedom in the choices,
we announce now that all we need, up to Section \ref{Sf}, is to take $M$ large enough,
and then $\varepsilon_1$ small enough, depending on $M$, $\delta_1$, and 
the other parameters. 

We can let $M$ depend on $\delta_1$ (in fact, we claim below that we could even let
$\delta_1$ be a large number).
 
The relation between $\delta_1$ and $\varepsilon_1$ is more delicate, and we announce it in
advance so that we cannot be suspected of cheating.
Both constants will be small in our argument, and correspond to stopping time conditions.
The basic reason for stopping in our geometric construction of a parameterization
is when the set starts being flat enough, and $\varepsilon_1$ corresponds to the minimal 
amount of flatness that we demand. 
We will choose $\varepsilon_1$ last, possibly depending on 
the other parameters. Now we also want to control the bilipschitz constants for our approximations,
and we use the Jones function to do this. The role of $\delta_1$ is to
control the Jones function, and then the bilipschitz constants for our mappings.
In a sense, $\varepsilon_1$ acts like the $L^\infty$ norm of some quantity (the $\beta$-numbers)
that needs to stay small, and $\delta_1$ like the $L^\infty$ norm of some integral, or sum, 
of some related (but different) quantity (the $\alpha$-numbers).

We promise that we will not let $\varepsilon_1$ depend on $\delta_1$, 
because this would contradict the spirit of stopping times, but we will nonetheless
do an offense to that spirit, because in some argument, and for the sake of laziness,
we will use $\delta_1$ to control some quantity that should be in fact be controlled by
$\varepsilon_1$ in a cleaner (but longer) argument. Because of this, we will
require $\delta_1$ to be small, but a real purist would allow it to be large too,
and this would create a more interesting parameterization when we only stop when this is really
needed. Formally speaking, we could also take $\delta_1$ much smaller than $\varepsilon_1$,
with the effect of stopping because of $\delta_1$ all the time and never because 
of $\varepsilon_1$; this would be allowed by our argument, but it would be a bad and 
confusing practice.

There is a second issue with $\delta_1$, which is that allowing $\delta_1$ to be small
(as we will do to simplify the proof) should have an advantage, which is that our bilipschitz 
mappings are actually bilipschitz with constants that are as close to $1$ as we want. 
We claim that this is true, but it is less easy to use because the estimates in \cite{DT} 
that prove this are rather well hidden, so we decided that we shall not use this extra information 
(other than saying that we have a uniform bound on the bilipschitz constants)
and merely add remarks along the proof that explain how
we could get and deal with this additional information.
\end{remark}

\section{The approximating surface $\Sigma$}
\label{S2}

We shall now describe the main lines of the construction of \cite{DT}, where one starts from a  stopping time region $\cF$ like the one above, and constructs an associated Reifenberg flat set 
$\Sigma$, parameterized by a mapping 
$f : \Sigma_0 =P_0=\RR^d \to \Sigma$, and even a global change of variable 
$g : \R^n \to \R^n$ (with $g_{\vert \Sigma_0} = f$).

For the construction to work, one needs to find what is called a 
\underline{coherent collection} \underline{of balls and planes} (in short, a CCBP),
which will be our first task here. This will involve choosing some collections
of $d$-planes, and let us first see what we have.

Recall that for each $Q \in \cF$, we are given a $d$-plane $P(Q)$ 
that satisfies \eqref{a9}-\eqref{a12}. In particular, \eqref{a9} says that
$d_{x_Q, M l(Q)}(E,P(Q)) \leq \varepsilon_1$.
This means that
\begin{equation} \label{b1}
\dist(y,P(Q)) \leq \varepsilon_1 M l(Q) 
\ \text{ for } y\in E \cap B(x_Q, M l(Q)),
\end{equation}
(and in particular $P(Q)$ passes within $2 \varepsilon_1 M l(Q)$ of $x_Q$), but also
\begin{equation} \label{b2}
\dist(y,E) \leq \varepsilon_1 M l(Q) 
\ \text{ for } y\in P(Q) \cap B(x_Q, M l(Q)).
\end{equation}
We do this also for $Q_0$ 
(which we have assumed to lie in $\cF$),
and call
$P_0 = P(Q_0)$ the plane that we get. We shall even assume, without loss
of generality, that 
\begin{equation} \label{b3}
P_0 = \R^d \ 
\text{ and } x_{Q_0} = 0
\end{equation}

In \cite{DT}, which we shall often refer to as ``there'', a CCBP starts with the 
choice of families $\{ B_{j,k} \}$, $j \in J_k$ of balls, where $k \geq 0$ still denotes
a generation. In fact $B_{j,k} = B(x_{j,k},r_k)$, where $r_k = 10^{-k}$ as above, so we
just need to choose the centers $\{ x_{j,k} \}$, $j \in J_k$. 
Recall the definition of $\cF(k)$ in \eqref{b4}, let
\begin{equation} \label{b5}
E(k) = \big\{ x\in E \, ; \, \dist(x,Q) \leq \frac{M r_k}{10}  \text{ for some } Q \in \cF(k) \big\},
\end{equation}
and finally pick a maximal family $\{ x_{j,k} \}$, $j \in J_k$, of points of $E(k)$
that lie at distances at least $r_k$ from each other. This defines our family of balls.
We need to check a coherence condition, (2.3) in \cite{DT}, that demands that for 
$k \geq 1$, each $x_{j,k}$ lies in $B(x_{i,k-1},2r_{k-1})$ for some $i\in J_{k-1}$. 
This comes from the heredity condition for $\cF$: since $x_{j,k} \in E(k)$, 
we know that $\dist(x_{j,k},Q) \leq \frac{M r_k}{10}$ for some $Q \in \cF(k)$; 
the parent $Q'$ of $Q$ lies in $\cF(k-1)$, and 
since $\dist(x_{j,k},Q') \leq \dist(x_{j,k},Q) \leq \frac{M r_k}{10}$, 
$x_{j,k} \in E(k-1)$ and we can find a point $x_{i,k-1}$, $i\in J_{k-1}$, 
that lies within $r_{k-1}$ of $x_{j,k}$.

We should also choose a nice surface $\Sigma_0$ with which we start the construction;
here we simply take $\Sigma_0 = P_0 = P(Q_0)$, and the properties (2.4)-(2.7) required
in \cite{DT} are easily satisfied; in particular (2.7) there follows from \eqref{b1} 
if $\varepsilon_1$ is chosen small enough, depending on $\varepsilon$ there.

Finally we need to associate a $d$-plane $P_{j,k}$ to each ball $B_{j,k}$, 
and this is easy to do: for each $j\in J_k$ we choose $Q_{j,k} \in \cF(k)$ such that
$\dist(x_{j,k},Q_{j,k}) \leq \frac{M r_k}{10}$, and then we set $P_{j,k} = P(Q_{j,k})$.
Notice that when $k=0$, we have many points $x_{j,0}$ (because $E(0)$ is rather large),
but all of them are associated to $P_0$.
 
There is an unfortunate little catch here, because it is also required in \cite{DT} 
that $P_{j,k}$ goes through $x_{j,k}$, but we really like $P(Q_{j,k})$ here,
in fact more than the precise location of $x_{j,k}$. So we modify the construction
a little bit. We start with a maximal collection of points $\wt x_{j,k} \in E(k)$, 
at mutual distances at least $\frac{11r_k}{10}$, define the $Q_{j,k}$ 
and $P_{j,k}$ as above, and then use \eqref{b1} (with $\varepsilon_1$ small enough) 
to find $x_{j,k} \in P(Q_{j,k}) \cap B(\wt x_{j,k}, \frac{r_k}{100})$, and use these in the
definition of $B_{j,k}$. This does not perturb our proof of (2.3) there, we lose
the fact that $x_{j,k} \in E$, which looked nicer, but this is not needed to apply
Theorems 2.4 and 2.5 there. Starting from (13.3) in \cite{DT}, another trick is explained,
which allows us to replace $\wt x_{j,k}$ with another point $x_{j,k} \in E$ that lies 
so close to $P_{j,k} = P(Q(j,k))$ that we could translate $P_{j,k}$ slightly and keep
\eqref{a10} with almost the same constant, but we don't need to do this.

We need to check that the $P_{j,k}$ satisfy compatibility conditions 
(namely, (2.7)-(2.10) there). We start with (2.8), which  demands that for 
$k \geq 0$ and all $i, j \in J_k$ such that $|x_{i,k}-x_{j,k}| \leq 100 r_k$, 
$P_{j,k}$ and $P_{i,k}$ are so close that
\begin{equation} \label{b6}
d_{x_{j,k},100r_k}(P_{i,k},P_{j,k}) \leq \varepsilon,
\end{equation}
for some $\varepsilon > 0$ that needs to be small enough for the construction
of \cite{DT} to work. We choose $\varepsilon_1$ small enough, depending on $\varepsilon$
and our constant $M$, and then this follows from our definitions, and in particular
\eqref{b1} and \eqref{b2}; the verification is fairly simple, and is essentially done
in \cite{DT}, below Lemma 12.2 on page 66, so we skip it.

In our case, (2.9) there is just a special case of (2.8) because $\Sigma_0 = P_0$, 
and (2.10) demands that for $k \geq 0$, $i\in J_k$, and $j \in J_{k+1}$ such that 
$|x_{i,k}-x_{j,k+1}| \leq 2 r_k$,
\begin{equation} \label{b7}
d_{x_{i,k},20r_k}(P_{i,k},P_{j,k+1}) \leq \varepsilon.
\end{equation}
The verification is almost the same as for \eqref{b6}, and we also
refer to the argument in \cite{DT}, below Lemma 12.2.

\ms
At this stage we are able to apply Theorem 2.4 in \cite{DT}, which provides us
with a $C\varepsilon$-Reifenberg flat set $\Sigma$ and biH\"older mappings 
$f : P_0 \to \Sigma$ and $g: \R^n \to \R^n$, with some good properties. 

\begin{theorem}[Theorem 2.4 from \cite{DT}]\label{t24DT} 
Let $(\Sigma_0, \{B_{j,k}\}, \{P_{j,k}\})$ be a CCBP as above and assume that $\eps$ is small enough depending on $n$ and $d$. Then there is a bijection $g:\RR^n\to \RR^n$ with the following properties; 
$$g(z)=z \quad\mbox{when}\quad \dist(z, \Sigma_0)\geq 2,  $$
$$g(z)-z \leq C\eps \quad\mbox{for}\quad z\in \RR^n, $$
$$\frac14 \, |z'-z|^{1+C\eps}\leq |g(z)-g(z')|\leq 3\,  |z'-z|^{1-C\eps} $$
for $z, z' \in \RR^n$ such that $|z-z'|\leq 1$, and $\Sigma=g(\Sigma_0)$ is a $C\eps$-Reifenberg flat set that contains the accumulation set defined as the 
collection of all $x\in \RR^n$ which can be written as $x=\lim_{m\to\infty} x_{j(m), k(m)}$ with $k(m)\in \NN$ such that $\lim_{m\to\infty} k(m)=\infty$ and  $j(m)\in J_{k(m)}$ for $m\geq 0$. The constant $C$ depends on $n$ and $d$ only.\end{theorem}

However, we are interested in more precise properties of $g$ (such as the fact that it
is bilipschitz), and we will also need some information that comes from the construction,
because in \cite{DT} no special attention was given to the specific form of the Jacobian
matrix of $g$, which we need to study for our application to degenerate elliptic
operators. We start with the bilipschitz part.

\begin{lem}\label{Lb1}
If $\varepsilon_1$ is small enough, the mapping $g : \R^n \to \R^n$ is bilipschitz, 
with a constant that depends only on $\delta_1$, $n$, $d$, and the different choices 
above (that depend on $M_0$, for instance).
\end{lem}

This will follow from Theorem 2.5 there and our additional constraint \eqref{a12}, 
once we decipher some additional definitions. But before we do this, let us check that
whenever $Q \in \cF(k)$ and $R \in \cF(k) \cup \cF(k-1)$ are such that
$\dist(Q,R) \leq \frac{M l(Q)}{2}$, then
\begin{equation} \label{b8}
d_{x_Q,Ml(Q)}(P(Q),P(R)) \leq C \beta(Q) + C \beta(R),
\end{equation}
where $C$ depends on $M$, $C_0$, $n$, and $d$, but this will not matter.

Let us rapidly prove \eqref{b8}. 
The argument is similar to what was done in \cite{DS}, below (13.25). 
We intend to use the fact that both $P(Q)$ and $P(R)$ are very close to $E$ in a common 
region to compare their positions. First choose an orthonormal basis 
$e_1, \ldots, e_d$ of the vector $d$-plane parallel to $P(R)$, and consider the 
points $\xi_0 = x_R$ and, for $1 \leq i \leq d$, $\xi_i = x_R + r_k e_i$. 
Notice that $B(x_R, 2r_k)$ lies well inside $B(x_R, M l(R))$, so by  
\eqref{b2} we can find points $x_i \in E$, $0 \leq i \leq d$, such that
$|x_i-\xi_i| \leq \varepsilon_1 M l(R) \leq 10^{-2} r_k$. Then we use \eqref{a10},
the Ahlfors regularity of $\mu$, and Chebyshev's inequality
to find that for more than half of the points $z \in E \cap B(x_i,10^{-2} r_k)$,
$\dist(z,P(R)) \leq C l(R) \beta(R)$. But also, $B(x_R, 2r_k)$ lies well inside $B(x_Q, M l(R))$,
so we can also apply \eqref{a10} to $P(Q)$, and find that for a majority of points
$z \in E \cap B(x_i,10^{-2} r_k)$, $\dist(z,P(Q)) \leq C l(Q) \beta(Q)$. For each $i$ we 
select a point $y_i$ with both properties, and this gives $z_i\in P(R) \cap B(\xi_i,10^{-1}r_k)$
and $w_i\in P(R) \cap B(\xi_i,10^{-1})$ such that $|z_i-w_i| \leq C r_k (\beta(Q)+\beta(R))$.
At this point, we have sufficiently many points of contact between $P(Q)$ and $P(R)$
to control their relative positions and prove \eqref{b12}; see also Lemma 12.7 there.

Now return to the lemma and Theorem 2.5 there.
Define the numbers $\varepsilon''_k(y)$, $k \geq 1$ and $y\in \R^n$, by 
\begin{equation} \label{b9}
\begin{split}
\varepsilon''_k(y) = \sup\big\{ &
d_{x_{i,l},100r_{l}}(P_{j,k},P_{i,l}) \, ; \,
j\in J_k, \, l \in \{ k-1, k \}, 
\\& \hskip 3.5cm
i\in J_{l}, \hbox{ and } y \in 11 B_{j,k} \cap 12 B_{i,l} 
\big\}  
\end{split}
\end{equation}
when $y \in  V_k^{11} =  \bigcup_{j\in J_k} B(x_{j,k}, 11 r_k)$, and simply by
$\varepsilon''_k(y) = 0$ when $y \notin V_k^{11}$. This is the same definition as 
in (2.17) and (2.18) there, and then Theorem 2.5 there says that $g$ is bilipschitz
as soon as there is a constant $M_3 \geq 0$ such that
\begin{equation} \label{b10}
\sum_{k \geq 0} \varepsilon''_k(g(z))^2 \leq M_3  \ \hbox{ for all } z \in P_0.
\end{equation}

Thus, in order to deduce the lemma from that result, we will just need to show that
the numbers $\beta(Q)$ of \eqref{a10} control the $\varepsilon''_k(y)$, $y \in V_k^{11}$. 

So let $z\in P_0$ be given, set $y = g(z)$, and let $k \geq 1$ be such that
$y \in V_k^{11}$ (we don't care about the other $k$, since $\varepsilon''_k(y)=0$). 
This last means that $y \in B(x_{j_0,k}, 11 r_k)$ for some $j_0 \in J_k$, but this 
will not really matter.

Next let $j$, $l$, and $i$ be as in \eqref{b9}, and follow the definitions: we picked
a cube $Q = Q_{j,k} \in \cF(k)$ such that $\dist(x_{j,k},Q) \leq \frac{M r_k}{10}$, 
and then we set $P_{j,k} = P(Q)$, and similarly we chose 
$R = Q_{i,l} \in \cF(l)$ such that $\dist(x_{i,l},R) \leq \frac{M r_l}{10}$
and then set $P_{i,l} = P(R)$. We record for later the fact that 
\begin{equation} \label{b11}
\dist(y,Q) \leq \frac{M r_k}{5} \ \text{ and } \ \dist(y,R) \leq \frac{M r_k}{5}.
\end{equation}
Obviously $\dist(Q,R) \leq \frac{M l(Q)}{2}$,
so \eqref{b8} says that 
\begin{equation} \label{b12}
d_{x_{i,l},100r_l}(P(Q_{j,k}),P(Q_{i,l})) = d_{x_{i,l},100r_l}(P(Q),P(R))
\leq C (\beta(Q)+ \beta(P)).
\end{equation}
For each scale $k \geq 0$, denote by $\cQ(y,k)$ the collection of cubes
$Q\in \cF(k)$ such that $\dist(y,Q) \leq \frac{M r_k}{5}$. Obviously
\begin{equation} \label{b13}
\sum_{k \geq 0} \varepsilon''_k(g(z))^2 
\leq C \sum_{k \geq 0} \sum_{Q \in \cQ(y,k)} \beta(Q)^2
\end{equation}
by \eqref{b11} and \eqref{b12}, and we shall use \eqref{a12} to control the right-hand side.
If the Jones function $J$ were only using 
the $\beta(Q)$, this may seem complicated;
here we can proceed as follow. Let $k_0$ be such that $\cQ(y,k_0)$ is not empty, and select
$Q \in \cQ(y,k_0)$; then for $0 \leq k \leq k_0$, denote by $Q_k$ the ancestor of $Q$
that lies in $\bD_k$; observe that $Q_k \in \cF(k)$ by heredity, and 
$Q_k \in \cQ(y,k)$ because $\dist(y,Q_k) \leq \dist(y,Q) \leq \frac{M r_k}{5}$.
Now all the other cubes $S$ of $\cQ(y,k)$ lie at distance less than $M r_k = M l(Q_k)$
from $Q_k$, so $\beta(S) \leq \wt \beta(Q_k)$ by \eqref{a11bis}. Thus
\begin{equation} \label{b14}
\sum_{0 \leq k \leq k_0} \sum_{Q \in \cQ(y,k)} \beta(Q)^2
\leq C \sum_{0 \leq k \leq k_0} \wt\beta(Q_k)^2 \leq C J(Q) \leq C \delta_1
\end{equation}
because each $\cQ(y,k)$ has at most $C$ elements, and by \eqref{a11} and \eqref{a12}.
Since this is true for every $k_0$ (with the same constant), we get \eqref{b10} and, as promised,
Lemma~\ref{Lb1} follows from Theorem 2.5 there.
\qed

\begin{rem} \label{r416}
For this lemma we do not need $\delta_1$ to be small, but the first author
claims that taking $\delta_1$ small would allow us to get a bound for the bilipschitz constant 
for $g$ which is as close to $1$ as we want. This would be reassuring, but apparently 
the authors of \cite{DT} were too busy controlling the large constants to make a clear remark, anywhere in that paper, to the effect that small bounds for $J$ yield small bilipschitz bounds
for $g$. We will 
manage not to
use this remark in this paper, so as not to make the reader feel too bad,
but will add some comments to this effect for the case when they  would be badly needed in the future.
\end{rem}

Lemma \ref{Lb1} will be quite useful to help us control other terms; for instance,
we will not need to worry about supremum norms for the derivatives of our mappings.
But we will need more information, typically on the structure of $Df$ and $Dg$,
so let us 
step back and recall the construction of $f$ and auxiliary functions $f_k$ and then 
we will pass to the construction of $g$ in Section~\ref{Sd}.

In Section 3 of \cite{DT}, one constructs a \ub{partition of unity} for each generation 
$k \geq 0$, composed of functions $\theta_{j,k}$, $j\in J_k$, supported in $10B_{j,k}$
plus a function $\psi_k$ supported away from $V_k^8 = \cup_{j\in J_k} 8B_{j,k}$. Thus
\begin{equation} \label{b15}
\psi_k + \sum_{j\in J_k} \theta_{j,k} = 1
\end{equation}
as in (3.13) there, and 
\begin{equation} \label{b16}
\sum_{j\in J_k} \theta_{j,k} = 1 \ \text{ on } V_k^8 = \cup_{j\in J_k} 8B_{j,k}.
\end{equation}
In addition, $|\nabla^m \theta_{j,k}| \leq C_m r_k^{-m}$, as expected (see (3.15) there).

\ms
Then we can \ub{define the mapping} $f$ on $\Sigma_0 = P_0 = \R^d$, as the limit of
functions $f_k$ defined by induction by 
\begin{equation} \label{b17}
f_0(y)=y \ \text{ and } \ f_{k+1} = \sigma_k \circ f_k,
\end{equation}
(as in (4.1) there), where $\sigma_k$ is a map that tries to move points in the direction of
$E$ (or rather, the local $P_{j,k}$), and is defined by
\begin{equation} \label{b18}
\sigma_k(y) = y + \sum_{j \in J_k} \theta_{j,k}(y) \, [\pi_{j,k}(y)-y]
=\psi_k(y) y + \sum_{j \in J_k}\theta_{j,k}(y) \, \pi_{j,k}(y)
\end{equation}
(as in (4.2) there), where $\pi_{j,k}$ denotes the orthogonal projection
from $\R^n$ onto $P_{j,k}$ and the equality comes from \eqref{b15}. 
It turns out that the $f_k$ converge quite fast to a limit mapping $f$,
which is our parameterization of the nice 
Reifenberg-flat surface $\Sigma = f(\Sigma_0)$. 
We remind the reader that $g$ will ultimately be defined so that $g=f$ on $\Sigma_0$ (see (10.13) in \cite{DT} and Section~\ref{Sd}) but the construction in \cite{DT} starts with the $f_k$ and $f$. 

Some observations will be useful concerning the local regularity of the intermediate surfaces
\begin{equation} \label{b19}
\Sigma_k = f_k(\Sigma_0),
\end{equation}
and the way each one maps to the next one. 
Proposition 5.1 in \cite{DT} gives a \ub{good local description} of $\Sigma_k$
in terms of Lipschitz graphs, which we can summarize as follows. For each
$j \in J_k$, there is a $C\varepsilon$-Lipschitz function 
$A_{j,k} : P_{j,k} \to P_{j,k}^\perp$, with $|A_{j,k}(x_{j,k})| \leq C \varepsilon r_k$,
such that inside $49B_{j,k}$, $\Sigma_k$ coincides with the graph
$\Gamma_{j,k}$ of $A_{j,k}$ over $P_{j,k}$. The same proposition also says
that $A$ is of class $C^2$, but does not record estimates on this, and this is a
part that we will need to complement.

We shall not use Proposition 5.1 there directly so much, but it is important in the description 
of trajectories that follows, and contains the estimate (5.11) there, 
which says that 
\begin{equation} \label{b35}
|\sigma_k(y)-y| \leq C \varepsilon r_k \ \text{ for $k \geq 0$ and } y\in \Sigma_k,
\end{equation}
which, after using \eqref{b17} repeatedly and summing a geometric series, yields
\begin{equation} \label{b36}
|f(x)-f_k(x)| \leq C \varepsilon r_k \ \text{ for $x\in \Sigma_0$ and $k \geq 0$.} 
\end{equation}

For the \ub{description of trajectories} that follows, 
we continue to
use the notation
\begin{equation} \label{b20}
V_k^{A} = \bigcup_{j \in J_k} B(x_{j,k},Ar_k) = \bigcup_{j \in J_k} AB_{j,k},
\end{equation}
when $A$ is an integer. When $y\in \Sigma_k \cap V_k^{8}$ 
(we call this \ub{the active region}), then $\psi_k(y) = 0$, and 
the formula \eqref{b18} becomes the simpler
\begin{equation} \label{bb23}
\sigma_k(y) =  \sum_{j \in J_k}\theta_{j,k}(y) \, \pi_{j,k}(y),
\end{equation}
where we know that in addition
\begin{equation} \label{bb24}
\sum_{j \in J_k}\theta_{j,k}(y) = 1.
\end{equation}
When on the opposite $y \in \R^n \sm V_k^{10}$ (we call this the dead region),
things are simple too, because all the $\theta_{j,k}(y)$ vanish, hence $\psi_k(y)=1$,
and \eqref{b18} says that
\begin{equation} \label{bb25}
\sigma_k(y) = y \ \text{ for } y\in \R^n \sm V_k^{10}.
\end{equation}
Things are a little more unpleasant in $V_k^{10} \sm V_k^{8}$, but fortunately 
the next lemma says that this 
never happens more than once along a given trajectory, 
and this will leave a reasonably small trace in our Carleson measure estimates.
We call $V_k^{10} \sm V_k^{8}$ the dying region.

\begin{lem}\label{Lb2}
Let $x\in  \Sigma_0$ be given, and denote by $y_k = f_k(x) \in \Sigma_k$ its successive
images. 
\begin{equation} \label{bb26}
\text{If $y_k \in \R^n \sm V_k^{10}$ for some $k\geq 0$, then 
$y_l=y_k \in \R^n \sm V_l^{10}$ for $l \geq k$;}
\end{equation}
\begin{equation} \label{bb27}
\text{If $y_k \in V_k^{10}$ for some $k\geq 1$, then 
$y_l \in V_l^{4}$ for $0 \leq l \leq k-1$.}
\end{equation}
\end{lem}

\ms 
This is Lemma 6.1 in \cite{DS}. Notice that if $y_k \in V_k^{10} \sm V_k^{8}$, 
\eqref{bb27} says that the previous images were in the active region, and also
(applying it to $y_{k+1}$) that $y_{k+1} \in \Sigma_{k+1}\sm V_{k+1}^{10}$
lies in the dead region, as well as all its successors (by \eqref{bb26}).

\ms
We need some \ub{estimates on $\sigma_k$} and its derivative that were not
necessarily recorded there. We claim that
\begin{equation} \label{bb28}
|f_{k+1}(x)-f_k(x)| \leq C \varepsilon''_{k}(f(x)) \, r_k
\ \text{ for $x\in \Sigma_0$ such that } 
f_k(x) \in V_k^8.
\end{equation}
Set $y = f_k(x) \in V_k^8$, and choose $j\in J(k)$ such that $|y-x_{j,k}| \leq 8r_k$. 
Then (7.8) there says that 
\begin{equation} \label{bb29}
|\sigma_k(y)-\pi_{j,k}(y)| \leq C \varepsilon_{k}(y) r_k,
\end{equation}
where $\pi_{j,k}$ denotes the orthogonal projection onto $P_{j,k}$ and the 
function $\varepsilon_k$ is defined in (7.7) there (we shall return to this soon). Thus 
\begin{equation} \label{bb30}
|f_{k+1}(x)-f_k(x)| = |\sigma_k(y)-y| \leq \dist(y, P_{j,k}) + C \varepsilon_{k}(y) r_k,
\end{equation}
and we now evaluate that distance. When $k=0$, $y = f_0(x) = x \in P_0$,
we actually took $P_{j,k} = P_0$, and the distance is $0$.
Otherwise, set $y' = f_{k-1}(x)$, observe that $y'\in V_{k-1}^4$ by Lemma \ref{Lb2},
and choose $i\in J(k-1)$ such that $|y'-x_{i,k-1}| \leq 4r_{k-1}$. This time (7.8) there 
says that 
\begin{equation} \label{bb31}
|\sigma_{k-1}(y')-\pi_{i,k-1}(y')| \leq C \varepsilon_{k-1}(y') r_{k-1},
\end{equation}
where $\pi_{i,k-1}$ denotes the orthogonal projection onto $P_{i,k-1}$.
Since $y = \sigma_{k-1}(y')$, we get that
\begin{equation} \label{bb32}
\dist(y, P_{i,k-1}) \leq C \varepsilon_{k-1}(y') r_{k-1}.
\end{equation}
Recall from \eqref{b36} that $|f(x)-y| \leq C\varepsilon r_k$, and similarly
$|f(x)-y'| \leq 10 C\varepsilon r_k$. Notice then that
$f(x) \in 11B_{j,k} \cap 12B_{i,k-1}$, so the definition \eqref{b9} says that
$P_{j,k}$ and $P_{i,k-1}$ are $100r_{k-1} \varepsilon''_k(f(x))$-close to each other
in $B(x_{i,k-1},100r_{k-1})$. In particular \eqref{bb32} implies that 
\begin{equation} \label{bb33}
\dist(y, P_{j,k}) \leq C \varepsilon_{k-1}(y') r_k + C \varepsilon''_k(f(x)) \, r_k.
\end{equation}
Now we compare the definition (7.7) there of $\varepsilon_{k}(y)$ with \eqref{b9}
and find out that $\varepsilon_{k}(y) \leq C \varepsilon''_k(f(x))$, because
if $y\in 10B_{i,k} \cap 10B_{j,k}$ for some $i, j\in J(k)$, then
$f(x) \in 11B_{i,k} \cap 11B_{j,k}$. Similarly, $\varepsilon_{k-1}(y') \leq C \varepsilon''_k(f(x))$,
with the only small difference that since this time we are comparing two planes of generation $k-1$, 
we need to go through our chosen plane $P_{j,k}$ of generation $k$.
Now our claim \eqref{bb28} follows from \eqref{bb30} and \eqref{bb33}.

\ms
We also need \ub{estimates on the derivatives} of $f_k$, and of course we first 
differentiate $\sigma_k$. We start in the active region (the open set $V_k^8$),
where we can use the simpler formulas \eqref{bb23} and \eqref{bb24}, and hence
\begin{equation} \label{b22}
D\sigma_k(y) =  \sum_{j \in J_k}\theta_{j,k}(y) \, D\pi_{j,k} +
\sum_{j \in J_k} D\theta_{j,k}(y) \pi_{j,k}(y),
\end{equation}
where the differential $D\pi_{j,k}$ of $\pi_{j,k}$ is the vector projection (which does not
depend on $y$), which we try not to mix with the affine projection $\pi_{j,k}$. 
In this sort of situation, we like to pull out a specific index $j(y) = j_k(y)$ such that 
$\theta_{j(y),k}(y) \neq 0$, and use the fact that $\sum D\theta_{j,k}(y) = 0$
by \eqref{bb24} to write that
\begin{equation} \label{b23}
D\sigma_k(y) =  \sum_{j \in J_k}\theta_{j,k}(y) \, D\pi_{j,k} +
\sum_{j \in J_k} D\theta_{j,k}(y) \, [\pi_{j,k}(y)-\pi_{j(y),k}(y)],
\end{equation}
and even
\begin{equation} \label{b24}
D\sigma_k(y) - D\pi_{j(y),k} = \sum_{j \in J_k}\theta_{j,k}(y) \, [D\pi_{j,k} - D\pi_{j(y),k}]
+\sum_{j \in J_k} D\theta_{j,k}(y) \, [\pi_{j,k}(y)-\pi_{j(y),k}(y)].
\end{equation}
We differentiate once more (but keep the same index $j(y)$ to do the computations near
$y$; we certainly don't want to differentiate $j(y)$) and get that
\begin{equation} \label{b25}
D^2\sigma_k(y) =  2\sum_{j \in J_k} D\theta_{j,k}(y) \, [D\pi_{j,k}-D\pi_{j(y),k}] +
\sum_{j \in J_k} D^2\theta_{j,k}(y) \, [\pi_{j,k}(y)-\pi_{j(y),k}(y)].
\end{equation}
Let us not pay too much attention on what we mean by multiplication in these formulas; 
the main thing is the size estimate that follows. 
In this sum the only terms that do not vanish come from balls such that $y \in 10 B_{j,k}$,
and there are at most $C$ of them. The size of $D \theta_{j,k}$ and $D^2 \theta_{j,k}$
is controlled below \eqref{b16}.
We look at the definition \eqref{b9} of $\varepsilon''_k$, and find out that 
for $y \in V_k^8$
\begin{equation} \label{b26}
|D^2\sigma_k(y)| \leq C \varepsilon''_k(z) \, r_k^{-1}
\ \text{ for any } z\in B(y,r_k),
\end{equation}
which we take as a good estimate. Here we shall just take $z = f(x)$ for the
point $x\in \R^d$ such that $y = f_k(x)$, and the fact that $z\in B(y,r_k)$
comes from \eqref{b36}. This was our better estimate for $y\in V_k^8$.

\ms
In the dead region $\R^n \sm V_k^{10}$ where all the $\theta_{j,k}$ vanish,
we have $\psi_k = 1$, $\sigma_k(y) = y$ (by \eqref{b18}), and hence $D\sigma_k = I$
and $D^2 \sigma_k = 0$ (see also (4.5) there). 

In the dying region $V_k^{10} \sm V_k^8$, we don't have 
very good estimates because $\psi_k$ is not identically $1$ near $y$.
This time we start from the first part of \eqref{b18}, which yields
\begin{equation} \label{b27}
D\sigma_k(y)- I = \sum_{j \in J_k} D\theta_{j,k}(y) \, [\pi_{j,k}(y)-y]
+ \sum_{j \in J_k} \theta_{j,k}(y) \, [D\pi_{j,k}-I]
\end{equation}
and then 
\begin{equation} \label{b28}
D^2\sigma_k(y) = \sum_{j \in J_k} D^2\theta_{j,k}(y) \, [\pi_{j,k}(y)-y]
+ 2\sum_{j \in J_k} D\theta_{j,k}(y) \, [D\pi_{j,k}-I];
\end{equation}
we observe that $|\pi_{j,k}(y)-y| \leq 10r_k$ 
when $D\theta_{j,k}(y) \neq 0$ or $D^2\theta_{j,k}(y)\neq 0$
because $P_{j,k}$ goes through $x_{j,k}$ and $\theta_{j,k}$ is supported in $10B_{j,k}$; 
this yields the brutal estimate
\begin{equation} \label{b29}
|D\sigma_k(y)| \leq C \ \text{ and } \ |D^2\sigma_k(y)| \leq C 10^{k}
\ \text{ for } y \in V_k^{10} \sm V_k^8. 
\end{equation}

\ms
Next we use this to estimate $Df_k$ and $D^2 f_k$. Recall from \eqref{b17}
that $f_{k+1} = \sigma_k \circ f_k$; thus (with probably very bad but yet understandable)
notation,
\begin{equation} \label{b30}
D f_{k+1}(x) = D\sigma_k(f_k(x)) \circ Df_k(x)
\end{equation}
(but we shall not always write the variables) and then 
\begin{equation} \label{b31}
D^2 f_{k+1}(x) = D^2\sigma_k(f_k(x)) [ Df_k(x), Df_k(x)] + D\sigma_k(f_k(x)) [D^2f_k(x)]
\end{equation}
with ugly notation, but we immediately put norms everywhere, forget the algebra, and get that
\begin{equation} \label{b32}
|D^2 f_{k+1}(x)| \leq 
C |D^2\sigma_k(f_k(x))| + 2 |D^2f_k(x)|
\end{equation}
also because we know that all the $f_k$ are bilipschitz with uniform constants 
(that may depend on $\delta_1$). 
We may rewrite this as 
\begin{equation} \label{b33}
r_{k+1} |D^2 f_{k+1}(x)| \leq 
C r_{k} |D^2\sigma_k(f_k(x))| + \frac15  |r_k D^2f_k(x)|
\end{equation}
because this is the proper scaling, and this way we insist on the fact that the second term 
contributes less. 

We start in the most interesting case when $y = f_k(x)$ lies in the 
active region $V_k^8$; then we use \eqref{b26} with $z = f(x)$ and get that
\begin{equation} \label{b34}
r_{k+1} |D^2 f_{k+1}(x)| 
\leq C \varepsilon''_k(f(x))  + \frac15  |r_k D^2f_k(x)|
\end{equation}
where $z$ is any point of $B(y,r_k)$. It should be noted that when $y = f_k(x)$ 
lies in the active region $V_k^8$, Lemma \ref{Lb2} says that this was the case
for all the previous images $f_l(x)$, $l < k$, so we also have estimates like \eqref{b34}
for these, that we can compose. We get that
\begin{equation} \label{b37}
|r_{k+1} D^2f_{k+1}(x)| \leq C \wt\varepsilon_k(x),
\ \text{ where } 
 \wt\varepsilon_k(x) = \sum_{l \leq k} 5^{l-k} \varepsilon''_l(f(x)).
\end{equation}

This was when $y\in  V_k^8$. 
When $y$ lies in the dying region $V_k^{10} \sm V_k^8$,
we use the bad estimate \eqref{b29} for $k$, but observe that the previous $f_l(x)$, 
$l < k$, were in the active region (by Lemma~\ref{Lb2}), so we can use the 
estimate \eqref{b37} for $D^2f_k$. Thus \eqref{b33} yields
\begin{equation} \label{b38}
|r_{k+1} D^2f_{k+1}(x)| \leq C + C \wt\varepsilon_k(x) \leq C,
\end{equation}
which is not a great estimate but should be enough. In the remaining case when
$y \in \Sigma_k \sm V_k^{10}$, we denote by $l < k$ the last time when $f_l(x)$
was in the active or dying region, and use \eqref{b37} or \eqref{b38} to prove that
\begin{equation} \label{b39}
|D^2f_{k+1}(x)| = |D^2f_{l+1}(x)| \leq C r_l^{-1}.
\end{equation}
These estimates on the second derivatives will be enough for the better control
that we want on the Jacobian matrix of our global mapping $g$.

\begin{rem}\label{r453}
Yet we feel bad about 
using such rude estimates, 
so let us rapidly say why \eqref{b29}, and then
\eqref{b38} can be improved. Our estimate $|\pi_{j,k}(y)-y| \leq 10r_k$ 
below \eqref{b28} was really lousy, because in fact $|\pi_{j,k}(y)-y| \leq C \varepsilon r_k$ 
when $y \in \Sigma_k$ (we have a good Lipschitz graph description of
$\Sigma_k$ near $y$). We need to be more careful about the terms
with $D\pi_{j,k}-I$, because of course it is not small. Yet, when we apply it
to a tangent vector $v$ to $\Sigma_k$, its effect is indeed
of size at most $C\varepsilon |v|$, because $T\Sigma_k(y)$ is nearly parallel to
$P_{j,k}$. This is good, because when we compose with $f_k$, we only
compute $D\sigma_k$ on vectors $v$ parallel to $T\Sigma_k(y)$.
So we can add a factor of size $\varepsilon$ in \eqref{b29}, and then
\eqref{b39}. This way, we get the not too surprising result that $f$
(and then $g$ later) is biLipschitz with a constant which can be taken as 
close to $1$ as we want, provided that we take $\delta_1$ small enough.
\end{rem}

\section{Tangent planes and fields of rotations}
\label{S3}

Once the mappings $f_k$ and the surfaces $\Sigma_k = f_k(P_0)$
are under control, \cite{DT} starts the construction of the mapping $g$.
The general idea is that for $(x,t) \in \R^d \times \R^{n-d}$, 
$g(x,t)$ should be obtained from $f(x)$ by going in the orthogonal 
direction, and at distance roughly $|t|$. Of course we need to organize this
in a coherent way, and also it actually makes more sense to start from $f_k(x)$,
and go in a direction orthogonal to the tangent direction of $\Sigma_k$,
because $\Sigma_k$ is smoother. This makes a difference because the limit
object $\Sigma$ may be spiraling at small scales.

So our first task will be to study a little the variations of the tangent plane to 
$\Sigma_k$. Here we roughly follow Chapter 9 of \cite{DT}.
We know, for instance from the first lines of Chapter 7 there, or more directly
the description in Proposition~5.1 there, that
each $\Sigma_k$ is (at least) of class $C^2$. Let us denote by
$T\Sigma_k(y)$ the tangent plane to $\Sigma_k$ at $y\in \Sigma_k$.
We also set
\begin{equation} \label{c1}
T_k(x) = T\Sigma_k(f_k(x)),
\end{equation}
and denote by $\pi_k(x)$ the orthogonal projection on the vector $d$-space
parallel to $T_k(x)$; it is easier to define the differential of $\pi_k$ than $T_k$, 
which is why we will often consider $\pi_k$. 
We claim that when $f_k(x) \in V_k^8$ (the active region),
\begin{equation} \label{c2}
|D \pi_{k+1}(x)| \leq C r_k^{-1} \wt\varepsilon_k(x).
\end{equation}
There is no problem with the existence of $D \pi_{k+1}(x)$, because $\Sigma_k$ is $C^2$,
so it is enough to show that for $x'\in P_0$, close enough to $x$,
\begin{equation} \label{c3}
|\pi_{k+1}(x)-\pi_{k+1}(x')| \leq C r_k^{-1} \wt\varepsilon_k(x) |x'-x|.
\end{equation}
Let us evaluate the distance between $T'_{k+1}(x)$, the vector space parallel to $T_{k+1}(x)$,
and its analogue $T'_{k+1}(x')$ for $T_{k+1}(x')$. Let $w\in T'_{k+1}(x)$ be given; we can write 
$w= Df_{k+1}(x)\cdot v$ for some vector $v\in \R^d$, and since we know 
(from Lemma \ref{Lb1} and the proof of Theorem 2.5 there) that the $f_k$ are bilipschitz 
with uniform bounds, we also get that $|v|\leq C|w|$. 
We know from \eqref{b37} that for $x'$ close to $x$,
\begin{equation} \label{c4}
|Df_{k+1}(x)\cdot v - Df_{k+1}(x')\cdot v| 
\leq C \wt\varepsilon_k(x) \, r_{k+1}^{-1} |x-x'| |v|;
\end{equation}
since $Df_{k+1}(x')\cdot v \in T'_{k+1}(x')$ by definition, we see that
$$
\dist(w,T'_k(x')) \leq |Df_{k+1}(x)\cdot v - Df_{k+1}(x')\cdot v| 
\leq C \wt\varepsilon_k(x) \, r_{k+1}^{-1} |x-x'| |w|.
$$
Similarly, $\dist(w',T'_{k+1}(x)) \leq C \wt\varepsilon_k(x) r_{k+1}^{-1} |x-x'| |w'|$
for $w'\in T'_{k+1}(x')$. It is easy to deduce \eqref{c3} from this, because both spaces
are $d$-dimensional. Our claim \eqref{c2} follows.

The estimates when $f_k(x) \in V_k^{10} \sm V_k^{8}$ are less glamorous; we use 
\eqref{b38} instead of \eqref{b37} and get that
\begin{equation} \label{c5}
|D \pi_{k+1}(x)| \leq C r_k^{-1}.
\end{equation}
When $f_k(x) \in \Sigma_k \sm V_k^{10}$, we use \eqref{b39} and get that
\begin{equation} \label{c6}
|D \pi_{k+1}(x)| \leq C r_l^{-1},
\end{equation}
where $l$ is the last index for which $f_l(x)$ lies in the active or dying zone.

\ms
Now we turn to the 
field of linear isometries which is constructed in
Proposition~9.3 of \cite{DT}. Let $\cR$ be the set of linear isometries of $\R^n$.
There exist $C^1$ functions $R_k : \Sigma_0 \to \cR$, with the following main properties:
\begin{equation} \label{c7}
R_k(\R^d) = T_k(x) \ \text{ for } x\in \Sigma_0;
\end{equation}
\begin{equation} \label{c8}
|R_{k+1}(x)-R_{k}(x)| \leq C \varepsilon
\ \hbox{ for $x \in \Sigma_0$ and } k \geq 0;
\end{equation}
\begin{equation} \label{c9}
|D  R_{k+1}(x)| \leq C_1 r_k^{-1} \varepsilon
\ \hbox{ for $k \geq 0$ and } x \in \Sigma_0.
\end{equation}
For \eqref{c9}, it comes from (9.33) there, which we can simplify because
$f_k$ is bilipschitz in the present situation.

We like these estimates, but want to improve them in many places to include Carleson bounds
that use the $\varepsilon''_k$. For this we want to use our bounds on the $D\pi_k$,
and the way the $\pi_k$ are used to produce the $R_k$.

Let us recall how this goes. We start with $R_0 = I$. Then we suppose that $R_k$ was 
already constructed, and start with a first approximation $S_k$, defined by 
\begin{equation} \label{c10}
S_k(x) = \pi_{k+1}(x) \circ R_k(x) \circ \pi_0 
+ \pi_{k+1}^\perp(x) \circ R_k(x) \circ \pi_0^\perp,
\end{equation}
where $\pi_{k+1}^\perp(x) = I - \pi_{k+1}(x)$ is the orthogonal projection in the 
direction orthogonal to $T_k(x)$, and we set $\pi_0 = p_0(x)$ (the projection on $P_0$)
and $\pi_0^\perp = I - \pi_0^\perp$.
This is the same formula as (9.34) there, with just minor changes in the notation. 

The now usual computation on composition, together with \eqref{c2}, yield that
\begin{equation} \label{c11}
|DS_k(x)| \leq 2 |D R_k(x)| + C r_k^{-1} \wt\varepsilon_k(x)
\ \text{ when } f_k(x) \in V_k^8.
\end{equation}
Then we look at (8.43) there, which says that $R_{k+1}(x)$
is obtained from $S_k(x)$ by the simple rule
\begin{equation} \label{c12}
R_{k+1}(x) = H(S_k(x)),
\end{equation}
where $R$ is a simple nonlinear projection from a set $U$ of linear transformations that
are almost isometries, to the set $\cR$ of linear isometries. This projection $R$ is given
by a reasonably simple formula, but the main point here is that by (9.45) there it is 
$(1+10^{-2})$-Lipschitz on $U$, where $S_k(x)$ takes its values. As a consequence,
\eqref{c11} implies that
\begin{equation} \label{c13}
|DR_{k+1}(x)| \leq 3 |D R_k(x)| + C r_k^{-1} \wt\varepsilon_k(f_k(x))
\ \text{ when } f_k(x) \in V_k^8.
\end{equation}
Recall that when $f_k(x) \in V_k^8$, Lemma \ref{Lb2} says that this happened
also for the previous indices. Then the same induction computations as for \eqref{b37}
yields
\begin{equation} \label{c14}
|r_{k+1} DR_{k+1}(x)| \leq C \wh\varepsilon_k(x),
\ \text{ where } 
 \wh\varepsilon_k(x) = \sum_{l \leq k} 2^{l-k} \varepsilon''_l(f(x)).
\end{equation}

This is a good complement to \eqref{c9}, and now let us see how we may improve 
the estimate \eqref{c8} on $|R_{k+1}-R_{k}|$. We claim that 
\begin{equation} \label{c15}
|R_{k+1}(x)-R_{k}(x)| \leq C\varepsilon''_k(f(x)) \leq C \wt\varepsilon_k(x)
\ \text{ when $x\in \Sigma_0$ and } f_k(x) \in V_k^8,
\end{equation}
where the second part follows at once from the definition \eqref{b37}.
So suppose that $y = f_k(x) \in V_k^8$. Choose $i \in J_k$ such that
$|y-x_{i,k}| \leq 10r_k$; then (7.19) there says that
\begin{equation} \label{c16}
{\rm Angle}(T_k(x),P_{i,k}) \leq C\varepsilon'_k(y),
\end{equation}
where $\varepsilon'_k(y)$ is defined by (7.16) there (recall that $T_k(x) = T\Sigma_k(y)$).
Similarly, (7.10) there 
says that
\begin{equation} \label{c17}
{\rm Angle}(T_{k+1}(x),P_{i,k}) 
= {\rm Angle}(T\Sigma_{k+1}(\sigma_k(y)),P_{i,k}) \leq C\varepsilon_k(y),
\end{equation}
where this time $\varepsilon_k(y)$ is defined in (7.7) there, but we can forget about this 
because it is noted in (1.17) there that $\varepsilon_k(y) \leq \varepsilon'_k(y)$. Thus 
\begin{equation} \label{c18}
{\rm Angle}(T_k(x),T_{k+1}(x)) \leq C\varepsilon'_k(y).
\end{equation}
To be honest, we did not define the angles above, and \cite{DT}
is not much more precise;
however all our angles here are small, and they are equivalent to, for instance, the norm
of the difference of orthogonal projections on the vector spaces parallel to the two spaces 
that we consider. That is,
\eqref{c18} can be taken to mean that
\begin{equation} \label{c19}
|\pi_k(x)-\pi_{k+1}(x)| \leq C\varepsilon'_k(y).
\end{equation}
When we compare $\varepsilon'_k$ and $\varepsilon''_k$, we see that the only difference
is that $\varepsilon''_k$ reaches $r_k$ further, which means that
$\varepsilon'_k(y) \leq \varepsilon''_k(z)$ for any $z\in B(y,r_k)$. This is very
convenient, because this allows us to take $z = f(x)$ 
(by \eqref{b36} and because $y=f_k(x)$); 
then \eqref{c19} implies that
\begin{equation} \label{c20}
|\pi_k(x)-\pi_{k+1}(x)| \leq C\varepsilon''_k(f(x)).
\end{equation}
Next we use the definition \eqref{c10} to estimate $|S_{k}(x)-R_{k}(x)|$.
Since $R_k$ sends $\R^d$ to $T_k(x)$ and hence its orthogonal complement $\R^{n-d}$
to $T_k(x)^\perp$, we see that 
$$ 
R_k(x) = \pi_{k}(x) \circ R_k(x) \circ \pi_0 + \pi_{k}^\perp(x) \circ R_k(x) \circ \pi_0^\perp,
$$
and then \eqref{c10} implies that 
\begin{eqnarray} \label{c21}
|S_k(x) - R_k(x)| &\leq& |\pi_{k+1}(x)-\pi_{k}(x)|+|\pi_{k+1}^\perp(x)-\pi_{k}^\perp(x)|
\nn\\
&=& 2 |\pi_{k+1}(x)-\pi_{k}(x)| \leq C\varepsilon''_k(f(x)).
\end{eqnarray}
This is good, because \eqref{c12} says that $R_{k+1}(x) = H(S_k(x))$
for a Lipschitz mapping $H$ such that $H(R_k(x) = R_k(x)$ because 
$R_k$ is a linear isometry (check with the definition (9.44) there); \eqref{c15} follows.

\section{The mapping $g$ and its Jacobian matrix}
\label{Sd}

\ms
We are now finally ready to define the mapping $g$. We shall keep $g=f$ on $\Sigma_0 = \R^d$,
and now we define $g$ on $\Omega_0 = \R^n \sm \Sigma_0$. Since $g$ will be a
bilipschitz mapping of $\R^n$, it will map $\Omega_0$ to $\Omega = \R^n \sm \Sigma$,
where $\Sigma = f(\Sigma_0)$.

In this section the generic point of $\Omega_0$ is denoted by $(x,t)$,
with $x\in \R^d$ and $t\in \R^{n-d}\sm \{ 0 \}$. We set
\begin{equation} \label{c22}
g(x,t) = \sum_{k\geq 0} \rho_k(t) \Big\{ f_k(x) + R_k(x)\cdot t \Big\}
\ \text{ for } (x,t) \in \Omega_0,
\end{equation}
where  the $\rho_k$ form a partition of $1$ that will be discussed shortly.
This is the same thing as (10.14) or (10.19) there, but some things were simplified, 
because here our initial surface $\Sigma_0$ is just $P_0 = \R^d$, so the projections 
$p$ and $q$ are just the projections $\pi_0$ and $\pi_0^\perp$ on $\R^d$ and $\R^{n-d}$.

The functions $\rho_k$ are defined by (10.15)-(10.17) there, They are radial, 
\begin{equation} \label{dd2}
\sum_{k \geq 0} \rho_k(t) = 1 \ \text{ for } t\in \R^{n-d} \sm \{ 0 \}
\end{equation}
(by (10.16) there), $\rho_0(t) = 0$ when $|t| \leq 1$, and (by (10.17) there) for $k \geq 1$,
\begin{equation} \label{dd3}
\rho_k(t) = 0 \ \text{ unless } r_k \leq |t| \leq 20 r_k. 
\end{equation}
Thus, for each $t$, there are at most three consecutive $k \geq 0$ such that
$\rho_k(t) \neq 0$.

Notice that $g$ does roughly what was announced at the beginning of the previous section:
we start from $f_k(x)\in \Sigma_k$ and go in the orthogonal direction for about $|t|$.
The fact that we actually use an average of up to three different $R_k(x)$ does not matter 
much, because \eqref{c8} and \eqref{c15} say that they are almost the same.
And we are happy that we do not need to take a limit this time.

We want to use the change of variable $g : \Omega_0 \to \Omega$ to reduce the
study of some degenerate elliptic operators $L$ on $\Omega$ to the study of
operators $L_0$ on $\Omega_0$, and because of this we are interested
in the structure of the the matrix of the differential mapping $Dg : \Omega_0 \to \Omega$.

As in \cite{DFM3}, we prefer to study the matrix $J(x,t)$ of $Dg(x,t)$ in a set of two
orthonormal bases of $\R^n$, where the first one is the canonical basis of $\R^n$,
and the second one its image by $R_k(x)$, where $k=k(t)$ is chosen such that
$\rho_k(t) \neq 0$. It does not really matter much which one, but for the sake of
definiteness, let us choose $k(t)$ as large as possible. 
Let us denote $\Jac(x,t)=Dg(x,t)$ in the usual Euclidean basis and $J(x,t):=Dg(x,t) Q(x, t)$ where $Q(x, t)$ is our matrix of isometry in the sense that $R_{k(t)}(x)(y, s)=(y,s) Q (x, t)$ for $(y, s)\in \RR^n$.
We know, just because $g$ is bilipschitz, that for small $\eps_1>0$ the matrix $J(x)$ is (uniformly) bounded and invertible, 
with a (uniformly) bounded inverse, and we are mostly interested now in the block structure of $J$
(when we cut $\R^n$ into $\R^d \times \R^{n-d}$).

\begin{pro}\label{Pd1}
We can write a decomposition of $J$ as a block matrix
\begin{equation} \label{d1}
J(x,t)= \left( \begin{array}{cc}
A^1(x,t) & C^2(x,t) \\ C^3(x,t) & 
I_{n-d}+C^4(x,t)  
\end{array} \right), 
\end{equation}
where the $d \times d$ matrix $A^1$ is bounded, $C^2$, $C^3$, and $C^4$ 
are bounded and satisfy Carleson measure conditions, and $I_{n-d}$ is our notation 
for an identity matrix of size $n-d$. Specifically, 
\eqref{d2} holds with a 
constant $C = C_1 (\varepsilon + \delta_1)$, where $C_1$ depends on 
$n,d,C_0, M_0, M$ (but obviously not on $\delta_1$ or $\varepsilon$), and also
we have the $L^\infty$ estimate
\begin{equation}\label{eq5.7a}
|C^2(x,t)| + |C^3(x,t)| + |C^4(x,t)| \leq C_1 (\varepsilon + \delta_1)
\end{equation}
for $(x,t) \in \R^n$.
\end{pro}

\begin{rem} \label{rk67}
We claim that with the help of Remark \ref{r453}, we can show that 
$A^1$ is as close to the identity matrix as we want, as long as we take $\varepsilon$
and $\delta_1$ small enough. 
This is because $g$ is biLipschitz with a constant close to $1$, and hence
$Dg$ is as close as we want to an isometry,
But in this paper we shall content ourselves with the fact that
$J(x,t)$ is uniformy biLipschitz because $g$ is. 
\end{rem}

\bp The proof of the proposition will keep us busy for some time. 
We first consider the $t$-derivatives of $g$. 
Let us compute $\d_1 g(x,t)$, where $\d_1$ is our notation for $\frac{\d}{\d t_1}$. 
Here we single out the first $t$-variable, because this way we do not have an extra index, 
but the other $t$-derivatives would be the same. From \eqref{c22} we deduce that
\begin{equation} \label{d4}
\d_{1}g(x,t) = D_1 + D_2
\end{equation}
where 
\begin{equation} \label{d5}
D_1 = \sum_l \d_1(\rho_l(t)) \Big\{ f_l(x) + R_l(x)\cdot t \Big\}
\end{equation}
and 
\begin{equation} \label{d6}
D_2 = \sum_{l} \rho_l(t) R_l(x)\cdot e_{d+1},
\end{equation}
where $e_{d+1}$ is the first element of the basis of $\R^{n-d}$.
In both term, the sum has at most $3$ terms, corresponding to $l=k$, $k-1$, 
by the comment below \eqref{dd3} and our choice of $k$ as the largest index
for which $\rho(t) \neq 0$.

We start with $D_1$, notice that $\sum_l \d_1(\rho_l(t)) = 0$ because 
$\sum_l \rho_l(t) = 0$ (see \eqref{dd2}), use this to subtract $f_k(x) + R_k(x)\cdot t$, 
and get that
\begin{eqnarray} \label{d7}
|D_1| &\leq& \sum_l |\d_1(\rho_l(t))| \Big\{ |f_l(x)-f_k(x)| 
+ \big|[R_l(x)-R_k(x)]\cdot t \big|\Big\} 
\nn\\
&\leq& C |t|^{-1} \sum_{k-2 \leq l \leq k-1} \Big\{ |f_{l}(x)-f_k(x)| + |t| \, |R_{l}(x)-R_k(x)| \Big\}
= D_{11} + D_{12},
\end{eqnarray}
where some of the terms may not exist. That is, if $k =0$, then there was only one term
in the initial sum, coming from $k=0$, we managed to kill it, and thus $D_1 = 0$.
Similarly, if $k = 1$, we are left with only one term, coming from $l=0$.
And it could be that even when $k \geq 2$, we don't need $l = k-2$, but the extra
term will not hurt.

\ms
We start with $D_{12}$ because the needed estimates are more recent. 
When $f_{k-1}(x) \in V_{k-1}^8$ (the active region), \eqref{c15} says that
\begin{equation} \label{d8}
|R_{k}(x)-R_{k-1}(x)| \leq C \wt\varepsilon_{k-1}(x).
\end{equation}
Then (by Lemma \ref{Lb2}) $f_{k-2}(x) \in V_{k-2}^8$ (if $k \geq 2$, otherwise there
is no term with $l=k-2$), and
\begin{equation} \label{d9}
|R_{k-1}(x)-R_{k-2}(x)| \leq C \wt\varepsilon_{k-2}(x) \leq C \wt\varepsilon_{k-1}(x).
\end{equation}
If $f_{k-1}(x) \in V_{k-1}^{10} \sm V_{k-1}^8$ (the dying region), 
we replace \eqref{d8} with \eqref{c8}, which says that
\begin{equation} \label{d10}
|R_{k}(x)-R_{k-1}(x)| \leq C \varepsilon.
\end{equation}
In this case Lemma \ref{Lb2} still says that $f_{k-2}(x) \in V_{k-2}^8$, and we can 
use \eqref{d9}.
We are left with the case when $y = f_{k-1}(x)$ lies in the dead region. Then 
$R_{k}(x)=R_{k-1}(x)$ (because $\sigma_{k-1}(y) = y$, $\Sigma_k$ and $\Sigma_{k-1}$ 
coincide at $y$, and the definitions give $S_{k-1}(x) = R_{k-1}(x)$ and then 
$R_{k}(x)=R_{k-1}(x)$. It is still possible that $f_{k-1}(x)$ lies in the dying or active
region, and then we use \eqref{d10} (for $k-1$) or \eqref{d9}. We summarize the cases
and find that
\begin{equation} \label{d11}
D_{12} \leq C \wt\varepsilon_{k-1}(x) + C \varepsilon \delta(x,t),
\end{equation}
where $\delta(x,t) = 1$ if $f_{k-1}(x)$ or $f_{k-2}(x)$ lies in their respective dying region,
and $\delta(x,t) = 0$ otherwise. 
We will see later that this leads to a good Carleson estimate for $D_{12}$.

\ms
Next consider $D_{11} = C |t|^{-1} \sum_{k-2 \leq l \leq k-1}  |f_{l}(x)-f_k(x)|$,
and first assume that $f_{k-1}(x) \in V_{k-1}^8$. Then by \eqref{bb28},
\begin{equation} \label{dd15}
|f_{k}(x)-f_{k-1}(x)| \leq C \varepsilon''_{k-1}(f(x)) \, r_{k-1}. 
\end{equation}
If $k \geq 2$, Lemma \ref{Lb2} says that $f_{k-2}(x) \in V_{k-2}^8$, and \eqref{bb28}
yelds
\begin{equation} \label{dd16}
|f_{k-1}(x)-f_{k-2}(x)| \leq C \varepsilon''_{k-2}(f(x)) \, r_{k-2}.
\end{equation}
Otherwise, we don't need this estimate because $D_{11}$ has only one term. 
Altogether,
\begin{equation} \label{dd17}
D_{11} \leq C \varepsilon''_{k-1}(f(x)) + C \varepsilon''_{k-2}(f(x))
\end{equation}
(because $t \geq r_k$ by \eqref{dd2}).

Next we assume that $f_{k-1}(x) \in \Sigma_{k-1} \sm V_{k-1}^8$.
If $f_{k-1}(x)$ or $f_{k-2}(x)$ lies in the dying region, we use the more brutal
estimate \eqref{b35} to see that $|f_k(x)-f_{k-1}(x)| \leq C \varepsilon$, or similarly for 
$k-1$, and get that $D_{11} \leq C$. Otherwise, $f_{k-1}(x)$ and $f_{k-2}(x)$ lie in their
dead regions, and $D_{11} = 0$. We summarize the estimates as we did above, by saying that
\begin{equation} \label{dd18}
D_{11} \leq C \varepsilon''_{k-1}(f(x)) + C \varepsilon''_{k-2}(f(x)) + C \varepsilon \delta(x,t)
\end{equation}
with the same definition for $\delta(x,t)$ and where we set $\varepsilon''_{l}(f(x)) = 0$
for $l < 0$.

Next we study $D_2$, which we write as $D_2 = D_{21} + D_{22}$, where
\begin{equation} \label{d19}
D_{21} = R_l(x)\cdot e_{d+1}
\end{equation}
and by \eqref{dd2} the rest is
\begin{equation} \label{d6}
D_{22} = \sum_{k-2 \leq l \leq k-1} \rho_l(t) [R_l(x)-R_k(x)] \cdot e_{d+1}.
\end{equation}
The main piece $D_{21}$ gives the first term of the identity matrix $I_{n-d}$
in \eqref{d1}, and of course its analogue for the other $t$-derivatives of $g$
give the rest of $I_{n-d}$. As for $D_{22}$, we use the same estimates 
\eqref{d8}-\eqref{d10} as above and find that
\begin{equation} \label{d21}
D_{22} \leq C \wt\varepsilon_{k-1}(x) + C \varepsilon \delta(x,t).
\end{equation}

\ms
We may now consider the $x$-derivatives of $g$, and 
which we feel couragous enough to consider all at the same time
and denote by $D_x g$. By \eqref{c22}, 
\begin{equation} \label{d22}
D_x g(x,t) = \sum_{j\geq 0} \rho_j(t) \Big\{ D f_j(x) + D R_j(x)\cdot t \Big\}
 = : D_3 + D_4.
\end{equation}
We start with $D_4$. Notice that when $k = 0$, $R_k = I$ and $DR_k = 0$, so
we may assume that $k \geq 1$. As usual begin with the case when 
$y = f_{k-1}(x) \in V_{k-1}^8$, apply \eqref{c14}, and find that 
$|r_{k} DR_{k}(x)| \leq C \wh\varepsilon_{k-1}(x)$, and similarly for 
the previous iterates if they are needed. This yields
\begin{equation} \label{d23}
|D_4| \leq C \wh\varepsilon_{k-1}(x) + C \wh\varepsilon_{k-2}(x) + C \wh\varepsilon_{k-2}(x)
\leq C \wh\varepsilon_{k-1}(x)
\end{equation}
because $|t| \leq 20r_k$ (by \eqref{dd3} and because $k \geq 1$), and 
where the last part comes from the definition of the $\wh\varepsilon_l$ in \eqref{c14}.
When $f_{k-1}(x) $, or one of its two predecessors, lies in the closure of its dying region, 
we use \eqref{c9} for all of them and find that $|D_4| \leq C \varepsilon$. 

We are left with the case where the three points lie in the interior of their dead region.
Denote by $l$ the smallest integer such that $f_l(x)$ lies in the the interior of its dead region; 
thus $l < k-2$. We know that all the $f_m(x)$, $m \geq l$, are equal to $f_l$ in a neighborhood
of $x$, and when we follow the computations we see that this means that the $R_m$,
$m \geq l$, also coincide with $R_l$ near $x$ (for instance, we check first that
since $\sigma_l(y) = y$ near $f_l(x)$, we get that $\Sigma_{l+1} = \Sigma_l$ near
$f_l(x)$, then $\pi_{l+1}=\pi_l$ near $x$, then $S_k = R_k$, etc.).
Now the $DR_j$ in the formula \eqref{d22} are all equal to $DR_l$, and \eqref{c9}
yields
\begin{equation} \label{d24}
|D_4| \leq C \varepsilon r_k r_l^{-1} \leq C \varepsilon 10^{l-k}.
\end{equation}

Again we claim that this decay will lead to a Carleson measure estimate, but let us
now concentrate on our last term
\begin{equation} \label{d25}
D_3 = \sum_{j\geq 0} \rho_j(t)  D f_j(x),
\end{equation}
which as usual we cut in two. The first part
\begin{equation} \label{d26}
D_{31} = \pi_k(x) \circ D_3,
\end{equation}
where we project on the vector plane $R_k(x)(\R^d)$ parallel to $T_k(x)$,
falls in the matrix $A^1$ of the decomposition \eqref{d1}, and we don't need any
special information about it, 
except that we know that $A^1$ is bounded (and even $J$ is bilipschitz). 
We are left with
\begin{equation} \label{d27}
D_{32} = \pi_k(x)^\perp \circ D_3 = \sum_{j\geq 0} \rho_j(t) [\pi_k(x)^\perp \circ D f_j(x)].
\end{equation}
Observe that the image of $D f_j(x)$ is contained in the tangent direction
$R_j(x)(\R^d)$ (parallel to $T_j(x)$), and $\pi_j(x)^\perp$ 
vanishes on this space.
Also, $D f_j(x)$ is bounded, so
\begin{equation} \label{d28}
|D_{32}| \leq C\sum_{k-2 \leq j \leq k} |\pi_j(x)^\perp - \pi_k(x)^\perp|
= C\sum_{k-2 \leq j < k} |\pi_j(x) - \pi_k(x)|.
\end{equation}
The simplest for us is to observe that $|\pi_j(x) - \pi_k(x)| \leq C |R_j(x) - R_k(x)|$.
Indeed $\pi_j(x)$ is the projection on $R_j(x)(\R^d)$, and an orthonormal basis
of that space is given by the $R_j(x)(e_l)$, $1 \leq l \leq d$, so
$\pi_j(x)(v) = \sum_{l=1}^d \langle v, R_j(x)(e_l) \rangle R_j(x)(e_l)$.
Of course using this is a little strange, because the estimates 
\eqref{c15} and \eqref{c8} that we are about to use come from estimates
on $|\pi_j(x) - \pi_k(x)|$, as in \eqref{c19}.
Anyway, $|D_{32}| \leq C \sum_{k-2 \leq j < k} |R_j(x) - R_k(x)|$ can now be estimated
exactly as $D_{22}$ and $D_{11}$, and we get that
\begin{equation} \label{d29}
D_{22} \leq C \wt\varepsilon_{k-1}(x) + C \varepsilon \delta(x,t).
\end{equation}
as in \eqref{d21}.

We completed the decomposition of $J$; now we need to show that 
the error terms in \eqref{d11}, \eqref{dd18}, \eqref{d21}, \eqref{d23},
\eqref{d24}, and \eqref{d29} are functions that satisfy a Carleson measure estimate.

We start with the function $\wh \varepsilon_{k-1}(x)$ that show up in \eqref{d23}.
Recall from \eqref{c14} that $\wh \varepsilon_{k-1}(x) = \sum_{l\leq k-1} 
2^{l-k+1} \varepsilon''_l(f(x))$. This is a function of $t$ as well, because since $k=k(t)$
(and $k \geq 1$), 
\eqref{dd3} says that $r_k \leq |t| \leq 20r_k$. We will use the fact that
$\wh \varepsilon^2_{k-1}(x) \leq C \sum_{l\leq k-1} 2^{l-k+1} \varepsilon''_l(f(x))^2$ 
by Cauchy-Schwarz.

In order to prove \eqref{d2} for this function, we have to estimate
\begin{equation} \label{d30}
I(X,R) = \int_{\Omega_0 \cap B(X,R)} |\wh \varepsilon_{k-1}(x)|^2 \frac{dydt}{|t|^{n-d}} 
\leq C \int_{\Omega_0 \cap B(X,R)} \sum_{0 \leq l\leq k-1} 2^{l-k} \varepsilon''_l(f(x))^2 
\frac{dydt}{|t|^{n-d}}.
\end{equation}

First observe that since $r_k \leq |t| \leq 20r_k$ and $|t| \leq R$ when $(y,t) \in B(X,R)$
(recall that here $X = (x,0)$ lies in $P_0$), we only sum over $k$ such that $r_k \leq R$.
Let us fix $x$, $k$, and $l$, and integrate in $t$ first. We integrate in the region
$A(k,l)$ where $r_k \leq |t| \leq 20r_k$, and $\int_{A(k,l)} 
\frac{dt}{|t|^{n-d}} \leq C$.
We are left with 
\begin{equation} \label{d30}
I(X,R) \leq  C \int_{y\in P_0 \cap B(X,R)} \sum_{0 \leq l\leq k-1} 
2^{l-k} \varepsilon''_l(f(x))^2  dy.
\end{equation}
We now sum over $k$, then $l$. The sum over $k$ disappears because of the converging 
factor $2^{l-k}$, and the sum over $l$ is 
less than $C \delta_1$, by \eqref{b13} and \eqref{b14} 
(recall that $g(y) = f(y)$ on $P_0$). We are left with
$I(X,R) \leq  C \delta_1 \H^{d}(P_0 \cap B(X,R)) \leq C \delta_1 R^d$, as needed.

The numbers $\wt \varepsilon_{k-1}(x)$ are just smaller than the $\wh \varepsilon_{k-1}(x)$,
so we don't need to worry about them; the same is true of $\varepsilon''_{k-1}(f(x)$
and $\varepsilon''_{k-2}(f(x)$, which are just two pieces of $\wh \varepsilon_{k-1}(x)$.
The next function to control is $\varepsilon \delta(x,t)$, which counts whether $f_{k-1}(x)$ or 
$f_{k-2}(x)$ lies in the dying zone, or its variant where we also include $f_{k-2}(x)$
which is implicit in the description below \eqref{d23}. We need to control
\begin{equation} \label{d31}
I'(X,R) = \int_{\Omega_0 \cap B(X,R)}  \varepsilon \,
\delta(x,t) \frac{dydt}{|t|^{n-d}} .
\end{equation}
We know from Lemma \ref{Lb2} that for a given $x$, there is at most one $l \geq 0$
such that $f_l(x) \in V_l^{10} \sm V_l^8$, and the only $t$ for which 
$l \in \{ k-1, k-2, k-3 \}$ are such that $r_k \leq |t| \leq 20r_k$ for 
$k \in \{ l+1, l+2, l+3 \}$. That is, $t \in [10^{-3}r_l, 2r_l/10]$. We integrate against
$|t|^{d-n} dt$, get at most 
$C \varepsilon$, then integrate against $y$ and get at most 
$C \varepsilon R^d$, as needed.

Our last contribution comes from \eqref{d24}, where for some earlier $l = l(x) < k-2$,
$f_l(x)$ lies in the dead region for the first time, and then we pay 
$D_4 \leq C \varepsilon 10^{l-k}$.
This yields the integral
\begin{equation} \label{d31}
I''(X,R) = \varepsilon
  \int_{y\in Z \cap B(X,R)}\int_{0 < |t| \leq R} 
\,\sum_{k(t) \geq l(x)} 100^{l-k}
\, \frac{dydt}{|t|^{n-d}} ,
\end{equation}
where $Z$ is the set of points $y \in P_0$ such that $l(y)$ exists.
As before, we integrate first against the $t$ such that $k = k(t)$ and get a constant,
then sum in $k$ and get another constant, and finally integrate in $y$ and get at most 
$C \varepsilon R^d$.

This completes out verification that the functions in our various estimates satisfy a Carleson
condition, as announced with a constant dominated by $C (\varepsilon + \delta_1)$. 
The $L^\infty$ bound \eqref{eq5.7a} is easier (we don't even have to sum the terms);
Proposition \ref{Pd1} follows.
\qed

\ms
It will be good to know that the class of matrices that have the special form given in 
Proposition \ref{Pd1} is stable under taking inverses, products, and transposes.
Indeed we start from our favorite operator $L=-\div D_\alpha(Y)^{d+1-n} \nabla$,
and then we use $g$ to change variables and get an operator $L_0$ on
the simpler domain $\Omega_0 = \R^n \sm \R^d$. A fairly standard computation, 
which the reader may find in \cite{DFM3}, Lemma 6.17, shows that 
$L_0 = {\rm div} A \nabla$, where the matrix of $A$ is
\begin{equation} \label{d3}
A(x,t) = (D_\alpha(g(x,t))^{-(n-d-1)} \, |\det J(x,t)| (J(x,t)^{-1})^T J(x,t)^{-1};
\end{equation}
maybe the reader expected the same formula with $J(x,t)$ replaced with $Dg(x,t)$,
but $J(x,t) = Dg(x,t) Q(x,t)$ and the isometry $Q(x,t)$ does not change 
$(J(x,t)^{-1})^T J(x,t)^{-1}$ nor the determinant (and $J$ has a simpler form!).
The next lemma will thus tell us that Theorem \ref{Itsf1} can be applied to
$L_0$ and $A$.

\begin{lem}\label{pRd3} 
Denote by $\M(M,\tau, K)$ the class of matrix-valued functions $J$ that have a decomposition \eqref{d1}, where $A^1$ is bounded and invertible, with a bounded inverse
such that  $\|A^1\|_\infty+\|(A^1)^{-1}\|_\infty\leq M$, each  $C^i$, $i=2,3,4$,
satisfies the Carleson bound \eqref{d2} with $C = K$, and 
$\|C^i\|_\infty\leq \tau$ for $i=2,3,4$. Also assume that $\tau < (6M)^{-1}$.
Then $\M(M,\tau, K)$ is stable under taking the transposed matrices,
$J^{-1} \in \M(M^2, 6 M^2 \tau, 36M^4 K)$ for $J \in \M(M,\tau, K)$,
and $JJ' \in \M(M^2, 8 M^2 \tau, 64M^4 K)$ for $J, J' \in \M(M,\tau, K)$.

Here, $\|\cdot\|$ is the norm of the associated operator acting on the Euclidean $\R^n$ and 
$\|\cdot\|_\infty$ is its supremum in $x,t$.
\end{lem}

We decided to compute invertibility on the block matrix $A^1$ rather than 
the full matrix $J$, because this is easier, but as soon as the $C^i$ are small enough, 
there is no real difference. 
That is, we know that $g$ is bilipschitz, so $J =J(x,t)$ is invertible, 
with some uniform bound $\wt M$ on $||J(x,t)^{-1}||$. Of course $\wt M$ does not depend
on $\varepsilon$ or $\delta_1$; taking these constants smaller only makes our assumptions on the stopping time region $\C$ harder to check. Set
\begin{equation} \label{dd37}
T = \left( \begin{array}{cc} A_1 & 0 \\ 0 & I_{n-d} \end{array} \right)
\text{ and } E = J-T; 
\end{equation}
then $||E||$ is as small as we want, by \eqref{eq5.7a}, and $T$ and $A_1$ are invertible too, with
$||T^{-1}|| = ||T^{-1}|| \leq 2\wt M$, so we can apply the lemma to $J$,
with $M = 2 \wt M$ and $\tau = C(\varepsilon + \delta_1)$. 

Because of the way we once used $\delta_1$ to control some geometric
quantities that should really have been controlled by $\varepsilon_1$,
we have to take $\delta_1$ small. Also, we decided not to use the various remarks 
leading to Remark \ref{rk67}, so we do not know officially that our changes of variable are in fact
bilipschitz with a constant near $1$.
So our argument is sound, but not optimal.

\ms

\bp
Let us now prove the lemma; the verification will be mostly a pointwise thing.
The fact that $J^T \in \M(M,\tau, K)$ when $J \in \M(M,\tau, K)$ is clear;
let us now consider the inverse of $J \in \M(M,\tau, K)$. Define $T$ and $E$
as in \eqref{dd37}, and observe that 
$T^{-1}$ is a block matrix like $T$
(the associated linear mapping acts as we want on $\R^d$ and $\R^{n-d}$)
\begin{equation} \label{g2}
J^{-1} = (T+E)^{-1} = [T(I+T^{-1}E)]^{-1} = (I+T^{-1}E)^{-1}T^{-1}
\end{equation}
(where the invertibility of $J$ follows from the computation below) and then 
use our assumption that $\tau \leq (6M)^{-1}$ to
write $(I+T^{-1}E)^{-1}$ as a Neuman series. This gives
\begin{equation} \label{g3}
||J^{-1}-T^{-1}|| \leq 2 ||T^{-1}|| ||T^{-1}E|| \leq 2M^2 ||E|| 
\leq 6M^2 \tau,
\end{equation}
which means that $J^{-1}$ has a nice block decomposition. 
This give pointwise bounds, and the $L^\infty$ bounds follow by taking supremums.
The Carleson estimate also follows directly from \eqref{g3}; 
we need squares because $C$ in Definition~\ref{cm} is quadratic.
Now we consider the product with another matrix 
$J' \in \M(M,\tau, K)$.
Write $J' = E'+T'$, with similar notation, notice that
\begin{equation} \label{g4}
JJ' = (E+T)(E'+T') = TT' + R
\end{equation}
where $TT'$ is a diagonal block matrix like $T$ and $T'$ (again look at the 
corresponding endomorphism), and 
$||R|| \leq 6\tau M + 9\tau^2 \leq 6\tau M$.
As before, the $L^\infty$ and Carleson bounds follow.
\ep

\section{Distance functions like $D_\alpha$}
\label{SDist}

The next string of estimates concerns the Carleson behavior of two things that are related.
The distance function $D = D_\alpha$, associated to our final set $\Sigma$, or in fact
any uniformly rectifiable set $E$, and a control function $\lambda(x,r)$ for the average density
of an Ahlfors-regular measure living on $E$. Later on, we shall study relations between
two distance functions, typically one coming from $E$ and one coming from our approximating surface.

\subsection{The function $\lambda$.}

Let $E$ be a uniformly rectifiable set $E$, and $\sigma$ any Ahlfors regular measure
supported on $E$. We want to define a $\lambda(x,r)$ that will measure,
in a reasonable smooth way, the density of $\sigma$.

Pick a smooth radial, nonnegative function $\eta$, supported in the unit ball of $\R^n$,
with $\int \eta = 1$, and set $\eta_t=\frac{1}{t^d}\,\eta \left(\frac xt\right)$ 
(notice the normalization adapted to $\R^d$). We will use the $\eta_t$ for different things.

For $x\in E$ and $r > 0$, we define a first, not too precise, measure of the density, namely
\begin{equation} \label{2.1}
\lambda_0(x,r) = \int_{E \cap B(x,r)} \eta_r(y-x) d\sigma(y) > 0,
\end{equation}
and then a center of mass
\begin{equation} \label{2.2}
\begin{aligned}
\Phi_r(x) &= \lambda_0(x,r)^{-1} \int_{E \cap B(x,r)} y \,\eta_r(y-x) d\sigma(y)
\cr&= x + \lambda_0(x,r)^{-1} \int_{E \cap B(x,r)} (y-x) \,\eta_r(y-x) d\sigma(y),
\end{aligned}
\end{equation}
where we write the second formula to insist on the translation invariance,
and finally the better density

\begin{equation} \label{2.3}
\lambda(x,r) = \int_{E \cap B(\Phi_r(x),r)} \eta_r(y-\Phi_r(x)) d\sigma(y).
\end{equation}

We prefer to 
use $\lambda$ rather than $\lambda_0$,
because maybe $x$ itself lies far from an optimal plane for $\alpha(x,10r)$ (defined below).

Anyway, we want to show that
\begin{equation} \label{2.4}
r |\nabla_{x,r} \lambda(x,r)| \in 
CM(E \times \R_+), 
\end{equation}
where we define $CM(E \times \R_+)$ as in Definition \ref{cm}, 
but with $\R^d$ replaced with $E$.
That is, we say that a function $F(x,t)$, defined on $E \times (0,+\infty)$,
satisfies the Carleson measure condition, and write $F\in CM(E \times \R_+)$,
when there is a constant $C \geq 0$ such that 
\begin{equation} \label{d2bis}
\int_{x \in E \cap B(X,R)}\int_{t\in(0,R)} |F(x,t)|^2 \, \frac{d\sigma(x)dt}{t}
\leq C R^d
\end{equation}
for $X \in E$ and $R>0$. We could replace $\sigma$ with $\H^d_{\vert E}$ without 
changing the class $CM(E \times \R_+)$.

The logical plan for proving \eqref{2.4} will work: for each $(x,r) \in E \times \R_+$, we 
will find a flat measure $\mu$ that approximates $\sigma$ well in 
$B(x,10r)$, and compare the three quantities above to the same ones with $\sigma$ replaced by $\mu$.
The good approximation will be in terms of the Tolsa numbers $\alpha(x,r)$, 
which we discuss now. 

We will use the same definition of $\alpha(x,r)$ is the same as in \cite{DFM3}. 

We first define flat measures and local Wasserstein distances. 
Denote by $\bP$ the set of affine $d$-planes in $\R^n$, and for each
plane $P \in \bP$, denote by $\mu_P=\H^d|_P$ the restriction of $\H^d$ to $P$
(in other words, the Lebesgue measure on $P$). By flat measure, we shall 
mean simply mean a measure $c \mu_P$, with $c > 0$ and $P\in \bP$.
The number $\alpha(z,r)$ will measure 
the distance between our measure $\sigma$ and flat measures, 
locally in the 
ball $B(z,r)$, which we shall often take centered on $E$ because this way we know 
that $\mu(B(z,r))$ is fairly large.

\begin{definition} \label{D6.1}
For $z\in \R^n$ and $r > 0$, denote by $Lip(z,r)$ the set of Lipschitz functions 
$f : \R^n \to \R$ such that $f(y)=0$ for $y\in \R^n \sm B(z,r)$ and $|f(y)-f(w)|\leq |y-w|$
for $y,w\in \R^n$. Then define the normalized 
(local)
Wasserstein distance between two measures
$\sigma$ and $\mu$ by
\begin{equation} \label{a6.2}
\dist_{z,r}(\sigma, \mu) = r^{-d-1} \sup_{f\in Lip(z,r)} \bigg|\int f d\sigma - \int f d\mu\bigg|
\end{equation}
and the local distance from $\sigma$ to flat measures by
\begin{equation} \label{2.5bis}
\alpha(z,r) = \inf_{c\geq 0, \, P\in \bP}\dist_{z,r}(\sigma,c\mu_P).
\end{equation}
\end{definition}

We normalized $\dist_{z,r}(\sigma, \mu)$ with $r^{-d-1}$ because this way,
if $\mu(B(z,r)) \leq C r^d$ and $\sigma(B(z,r)) \leq C r^d$, then 
$\dist_{z,r}(\sigma, \mu) \leq 2C$ because
\begin{equation} \label{a6.4}
\|f\|_\infty \leq r \ \text{ for } f\in Lip(z,r).
\end{equation}
Also observe that if $B(y,s) \subset B(z,r)$, then $Lip(y,s) \subset Lip(z,r)$; it follows that
$\dist_{y,s}(\sigma, \mu) \leq (r/s)^{d+1} \dist_{z,r}(\sigma, \mu)$, and hence
\begin{equation} \label{a6.5}
 \alpha(y,s) \leq (r/s)^{d+1} \alpha(z,r) 
\ \text{ when } B(y,s) \subset B(z,r).
\end{equation}

Return to the proof of \eqref{2.4}. We want to show that 
\begin{equation} \label{2.5}
r |\nabla_{x,r} \lambda(x,r)| \leq C \alpha(x,10r)
\text{ for $x\in E$ and $r > 0$,}  
\end{equation}
because then \eqref{2.4} will follow from Theorem~1.2 in\cite{To},
which says that when $E$ is uniformly rectifiable and $\sigma$ is any
Ahlfors regular measure on $E$, then 
\begin{equation} \label{2.4t}
(x,r) \mapsto \alpha(x,r)  \in CM(E \times \R_+).
\end{equation}
Strictly speaking, \cite{To} defines the function $\alpha$ slightly differently, 
on the set of dyadic cubes in $\R^n$ rather than balls centered on $E$.
But the diffefence is really minor, in the sense that one quantity controls the other,
and we refer to Lemma 5.9 in \cite{DFM3} for the verification.

It is easy to see that $|r \nabla_{x,r} \lambda(x,r)| \leq C$, so 
\eqref{2.5} is trivial when $\alpha(x,10r) \geq C^{-1}$. Therefore we may assume that
$\alpha(x,10r) \leq C^{-1}$, with $C$ as large as we want.
Choose an almost optimal flat measure $\mu = a \H^d_{\vert P}$ in the definition of 
$\alpha(x,10r)$, where of course $P$ is a $d$-plane and $a > 0$. 
We do not intend to use the fact that $\mu$ is nearly optimal here, just 
that its distance to $\sigma$ is small. That is, if this distance was some small $\beta > 0$, 
we would just get \eqref{2.5} with $C\beta$.
The most trivial application of this (obvious) remark is that we may use $\alpha(x,10^7r)$
instead, or use some other numbers and planes. 

With the assumption that $\alpha(x,10r) \leq C^{-1}$
(and by testing for example the definition against a multiple of  $\eta_r$), 
we get that $C^{-1} \leq a \leq C$, for some $C$ that depends on $\eta$ 
and the regularity constant for $\sigma$.

Set $\alpha = \alpha(x,10r)$ to save energy. 
Also write $\wt\lambda_0(x,r)$ for the analogue of $\lambda_0(x,r)$, but with $\mu$,
and do the same thing for $\wt \Phi_r(x)$ and $\wt\lambda(x,r)$.
We use the definition of $\alpha$ and find easily that 
\begin{equation} \label{2.6}
|\wt\lambda_0(x,r)-\lambda_0(x,r)| \leq C \alpha
\end{equation}
and  then, with just a bit more of computation, that
\begin{equation} \label{2.7}
|\wt\Phi_r(x)-\Phi_r(x)| \leq C r \alpha.
\end{equation}
Then we first try to differentiate with respect to $x$, i.e. estimate
\begin{equation} \label{2.8}
\nabla_x \lambda(x,r) =  \int_{E} \nabla_x[\eta_r(y-\Phi_r(x))] d\sigma(y),
\end{equation}
and the first thing to do is differentiate $\Phi_r$. Thus we first differentiate 
(with respect to $x$) 
the quantity $\eta_r(y-x) = r^{-d} \eta((y-x)/r)$ 
and we get $- r^{-d-1} \nabla\eta((y-x)/r)$. So, for instance,
\begin{equation} \label{2.9}
\nabla_x 
\lambda_0(x,r) = - r^{-d-1} \int_{E} \nabla\eta((y-x)/r) d\sigma(y)
\end{equation}
and, using the second part of \eqref{2.2},
\begin{eqnarray} \label{2.10}
\nabla_x [\Phi_r(x) -x] &=& 
- \frac{\nabla_x \lambda_0(x,r)}{\lambda_0(x,r)} \, [\Phi_r(x)-x]
- \lambda_0(x,r)^{-1}  \int_{E} \eta_r(y-x) d\sigma(y)
\nn\\
&\,&\hskip0.5cm
 - \lambda_0(x,r)^{-1} r^{-d-1} \int_{E} (y-x) \nabla\eta(y/r - x/r)d\sigma(y).
\end{eqnarray}
We prefer to subtract $x$ because despite the more complicated formula, the 
flawless homogeneity makes it easier to check that $|\nabla_x \Phi_r(x)| \leq C$.
We will also need to know that 
\begin{equation} \label{2.11}
\dist(\Phi_r(x),P) \leq C \alpha r.
\end{equation}
Indeed, it is obvious that $\wt \Phi_r(x) \in P$ because $\mu$ is supported in $P$, so
$\dist(\Phi_r(x),P) \leq |\Phi_r(x) - \wt \Phi_r(x)| \leq C r \alpha$
by \eqref{2.7}.

Let us return to \eqref{2.8}, set $V(x) = \nabla_x \Phi_r(x)$ to save notation, and notice that
$$
\nabla_x[\eta_r(y-\Phi_r(x))]
= - r^{-d-1} \nabla\eta((y - \Phi_r(x))/r) \cdot V(x),
$$
so that \eqref{2.8} becomes
\begin{equation} \label{2.12}
\nabla_x \lambda(x,r) = - r^{-d-1}
 \Big(\int_{E} \nabla\eta((y - \Phi_r(x)) d\sigma(y)\Big)  \cdot V(x),
\end{equation}
where we pulled $V(x)$ out of the integral to stress the fact that it does not depend on $y$.
Since $V(x)$ is bounded, we see that $|\nabla_x \lambda(x,r)| \leq CA$, where
\begin{equation} \label{2.12a}
A = r^{-d-1} \int_{E} \nabla\eta((y - \Phi_r(x)) d\sigma(y).
\end{equation}
Let us compare $A$ with the same expression $A_1$, where we just replace 
$\sigma$ by $\mu$. Notice that the integrand $f(y) = \nabla \eta((y - \Phi_r(x)/r)$
is a nice Lipschitz function supported on $B(x,r)$, with Lipschitz norm less than
$C r^{-1}$, so 
$$
|A-A_1| =  r^{-d-1} \Big| \int_{E} f(y)(d\sigma(y) - d\mu(y))\Big|
\leq C \alpha r^{-1}
$$
by \eqref{a6.2} and the definition of $\mu$, and where we get an extra $r^{-1}$
coming from the Lipschitz norm of $f$. Next set denote by $\xi$ the 
orthogonal projection of $\Phi_r(x)$ on $P$, and consider
$$
A_2 = - r^{-d-1} \int_{P} \nabla\eta((y-\xi)/r) d\mu(y),
$$
where we just replaced $\Phi_r(x)$ by $\xi$ in the definition of $A_1$. 
We claim that $A_2 \cdot V =0$ for every vector $V$ (and hence $A_2 = 0$). 
When $V$ is parallel to $P$, $A_2 \cdot V =0$ because we integrate the partial derivative 
in the direction of $V$ of a function with compact support. 
When instead $V$ is orthogonal to $P$, 
$V(x) \cdot \nabla \eta((y- \xi)/r) = 0$
for every $y\in P$, because $\eta$ is radial and $V$ is orthogonal to the direction of
$y - \xi$. So $A_2(V) = 0$. Finally,
$|A_1(V) - A_2(V)| \leq C r^{-1} |\Phi_r(x)/r-\xi/r| \leq C r^{-1}\alpha$, 
by differentiating again under the integral, between $\Phi_r(x)$ and $\xi$.
Altogether $|\nabla_x \lambda(x,r)| \leq CA \leq C r^{-1}\alpha$; 
this proves the $x$-derivative part of \eqref{2.5}, and indeed the only important properties 
of $\Phi_r(x)$ are that
\begin{equation} \label{2.13}
|\nabla_x \Phi_r(x)| \leq C \ \text{ and  } \dist(\Phi_r(x),P) \leq C \alpha r.
\end{equation}

\smallskip
We still need to take care of $r$-derivatives, and this will work the same way.
This time we need to compute
\begin{equation} \label{2.14}
\d_r \lambda(x,r) =  \int_{E} \d_r[\eta_r(y-\Phi_r(x))] d\sigma(y),
\end{equation}
and we start with the derivatives of $\lambda_0$ and $\Phi_r(x)$ with respect to $r$.
The derivative of $\eta_r(y-x) = r^{-d} \eta((y-x)/r)$ is 
\begin{equation} \label{a721}
\frac{\d}{\d r}\big(\eta_r(y-x)\big) = 
-d r^{-d-1}\eta((y-x)/r) - r^{-d-2} \nabla\eta((y-x)/r) \cdot (y-x),
\end{equation}
which means that for instance
\begin{equation} \label{a722}
\d_r \lambda_0(x,r) = - r^{-d-1} \int_{E} \Big[ d \eta((y-x)/r) 
+ \nabla\eta((y-x)/r) \cdot (y-x)/r  \Big] d\sigma(y).
\end{equation}
Then $r \d_r \lambda_0(x,r)$ is also bounded, as for $\nabla_x \lambda(x,r)$.

Next we study 
$W = \d_r \Phi_r(x) = \d_r [\Phi_r(x)-x]$,
with 
$\Phi_r(x)-x =
\lambda_0(x,r)^{-1} \int_{E \cap B(x,r)} (y-x) \,\eta_r(y-x) d\sigma(y)$.
Recall that $C^{-1} \leq \lambda_0(x,r) \leq C$ and
$\int_{E \cap B(x,r)} (y-x) \,\eta_r(y-x) d\sigma(y) \leq Cr$.
The part $W_1$ where we differentiate $\lambda_0(x,r)^{-1}$
is thus at most $C r |\d_r \lambda_0(x,r)| \leq C$; we are left with
$W_2 = \int_{E \cap B(x,r)} (y-x) \frac{\d}{\d r}\big(\eta_r(y-x)\big) d\sigma(y)$.
We use \eqref{a721} again, get one more power of $r$ than in \eqref{a722},
and it follows that $W$ is bounded.
Finally we return to \eqref{2.14}; compared to the computation for $\lambda_0$,
we get an extra term coming from $W$. That is, 
\begin{eqnarray} \label{2.15}
\d_r \lambda(x,r) &=& - r^{-d-1} \int_{E} \Big[ d \eta((y-\Phi_r(x))/r) 
+ \nabla\eta((y-\Phi_r(x))/r) \cdot (y-\Phi_r(x))/r 
\nn \\ &\,& \hskip 3cm
+ \nabla\eta((y-\Phi_r(x))/r) \cdot W
\Big] d\sigma(y)
\end{eqnarray}
and we just need to estimate
$A = r^{-d-1} \int_{E} \nabla\eta((y-\Phi_r(x))/r)$ because $W$ is bounded.
We are lucky; $A$ is the same as in \eqref{2.12a}, and we proved that 
$|A| \leq C r^{-1} \alpha$, so $|\d_r \lambda(x,r) C r^{-1} \alpha$ as well, and
the full \eqref{2.5} follows. This also completes our proof of \eqref{2.4} (because of \eqref{2.4t}).

\subsection{The distance function $D_\alpha$ versus the distance to a good plane}

Now we take a distance $D=D_{\sigma,\alpha}$ related to $\sigma$, 
and use the $\alpha$ numbers to compare it locally to the distance to a plane.

\begin{lem} \label{l72}
Let $\sigma$ be any Ahlfors-regular measure on any AR set, and define 
$D=D_{\sigma,\alpha}$ by \eqref{eq9.40}. For $x\in E$, $r > 0$, 
any $d$-plane $P = P(x,r)$ that almost minimizes in the definition of $\alpha(x,16r)$, 
$z\in B(x,2r)$ such that 
\begin{equation} \label{727}
\min(\dist(z,P), \dist(z,E)) \geq 10^{-2} r, 
\end{equation}
 we have 
\begin{equation} \label{2.16}
\Big|\frac{D(z)}{\dist(z,P(x,r))} - C_\alpha \lambda(x,r)^{-1/\alpha} \Big| 
\leq C \sum_{l \geq 4}  2^{-\alpha l}\alpha(x,2^l r), 
\end{equation}
where $C$ depends on $n$ and the AR constants for $\sigma$, and
$C_\alpha$ is a dimensional constant that does not depend on $E$ or $\sigma$. 
\end{lem}

We could stop the sum when $2^l r \geq 10$ if we are really talking about 
$\sigma$ and the approximating surface $\Sigma$, but we continue it forever
because we are talking about an arbitrary Ahlfors-regular set $E$ with an 
Ahlfors regular measure $\sigma$ on it. 
We did not require $E$ to be uniformly rectifiable in the statement, but this
assumption will be very useful to control the right-hand side through Tolsa's theorem.

We like to keep some choice on which good plane $P = P(x,r)$ to use, because some
different constraints may show up. 

\smallskip\bp
This statement looks like Lemma 6.57 in \cite{DFM3}, but since the notation may be confusing we
give a proof here. This will allow us to think at the same time
about a similar control on the difference between the quantities $\nabla_z D$ and 
$C_\alpha \lambda(x,r)^{-1/\alpha} \nabla_z \dist(z,P)$, which of course
is $C_\alpha \lambda(x,r)^{-1/\alpha}$ times the unit vector that points in the
direction opposite (and orthogonal) to $P$.

Notice that 
with our assumption \eqref{727},
both $D(z)$ and $\dist(z,P)$ are both of the order of $r$,
and $\lambda(x,r)$ is bounded from above and below, we may instead check that
\begin{equation} \label{2.17}
\big|D(z)^{-\alpha}\dist(z,P)^{\alpha} - C'_\alpha \lambda(x,r) \big| 
\leq C \sum_{l \geq 1} 2^{-\alpha l}\alpha(x,2^l r),
\end{equation}
where by \eqref{eq9.40} 
\begin{equation} \label{2.18}
D(z)^{-\alpha} = \int_E |z-y|^{-d-\alpha} d\sigma(y)
\end{equation}
which is easier to compute. And in the gradient variant, we would compare the gradient of
$D^{-\alpha}$ to $C'_\alpha \lambda(x,r)$ times the gradient of $\dist(z,P(z))^{-\alpha}$.

When we say that $P = P(x,r)$ that almost minimizes in the definition of $\alpha(x,16r)$,
we mean that there is a flat measure $\mu_0$ on $P$ such that, say, 
\begin{equation} \label{731}
\dist_{x,16r}(\sigma, \mu_0) \leq 2 \alpha(x,16r).
\end{equation}

We proceed as in  Lemma~6.57 in \cite{DFM3}, and cut $D(z)^{-\alpha}$ into pieces
\begin{equation} \label{732}
I_k = \int |z-y|^{-d-\alpha} \theta_k(y) d\sigma(y),
\end{equation}
where the $\theta_k$ form a smooth partition of $1$ such that $\theta_k$ is supported
in the annulus $A_k = B(x, 2^{k+4}r) \sm B(x, 2^{k+2}r)$ 
(but just $A_k = B(x, 16r)$ for $k=0$). We also set
\begin{equation} \label{733}
I'_k = \int |z-y|^{-d-\alpha} \theta_k(y) d\mu_0(y),
\end{equation}
Next write $\mu_0 = \lambda \H^d_{\vert P}$, and observe that
\begin{equation} \label{734}
\sum_k I'_k = \lambda \int_P |z-y|^{-d-\alpha} d\H^d(y) 
= C''_\alpha \lambda \dist(z,P)^{-\alpha},
\end{equation}
by rotation and dilation invariance.
So we want to estimate $\sum |I_k-I'_k|$. Also, for $k$ large,
$\alpha(x,16r)$ does not control the difference between $\sigma$ and $\mu_0$,
so we will need a flat measure $\mu_k = \lambda_k \H^d_{\vert P_k}$ that nearly
minimizes in the definition of $\alpha(x, 2^{k+4}r)$, as in \eqref{731} but at a larger scale;
we also set 
\begin{equation} \label{735}
I''_k = \int |z-y|^{-d-\alpha} \theta_k(y) d\mu_k(y) 
= \int f_k(y) d\mu_k(y),
\end{equation}
with $f_k(y) = |z-y|^{-d-\alpha} \theta_k(y)$. Obviously we want
to use the definition of $\dist_{x,2^{k+4}r}(\sigma, \mu_k)$ to the function $f_k$.
Notice that $f_k$ is supported in $B_k = B(x, 2^{k+4}r)$, but the reader may be afraid that
it is not smooth near $z$. 

When $k \geq 1$, we know that $z\in B(x,2r)$ and $y \in A_k$, so 
$|z-y| \geq 2^{k+1}r$, $f_k$ is Lipschitz with a constant $C (2^k r)^{-(d+\alpha+1)}$,
and \eqref{a6.2} yields
\begin{equation} \label{736}
|I''_k - I_k| \leq C ||f_k||_{lip} \, (2^{k+4}r)^{d+1} \dist_{x,2^{k+4}r}(\sigma, \mu_k) 
\leq C  (2^k r)^{-\alpha} \alpha(x,2^{k+4}r).
\end{equation}
When $k=0$, the function $f_k$ as it is defined has a singularity at $z$, 
but our assumption \eqref{727} says that it lies at distance at least
$10^{-2} r$ from both $E$ and $P$. So we may modify $\theta_0$,
so that $\theta_0$ and $f_0$ take the same values as before
on $E$ and $P$, but now $f_k$ is smooth, with $||f_0||_{lip} \leq C r^{-(d+\alpha+1)}$;
then \eqref{736} is also valid with $k=0$.

Note that if we wanted to estimate a derivative of order $m$ of $D_{\sigma,\alpha}$
we could just apply the same argument, with a function $f_k$ coming from
a derivative of $|z-y|^{-d-\alpha}$, with the effect of merely adding $C 2^{-km} r^{-m}$ in the the right-hand side of \eqref{736}. 
The same remark will apply to the computations and estimates that follow.

Next we estimate $|I''_k - I'_k|$, where we go from $\mu_0$ to $\mu_k$;
we write this as a telescopic sum, i.e., say that 
$|I''_k - I'_k| \leq \sum_{1 \leq j \leq k} \delta_{j,k}$, where 
\begin{equation} \label{737}
\delta_{j,k} = \int f_k(y) (d\mu_j(y)-d\mu_{j-1}(y)).
\end{equation}
The difference between $\mu_j$ and $\mu_{j-1}$ is controlled by
$\alpha(x, 2^{j+4}r) + \alpha(x, 2^{j+3}r)$ (compare both measures to
$\sigma$ and use the triangle inequality in \eqref{a6.2}). Since we are talking about flat measures here, this has two contributions on $\delta_j$. The first one is from the difference of densities
$|\lambda_j - \lambda_{j-1}| \leq C \alpha(x, 2^{j+4}r) + \alpha(x, 2^{j+3}r)$,
which we need to multiply by $C (2^kr)^d ||f_k(y)||_\infty \leq C (2^kr)^{-\alpha}$.
The second one is from the distance between the planes in the region $A_k$,
which is less than $C(\alpha(x, 2^{j+4}r) + \alpha(x, 2^{j+3}r)) (2^k r)$,
which we need to multiply by $C (2^kr)^d ||f_k(y)||_{lip} \leq C (2^kr)^{-1-\alpha}$.
We sum and get that 
\begin{equation} \label{738}
\delta_{j,k} \leq C (2^kr)^{-\alpha} (\alpha(x, 2^{j+4}r) + \alpha(x, 2^{j+3}r)),
\end{equation}
and then
\begin{equation} \label{739}
\sum_{k \geq 1} |I''_k - I'_k| \leq 
\sum_{k \geq 1} \sum_{1 \leq j \leq k} \delta_{j,k}
\leq \sum_{j \geq 1} (\alpha(x, 2^{j+4}r) + \alpha(x, 2^{j+3}r)) (2^jr)^{-\alpha}.
\end{equation}

Finally, we need to evaluate $\lambda - \lambda(x,r)$. Let us compute 
$I = \int \eta_r(y-\Phi_r(x)) d\mu_0$, where $\eta_r(y-\Phi_r(x))$ is the same
function that was used in the definition \eqref{2.3} of $\lambda(x,r)$.

We finally evaluate $\lambda$, by computing 
$I = \int \eta_r(y-\Phi_r(x)) d\mu_0(y)$, where $\eta_r(y-\Phi_r(x))$ is the same
function that was used in the definition \eqref{2.3} of $\lambda(x,r)$.
This way, replacing $\mu_0$ with $\sigma$ in $I$ would yield $\lambda(x,r)$,
and so
\begin{equation} \label{740}
|I - \lambda(x,r)| \leq C r^{d+1} ||\eta_r(\cdot-\Phi_r(x))||_{lip} \, \alpha(x,16r)
\leq C  \alpha(x,16r)
\end{equation}
by \eqref{a6.2} and the definition of $\mu_0$.
If $\Phi_r(x)$ were luckily lying on $P$, we would get $I = \lambda$ immediately,
because $\eta$ is radial and $\int \eta_r =1$ on $\R^d$; this is not necessarily true,
but we will check in a moment that
\begin{equation} \label{2.21}
\dist(\Phi_r(x),P) \leq C r \alpha(x,16r),
\end{equation}
and then it will follow, by the usual argument where we estimate the derivative of
$\xi \to \int_{P(x,r)} \eta_r(z-\xi) d\mu_0(z)$ along a segment from $\Phi_r(x)$ to $P$, that 
\begin{equation} \label{742}
|\lambda - \lambda(x,r)| \leq  C \alpha(x,16r).
\end{equation}
Incidentally, this is the reason why we decided to use $\lambda(x,r)$ rather
than $\lambda_0$: it could happen that $\dist(x,P)$ is much larger than $C r \alpha(x,16r)$.
To check \eqref{2.21} we return to the definition of $\Phi_r(x)$ by \eqref{2.2},
project on the $(n-d)$-space orthogonal the direction of $P(x,r)$, and then use the triangle inequality to find that
\begin{equation} \label{2.22}
\begin{split}
\dist(\Phi_r(x),P)  
&\leq \lambda_0(x,r)^{-1} \int_{E \cap B(x,r)} \dist(y,P(x,r) \, \eta_r(y-x) d\sigma(y)
\\ &\leq C \lambda_0(x,r)^{-1} r  \alpha(x,16r)
\leq C r  \alpha(x,16r),
\end{split}
\end{equation}
where the last inequalities come again from \eqref{a6.2} and the definition of $\mu_0$,
because the same integral, but against $d\mu_0$, would give $0$ because we would 
integrate on $P(x,r)$; \eqref{2.21} follows.

We may now summarize. We have seen that $D(z)^{-\alpha} = \sum_k I_k$
is quite close to $\sum_k I''_k$, by \eqref{736}, and then to $\sum_k I'_k$,
by \eqref{739}; then by \eqref{734} $\sum_k I'_k = C''_\alpha \lambda \dist(z,P)^{-\alpha}$
Thus by \eqref{742}
\begin{equation} \label{744}
\big| D(z)^{-\alpha} - C''_\alpha \lambda(x,r) \dist(z,P)^{-\alpha} \big|
\leq C r^{-\alpha} \sum_k 2^{-k \alpha} \alpha(x,2^{k+4}r)
\end{equation}
which is the same as \eqref{2.17} and implies \eqref{2.16}. 
Lemma \ref{l72} follows.
\ep

\subsection{ The distance $D_\Sigma(g(x,t))$}

In what follows, we return to the construction of a bilipschitz change
of variable associated to a stopping time region $\cF$, and give a good evaluation
of the distance $D_\Sigma(g(x,t))$ 
associated to $\Sigma = g(\R^d)$,
first compared to the distance to a good plane. 

We shall use the same notation as in the first part, and in particular work
on $\R^n$, except that we use the coordinates $y\in \R^d$ and $t \in \R^{n-d}$
to avoid some confusion with the previous subsection.
In the present subsection, we use the distance $D_\Sigma$ associated with the 
surface $\Sigma = g(\R^d)$, where we put any Ahlfors-regular measure 
$\sigma$, and we relate this to our change of variable. 
We no longer mention the exponent $\alpha$ in the notation, 
both because it is fixed and we want to avoid extra confusion with the Tolsa numbers.

\begin{pro} \label{Pe1}
The function $\Phi$ defined by 
\begin{equation} \label{e28}
\Phi(y,t) = \Big|\frac{D_\Sigma(g(y,t))}{|t|} - C_\alpha \lambda_\sigma(f(y),|t|)^{-1/\alpha} \Big| 
\end{equation}
satisfies a Carleson condition on $\Omega_0$. 
\end{pro}

As the reader guessed,
$\lambda_\sigma$ is the same function $\lambda$ as above, but associated
to the measure $\sigma$ and the set $\Sigma$. The constant $C_\alpha$
is the same as above.

Later in the subsection, we will manage to apply \eqref{2.16}, but for the moment
we first estimate the distance from $g(y,t)$ to some other plane that we define now.
Let $(y,t) \in \Omega_0$ be given, and as always set $r = |t|$.
Also set $x = f(y)$ and $z = g(y,t)$. 
Choose $k = k(t)$, as we did above, to be the largest integer such that 
$\rho_k(t) \neq 0$ (see near the definition \eqref{c22} of $g$). 
Let as before $T_k(y)$ denote the tangent plane to $\Sigma_k$ at $f_k(y)$, 
and recall from \eqref{c22} that 
\begin{equation} \label{e29}
z = g(y,t) = \sum_{j\geq 0} \rho_k(t) \Big\{ f_j(y) + R_j(y)\cdot t \Big\}.
\end{equation}
In this  sum there are at most $3$ terms, corresponding to $j = k, k-1, k-2$
(when they are nonnegative), and we will see that these terms are almost the same.
Set $z' = f_k(y) + R_k(y)\cdot t$; notice that since 
$z' - f_k(y) = R_k(y)\cdot t$ is orthogonal to
$T_k(y)$ 
(because $R_k(y)$ maps (the orthogonal complement of) $\R^d$ to 
(the orthogonal complement of) $T_k(y)$,
\begin{equation} \label{e30}
\dist(z',T_k(y)) = |z' - f_k(y)| = |t|.
\end{equation}
We want to compare $D_\Sigma (g(x,t)) = D_\Sigma(z)$
to $\dist(z', T_k(y))$. 
Some error terms will come from $ |z-z'|$, but observe that
\begin{equation} \label{e31}
\psi_1(y,t) := |t|^{-1} |z-z'| \leq |t|^{-1}
\sum_{k-2 \leq j \leq k-1} |f_j(y) - f_k(y)| + |t| \, |R_j(y)-R_k(y)|; 
\end{equation}
(by \eqref{e29}). Notice that $\psi_1$ satisfies a Carleson condition by our treatment of 
$D_1$ from \eqref{d7}.

Next we want to use \eqref{e30} to estimate $\dist(z',P(x,r))$ where $P(x,r)$
is a good plane to apply \eqref{2.16}. More precisely, since we want to apply \eqref{2.16}
to various points, we choose for each $x\in \Sigma$ and $r > 0$ a nearly optimal 
flat measure $\mu_{x,r}$ for $\alpha(x,16r)$ (where the $\alpha$-numbers 
are associated to $\Sigma$ and $\sigma$), and then let $P(x,r)$ be the support of $\mu_{x,r}$.

We return to our initial pair $(y,t)$, and try to estimate the distance between $P(x,r)$
and $T_k(y)$. This will take some time, but we shall remember that the main property of 
$P(x,r)$ in this respect is that
\begin{equation} \label{e32}
\int_{\Sigma \cap B(x,8r)} \dist(w,P(x,r)) \, d\sigma(w) \leq C r^{d+1} \alpha(x,16r),
\end{equation}
which as usual we obtain by testing the product of $\dist(w,P(x,r))$ by a bump
function against the difference $\sigma - \mu_{x,r}$. Thus it makes sense to estimate
the distance from points of $\Sigma \cap B(x,8r)$ to $T_k(y)$ too.

We start with the distance from points of $\Sigma_k$ to $T_k(y)$. 
Let $L$ be a bound for the biLipschitz constant for $f$ and the $f_k$.
Such a uniform bound for the $f_k$ comes from the proof of Lemma~\ref{Lb1},
but if the reader does not want to believe this, there is an easy fix explained below.
Set 
\begin{equation} \label{e33}
\psi_2(y,t) = r^{-1} \sup_{w \in \Sigma_k \cap B(x, (10L)^{-2} r)} \dist(w,T_k(y));
\end{equation}
we want to show that this is a Carleson function. 
For $t$ so large that $k=k(t) = 0$, $\Sigma_k = T_k(y) = P_0$ and so $\psi_2(y,t) = 0$.
Hence we can restrict our attention to the pairs $(y,t)$ such that $k = k(t) \geq 1$.

We start with the simpler function
\begin{equation} \label{e34}
\psi_3(y,t) = \sup_{y'\in P_0 \cap B(y,((2L)^{-1}r_k)} r_k |DR_{k}(y')|,
\end{equation}
where in fact we restrict to $t$ such that $k(t) \geq 1$ (otherwise set
$\psi_3(y,t)=0)$.

To estimate $\psi_3(y,t)$, let $y'\in P_0 \cap B(y,(2L)^{-1}r_k)$, 
and first assume that $f_k(y') \in V_k^8$. Then \eqref{c14} says that 
$r_k |DR_{k}(y')| \leq C \wh\varepsilon_k(y')$, 
where $\wh\varepsilon_k(y') = \sum_{l \leq k} 2^{l-k} \varepsilon''_l(f(y')$.
But recall that when we chose $f(y')$ to evaluate the $\varepsilon''_l$,
we could in fact have chosen any point $w$ such that $|w-f(y')| \leq r_k/2$,
and in particular, since $f$ is $M$-biLipschitz and $|y'-y| < (2L)^{-1}r_k$, 
$w = f(y)$. Thus $|DR_{k}(y')| \leq C \wh\varepsilon_k(y)$ in this case.
Notice also that if $f_k(y) \in V_k^7$ and since $f_k$ is biLipschitz,
this holds for all $y'\in P_0 \cap B(y,(2L)^{-1}r_k)$. If you do not trust this,
use \eqref{b36} to go through $|f(y)-f(y')|$.

A second case is when $f_k(y) \in V_k^{11} \sm V_k^7$; this is a little larger
than the usual dying zone, but Lemma \ref{Lb2} still says that for a given $y$ this
happens for at most one $k$. Then we use \eqref{c9}, which is valid
everywhere, to get that $r_k |DR_{k}(y')| \leq C\varepsilon$ on $P_0 \cap B(y,(2L)^{-1}r_k)$.

When $f_k(y) \in \Sigma_k \sm V_k^{11}$, 
we return to the largest $l$ such that
$f_l(V_k^{11})$, find that $R_k(y') = R_{l+1}(y')$, use the estimate above, and
find that $r_k |DR_{k}(y')| \leq C 10^{l-k}\varepsilon$ on $P_0 \cap B(y,(2L)^{-1}r_k)$.
Now we can follow our estimates for $D_4$ (see near \eqref{d23}-\eqref{d24}) and find that
$\psi_3$ satisfies the Carleson condition. 

Next we use $\psi_3$ to control $\psi_2$. Let $w \in \Sigma_k \cap B(x, (10L)^{-2} r)$
be given. A way to find out where $w$ lies is to return to $u\in \Sigma_0$ such that
$f_k(u) = w$, take the line segment $[y,u]$, and follow its image by $f_k$. 
Recall that $x = f(y)$ and $|f(u)-w| = |f(u)-f_k(u)| \leq C r_k \varepsilon$ by \eqref{b36}, 
so $|y-u| \leq (90L)^{-1} r \leq (4L)^{-1} r_k$ because $r \leq 20 r_k$ by \eqref{dd3}.
Then by \eqref{e33}  
\begin{equation} \label{e35}
|R_k(s)-R_k(y)| \leq |s-y| \, r_k^{-1} \psi_3(y,t) \ \ \text{ for } s\in [y,u].
\end{equation}
But $R_k(y)$ maps $\R^d$ to the vector space parallel to $T_k(y)$,
so the derivative in $s \in [y,u]$ of $\dist(f_k(s),T_k(y))$ is at most
$C |R_k(s)-R_k(y)| \leq C \psi_3(y,t)$. Of course this distance is null for $f_k(y)$,
hence it is at most $C r_k \psi_3(y,t)$ at the end of the path, for $w = f_k(u)$.
That is, $\dist(w,T_k(y)) \leq C r_k \psi_3(y,t)$, and this proves that $\psi_2 \leq C \psi_3$,
hence $\psi_2$ satisfies a Carleson measure estimate.

\smallskip
This is not over yet; now want to control the average distance
from $\Sigma$ to $\Sigma_k$, i.e.,
\begin{equation} \label{e36}
\psi_4(y,t) = r^{-d-1} \int_{\Sigma \cap B(x, (20L)^{-1} r)} \dist(w,\Sigma_k) d\sigma(w),
\end{equation}
and show that
\begin{equation} \label{e37}
\psi_4 \ \text{ satisfies a Carleson estimate.} 
\end{equation}
We start when $r \geq 10^{-2}$ use the fact that for $w\in \Sigma$,
we can write $w = f(u)$ for some $u\in P_0$, and then
$\dist(w,\Sigma_k) \leq |w-f_k(u)| = |f(u)-f_k(u)| 
\leq C \varepsilon$
by \eqref{b36}, so that
\begin{equation} \label{e38}
\psi_4(y,t) \leq C \varepsilon  
 \ \text{ for } |t| \geq 10^{-2}.
\end{equation}
It is easy to see that \eqref{e38} gives a bounded contribution to the Carleson norm
of $\psi_4(x,t)^2 \frac{dx dt}{|t|^{n-d}}$.
Otherwise, when $|t|\leq 10^{-2}$, $k \geq 1$ and, since $r_k \leq r \leq 20r_k$ by \eqref{dd3},
$$
\begin{aligned}
\psi_4(y,t) &\leq C r_k^{-d-1} \int_{\Sigma \cap B(x, L^{-1} r_k)} \dist(w,\Sigma_k) d\sigma(w)
\leq C r_k^{-d-1} \int_{P_0 \cap B(y, r_k)} \dist(f(u),\Sigma_k) du
\cr&\leq C r_k^{-d-1} \int_{P_0 \cap B(y, r_k)} |f(u)-f_k(u)| du
\leq C r_k^{-d-1} \int_{P_0 \cap B(y, r_k)} \sum_{l \geq k}  |f_{l+1}(u)-f_l(u)| du.
\end{aligned}
$$
We take our earlier estimate for $|f_{l+1}(u)-f_l(u)|$, which we did when we estimated
$D_{11}$ in \eqref{dd18}, and which writes
\begin{equation} \label{e40}
|f_{l+1}(u)-f_l(u)| \leq C r_l (\varepsilon''_l(u) + \varepsilon \delta_l(u))
\end{equation}
where $\delta_k(u) = 1$ when $f_l(u) \in V_l^{10} \sm V_l^8$, and $\delta_k(u) = 0$
otherwise. Thus by Cauchy-Schwarz
\begin{equation} \label{e41}
\psi_4(y,t)^2 \leq C r_k^{-d} \sum_{l \geq k} 10^{l-k}
\int_{P_0 \cap B(y, r_k)} \varepsilon''_l(u)^2 + \varepsilon^2 \delta_l(u)^2 du,
\end{equation}
where $10^{l-k} = r_k^{-1} r_l$. Set $\psi'_4(y,t) = \psi_4(t,y) \1_{k(t) \geq 1}$
(the piece that we estimate now). We integrate \eqref{e41} on a Carleson box 
$B(X,R) \subset \R^n$ centered at $X = (x_0,0) \in P_0$, and get 
$$
\int_{B(X,R)} \psi'_4(y,t)^2 \frac{dydt}{|t|^{n-d}} 
\leq C \sum_{k \geq 1} \sum_{l \geq k} 10^{l-k} \int_{(y,t) \in B(X,R); k(t) = k} 
\int_{P_0 \cap B(y, r_k)} \varepsilon''_l(u)^2 + \varepsilon^2 \delta_l(u)^2 
\frac{dudydt}{r_k^d |t|^{n-d}}.
$$
 Given $y$, $u$, $l \geq 1$, and $k \in [1,l]$, we integrate in the region 
 where $r_k \leq |t| \leq 20r_k$,
 where
 $\int_t \frac{dt}{|t|^{n-d}} \leq C$.
Then we sum over $y \in B(u, r_k)$ and make the $r_k^d$ disappear.
Then we sum the geometric series in $k$, and are left with
\begin{equation} \label{e42}
\int_{B(X,R)} \psi'_4(y,t)^2 \frac{dydt}{|t|^{n-d}} 
\leq C \sum_l \int_{P_0 \cap B(x_0, 2R)} \varepsilon''_l(u)^2 + \varepsilon^2 \delta_l(u)^2.
\end{equation}
where we used the fact that $|u-x_0| \leq |u-y|+|y-x_0| \leq r_k + R$ and $r_k \leq |t| \leq R$.
Notice that $r_l \leq r_k \leq |t| \leq R$ in the sum above; we use the Carleson measure estimates
proved in Section \ref{Sd} and get less than $CR^d$. This completes our proof of \eqref{e37}.

\ms 
We are now ready to compare $P(x,r)$ 
(from \eqref{e32})
and $T_k(y)$, and we start with the most
interesting case when $r \leq 100$ (so that $r_k \sim r$).
Since \eqref{e32} and \eqref{e36} are merely averages, we start with a Chebyshev 
argument to select good points $\xi_j$ of $\Sigma \cap B(x,r)$. 
We assume that $\psi_2(y,t)$, $\psi_4(y,t)$ and $\alpha(x,16r)$ are small, 
otherwise we will be happy with a trivial estimate. 

We first choose points near which we want to select these points.
Lemma 6.2 in \cite{DT} 
says that in $B(x,19r_k)$, $\Sigma_k$ coincides with
the graph over a plane $P$ of some $C \varepsilon$-Lipschitz function.
Let $P'$ be the vector $d$-plane parallel to $P$, and choose an orthonormal basis
$e_1, \ldots, e_d$ of $P'$. Then set 
$w_j = f_k(y) + (20L)^{-2} e_j r_k$ 
for $j \geq 1$ and $w_0 = f_k(y)$. By the Lipschitz graph description, 
we can find points $w'_j \in \Sigma_k$ 
such that $|w'_j-w_j| \leq (200L)^{-2} r_k$. 
Then by \eqref{e32}, \eqref{e36}, the Ahlfors regularity of $\Sigma$, and Chebyshev, we can find 
points $\xi_j \in \Sigma_k$, such that $|\xi_j-w_j| \leq (200L)^{-2} r_k$, and for which
\begin{equation} \label{e43}
\dist(\xi_j,P(x,r)) \leq C \alpha(x,16r) r_k \ \text{ and } \ 
\dist(\xi_j,\Sigma_k) \leq C \psi_4(y,t).
\end{equation}
By the definition \eqref{e33} of $\psi_2$ (and if $\psi_4(y,t)$ is small enough to guarantee
that $\dist(\xi_j,\Sigma_k) \leq (4L)^{-1}$), we even get that
\begin{equation} \label{e44}
\dist(\xi_j,T_k(y)) 
\leq C \psi_4(y,t) + \psi_2(y,t).
\end{equation}
Thus we manage to find $d+1$ points $\xi_i$ of $\Sigma$, that are sufficiently far from
each other, and that all lie very close to both $P(x,r)$ and $T_k(y)$.
With a little bit of geometry, we get that
\begin{equation} \label{e45}
d_{x,10r_k}(P(x,r),T_k(y)) \leq C \alpha(x,16r) + C \psi_2(y,t) + C\psi_4(y,t);
\end{equation}
see the discussion below \eqref{b8} and Lemma 12.7 on page 74 
of \cite{DT}.
Of course this estimate is still valid when the right-hand side of \eqref{e45} is large,
but it is useless.

Also, we forgot the case when $r \geq 100$, but then $k=0$, $T_k(y) = P_0$, 
the construction of $\Sigma_k$ gives $\dist(w,P_0) \leq C \varepsilon$ for every 
$w \in \Sigma_k$, and the same argument as above yields
\begin{equation} \label{e46}
d_{x,10r}(P(x,r),T_k(y)) = d_{x,10r}(P(x,r),P_0) \leq C \alpha(x,16r) + C \varepsilon r^{-1}.
\end{equation}

\ms
Finally we return to the point $z'$ of \eqref{e30}. Notice that $z'\in B(x,2r)$ because 
$|z'-f_k(y)| = |t| = r$ 
(by \eqref{e30}) and  
$|f_k(y) - x| = |f_k(y) - f(y)| \leq C \varepsilon r_k \leq C \varepsilon r$ 
(by \eqref{b36}).  In addition,
\begin{equation} \label{e47}
r^{-1} |\dist(z',P(x,r)) - r| \leq C \alpha(x,16r) r + C \psi_2(y,t) + C \psi_4(y,t) r
\end{equation}
when $r \leq 100$ 
(by  \eqref{e30} and \eqref{e45}), 
and
\begin{equation} \label{e48}
r^{-1} |\dist(z',P(x,r)) - r| \leq C \alpha(x,16r) r + C \varepsilon r^{-1}
\end{equation}
otherwise (by \eqref{e46}). 
Let us first assume for the moment that these numbers are small, and 
try to apply Lemma \ref{l72} to $z'$ and $P(x,r)$. 
Notice that $\dist(z', P(x,r)) \geq r/2$ directly by \eqref{e48}, 
but we also need to show that $\dist(z', P(x,r)) \geq 10^{-2} r$.
Recall that $\int_{\Sigma \cap B(x,8r)} \dist(w,P(x,r)) \, d\sigma(w) 
\leq C r^{d+1} \alpha(x,16r)$, by \eqref{e32}; it then follows from the Ahlfors regularity
of $\Sigma$ that $\dist(w,P(x,r)) \leq 10^{-1}r$ for $w \in E \cap B(x,7r)$,
and since $z'\in B(x,2r)$, that $\dist(z', P(x,r)) \geq 10^{-2} r$, as needed for
\eqref{727}. So \eqref{2.16} holds, hence
\begin{equation} \label{e49}
|D_\Sigma(z') - C_\alpha \lambda_\sigma(x,r)^{-1/\alpha} \dist(z,P(x,r))|
\leq C r \sum_{l \geq 4}
2^{-\alpha l}\alpha(x,2^l r)
\end{equation}
(where we also used the fact that $\dist(z,P(x,r)) \leq 3r$
to multiply the estimate). Recall from \eqref{e31} that $|z-z'| = |t| \psi_1(y,t) = r\psi_1(y,t)$, 
and since it is easy to check that $D_\Sigma$ is Lipschitz, we also have
\begin{equation} \label{e50}
|D_\Sigma(z') - |D_\Sigma(z)| \leq C r \psi_1(y,t).
\end{equation}
We add \eqref{e50}, \eqref{e49}, and \eqref{e48} or \eqref{e47} and get a 
good control on $r^{-1} |D_\Sigma(z) - C_\alpha \lambda_\sigma(x,r)^{-1/\alpha} |t||$, 
which is the same as $\Phi(y,t)$ in \eqref{e28}. 

When the controlling numbers in \eqref{e48} are large, we just say that $\Phi(y,t) \leq C$.

At this point, we have a good control of $\Phi(y,t)$ by various quantities, 
which are functions of $(y,t)$, and we just need to check that they
satisfy Carleson measure estimates on $\Omega_0$.

For the functions $\psi_i$, this was proved along the way. 
For the $\alpha$-function, it is a function of $(x,r)$ that satisfies a Carleson 
measure estimate in $\Sigma \times (0,+\infty)$, by the theorem of Tolsa \cite{To} 
(also see Lemma 5.89 in \cite{DFM3} for the control of the geometric series), and
it is easy to see that when we compose it with the mapping $(y,t) \to (x,r) = (f(x),|t|)$, 
we get a function of $(x,t)$ that satisfies a Carleson estimate. 
We are left with the last term  $C \varepsilon r^{-1} \1_{r \geq 100}$, from \eqref{e48}, 
which also satisfies a Carleson measure estimate by direct computation. This 
completes our proof of Proposition \ref{Pe1}.
\qed

\section{Surgery with $D_{\Sigma}$ and $D_E$}
\label{Sf}

In the previous section we managed to control reasonably well the effect
of our change of variable $g$, provided that we consider the distance 
$D_\Sigma = D_{\Sigma,\sigma}$ associated to an (in fact, any) 
Ahlfors regular measure $\sigma$ on $\Sigma$. 

In the larger picture, we started from a set $E$, with its own Ahlfors regular measure
(we shall now call it $\mu$), and we would like to use the corresponding distance
function $D_E = D_{E,\mu}$. Notice that both measures also
depend on $\alpha$, but we shall not mention this in the notation.

It does not make sense to compare our two measures in the places where
$E$ and $\Sigma$ have nothing to do with each other, so we will only compare
them in the same region 
\begin{equation} \label{f1}
\Omega_\cF = \bigcup_{Q \in \cF} W(Q),
\end{equation}
where for each $Q \in \cF$, $W(Q)$ is the Whitney box 
defined by
\begin{equation} \label{f2}
W(Q) = \big\{ x\in B(x_Q,M_0 l(Q)) \, ; \, \dist(x, E) \geq M_0^{-1} l(Q) \big\}.
\end{equation}
As the reader may have guessed, the precise shape of $\Omega_\cF$ does not matter
so much, but we probably don't want it to be too 
small because this is the region where we can play. We state the main result of this section, and then discuss.

\begin{pro}\label{Pf1}
Suppose $E$ is  
uniformly rectifiable, $\mu$ is an Ahlfors regular measure on $E$,
$\cF$ satisfies the conditions of Section \ref{S2} (with $M$ large enough and $\varepsilon_1$
small enough), and $\Sigma = f(P_0)$ denote the surface constructed above. 
Then there is an Ahlfors regular measure $\sigma$ on $\Sigma$ such that
\begin{equation} \label{f3}
\Psi := \1_{\Omega_\cF}  \,\Big| \frac{D_{E,\mu}}{D_{\Sigma,\sigma}}  -1 \Big|
\ \text{ satisfies a Carleson measure condition on $\R^n \sm \Sigma$.}  
\end{equation}
\end{pro}

Here it is more convenient (or just safer) to let $M_0$ be as large as possible,
then choose $M$ and $\varepsilon_1$, and do the stopping time construction
accordingly. We may pay a huge price (depending on $M_0$), but this is more 
transparent. 

We decided to require a Carleson measure estimate relative to $\Sigma$
because $\R^n \sm \Sigma$ is the place the where distance function $\wh D$ below will live, 
and this is also the region where we hope to use our change of variable to control
operators. 
We could 
equally prove 
a Carleson measure estimate relative to $E$, in fact with the same proof; 
see the remark below \eqref{f38}.

The goal of this is to control a degenerate elliptic operator $L = L_E$ associated to 
$D_{E,\mu}$ on $\Omega_E = \R^n \sm E$, and we hope to compare it to an
operator $\wh L$ on $\Omega_\Sigma = \R^n \sm \Sigma$, but we would like to keep
the same formula on a set which is as large as possible, so we let $\wh L$ be associated to 
the distance function $\wh D$ defined by
 \begin{equation} \label{f4}
\wh D(z) = D_{E,\mu}(z) \ \text{ for } z\in \Omega_\cF
\end{equation}
and 
\begin{equation} \label{f5}
\wh D(z) = D_{\Sigma,\sigma}(z) \ \text{ for } z\in \R^n \sm [\Sigma \cup\Omega_\cF].
\end{equation}
We do not fear a discontinuity between the two regions; our elliptic conditions allow this.
But we will
need to make sure that 
\begin{equation} \label{f6}
C^{-1} \dist(z,\Sigma) \leq \wh D(z) \leq C \dist(z,\Sigma)
\end{equation}
(we will do this after \eqref{ff10}), to make sure that $\wh L$ lies in the class of
acceptable operators studied in \cite{DFM2}. This estimate is also reassuring,
because it says that $D_{E,\mu}$ and $D_{\Sigma,\sigma}$ are equivalent on $\Omega_\cF$,
which implies in particular that $\dist(z,\Sigma) \geq C^{-1} l(Q)$ on $W(Q)$.

Then 
Proposition \ref{Pf1} will allow us to prove, via the change of variable
$g$, that the elliptic operator associated to $\wh L$ has an absolutely continuous 
harmonic measure, because of a Carleson control on $|t|^{-1} \wh D \circ g$
that comes from Proposition \ref{Pf1} and Proposition \ref{Pe1}. So what will be
left to do is use the fact that $\wh D = D_E$ on the hopefully sufficiently large
region $\Omega_\cF$, to get some control on $L_E$ itself.

We start our proof with some basic geometric information about $E$ and $\Sigma$. 

\begin{lem} \label{Lf2}
Set $M_1 = 10^{-2} M$. Then for each $Q \in \cF$,
\begin{equation} \label{f13}
M_1 d_{x_Q, M_1l(Q)}(E, \Sigma) \leq C \varepsilon.
\end{equation}
\end{lem}

See \eqref{a8} for the definition of the normalized local Hausdorff distance $d$.
As usual, this is true if we assume that $M$ is large enough, and $\varepsilon_1$
is small enough, depending on $M$. Also we added $M_1$ on the left-hand side
just not to lose an additional $M$ uselessly, but this does not matter because we 
always choose $\varepsilon$ and $\varepsilon_1$ last.
Let $Q\in \cF$ be given, and let $P(Q)$ be as 
in \eqref{b1} and \eqref{b2}. That is,
\begin{equation} \label{f14}
d_{x_Q, M l(Q)}(E,P(Q)) \leq 2\varepsilon_1.
\end{equation}
We now want to prove that
\begin{equation} \label{f15}
M_1 d_{x_Q, 2M_1 l(Q)}(\Sigma,P(Q)) \leq C\varepsilon,
\end{equation}
and the lemma will follow. Set  $B(Q) = B(x_Q, 3M_1l(Q))$ and let
$p \in P(Q) \cap B(Q)$ be given. First use \eqref{f14} to find
$x\in E$ such that $|x-p| \leq C M l(Q) \varepsilon_1$. Set $k = k(Q)$
(the generation of $Q$ and observe that $x$ lies in the set $E(k)$ of \eqref{b5}.
Hence (by the line below \eqref{b5}) we can find $j\in J_k$ such that $|x-x_{j,k}| \leq 2r_k$.
This is good, because then Proposition 5.1 in \cite{DT}  
gives a good description of
$\Sigma_k \cap B(x_{j,k},49r_k)$ as a piece of a $C\varepsilon$-Lipschitz graph over
$P_{j,k}$ that passes within $C\varepsilon r_k$ from $x_{j,k}$. 
Recall that we even managed to pick planes $P_{j,k}$ that contain $x_{j,k}$,
but if $P_{j,k}$ was only $C\varepsilon r_k$-close, what we are going to say would work
too. The small Lipschitz graph description implies that every point of 
$P_{j,k} \cap B(x_{j,k},48r_k)$ lies within $C \varepsilon r_k$ of $\Sigma_k$.
By \eqref{b36}, it also lies within $\varepsilon r_k$ of $\Sigma$ (recall that
$\Sigma_k = f_k(P_0)$ and $\Sigma = f(P_0)$). Now $P_{j,k}$ was chosen
to be equal to $P(Q_{j,k})$ for some $Q_{j,k} \in \cF(k)$ such that 
$\dist(x_{j,k},Q_{j,k}) \leq \frac{Mr_k}{10}$ (and we even had to move $x_{j,k}$
slightly so that $P(Q_{j,k})$ goes through $x_{j,k}$, but this is not the point here),
and so $P(Q_{j,k})$ is quite close to $P(Q)$ near $x_{j,k}$ (by \eqref{b1} and \eqref{b2}
for both cubes). Consequently, $p$ lies within $C \varepsilon M r_k$ of
$P_{j,k} \cap B(x_{j,k},47r_k)$, and we find $\xi \in \Sigma$ such that 
$|\xi-p| \leq C \varepsilon M r_k$. 

Now we 
take $\xi \in \Sigma \cap B(Q)$ and try to find $p\in P(Q)$ near $\xi$.
Let $x\in P_0$ be such that $\xi = f(x)$. The easiest case is when $f_k(x) \in V_k^{10}$,
because this means that we can find $j \in J_k$ such that $|x_{j,k}-f_k(x)| \leq 10r_k$,
and we can use the same Lipschitz graph description of $\Sigma_k \cap B(x_{j,k},49r_k)$
as above. We find a point $p' \in P_{j,k} \cap B(x_{j,k},49r_k)$ such that
$|p'-f_k(x)| \leq C \varepsilon r_k$, 
and use the fact that $P(Q_{j,k})$ is quite close to $P(Q)$ near $x_{j,k}$ to find
$p\in P(Q)$ such that $|p-p'| \leq C M \varepsilon_1 r_k \leq C \varepsilon r_k$.
Then $|p-\xi| \leq |p-p'|+|p'-f_k(x)| + |f_k(x)-\xi| \leq C \varepsilon r_k$ by
\eqref{b36}. 

We are left with the case when $f_k(x) \in \Sigma_k \sm V_k^{10}$.
First assume that $f_l(x) \in V_l^{10}$ for some $l \in [0,k-1]$, and take
$l$ as large as possible. Then take $j \in J(l)$ so that $|f_l(x)-x_{j,l}| \leq 10r_l$.
Again use the good Lipschitz description of $\Sigma_l \cap B(x_{j,k},49r_l)$
provided by Proposition 5.1 in \cite{DS}, or the case $k=l$ of the description
above: there are points of $P_{j,l}$, and then points of $E$, that lie at distance
less than $C \varepsilon r_l$ from $f_k(x)$. Use this to pick $w\in E \cap B(f_k(x),r_{l+1})$.
By definition of $l$, $w\in \R^n \sm V_{l+1}^{10}$, which implies that 
$x \notin E(l+1)$. In other words, $\dist(x,R) \geq \frac{M r_l}{10}$ for
every $R \in \cF(l+1)$. Since $|w-\xi| = |w-f(x)| \leq |w-f_l(x)| + C \varepsilon r_l
\leq 2r_{l+1}$, we see that $\dist(\xi,R) \geq \frac{M r_l}{11}$ for every $R \in \cF(l+1)$.
We apply this to the ancestor of $Q$ of generation $l$, and find that
$\dist(\xi,R) \geq \frac{M r_l}{11} \geq \frac{10 M r_k}{11}$, which contradicts
our assumption that $\xi \in \Sigma \cap B(Q)$. 

If we cannot find $l < k$ such that $f_l(x) \in V_l^{10}$, then 
$f_k(x) = f_0(x) = x \in P_0 \sm V_0^{10}$, and as before 
$\dist(x,Q_0) \geq M/11$. This is also impossible because $\xi \in B(Q)$.
This completes our proof of \eqref{f15}, and the lemma follows.
\qed

\ms
A simple consequence of this is the following improvement of \eqref{f6}.
We claim that
\begin{equation} \label{ff10}
(1-CM_0\varepsilon) \dist(z,E) \leq \dist(z,\Sigma) \leq (1+CM_0\varepsilon) \dist(z,E)
\ \text{ for } z\in \Omega_\cF.
\end{equation}
Indeed let $z\in \Omega_\cF$ and let $Q\in \cF$ be such that $z\in W(Q)$.
Let $w\in E$ be such that $|w-z| = \dist(z,E)$; observe that $\dist(z,E) \leq |z-x_Q| \leq M_0l(Q)$
by \eqref{f2}, hence $|w-x_Q| \leq 2M_0l(Q) < M_1l(Q)$ if $M$ is large enough, 
Lemma \ref{Lf2} applies to $w$ and gives $\dist(w,\Sigma) \leq C \varepsilon l(Q)$,
so $\dist(z,\Sigma) \leq \dist(z,E) + C \varepsilon l(Q) \leq (1+CM_0\varepsilon) \dist(z,E)$
because $\dist(z,E) \geq M_0^{-1} l(Q)$ by \eqref{f2}. Now let $\xi \in \Sigma$ be such that
$\dist(z,\Sigma) = |z-\xi|$. If $|z-\xi| \geq M_0 l(Q)$, the first inequality in \eqref{ff10} is
trivial. Otherwise, we can apply Lemma \ref{Lf2} to $\xi$ and get that 
$\dist(\xi,E) \leq C \varepsilon l(Q)$, hence 
$\dist(z,E) \leq \dist(z,\Sigma) + C \varepsilon l(Q) \leq \dist(z,\Sigma) +CM_0\varepsilon \dist(z,E)$
and the first part of \eqref{ff10}  follows.

Notice that \eqref{f6} follows from this, since $D_E(z)$ is equivalent to $\dist(z,E)$
and $D_\Sigma(z)$ is equivalent to $\dist(z,E)$.

\ms
Our next task is to construct Whitney cubes (in fact, pseudocubes) in $E$,
which we will use to define the measure $\sigma$ on $\Sigma$ that approximates $\mu$.
For this we will use the somewhat classic distance to small cubes of $\cF$, defined by
\begin{equation} \label{f7}
d_\cF(z) = \inf_{Q \in \cF}\, (\dist(z,Q) + l(Q))
\end{equation}
for $z\in \R^n$. Notice that 
$d_\cF(z) \leq 1 + \dist(z,Q_0)$ (because we can try $Q = Q_0 \in \cF$), and $d_\cF$ is $1$-Lipschitz.
Associated to $d_\cF$ are a closed set
\begin{equation} \label{f8}
F = \big\{ z\in E \, ; \, d_\cF(z) = 0 \big\}
\end{equation}
and a decomposition of $E \sm F$ into \ub{Whitney cubes} that we describe now. 
We give ourselves a small constant $\tau \in (0,10^{-2})$, 
and we denote by $\cR$ the collection of maximal cubes $R \in \bD$ 
(for the inclusion as a first criterion, and then the smallest generation if a same set 
corresponds to cubes of different generations), with the property 
\begin{equation} \label{f9}
l(R) \leq \tau d_\cF (x_R).
\end{equation}
These cubes are disjoint (by maximality), they do not meet $F$ because
it is easy to see that for $R \in \cR$,
\begin{equation} \label{f10}
d_\cF (x) \geq (2 \tau)^{-1} l(R) \ \text{ for } x\in R,
\end{equation}
because $d_\cF$ is $1$-Lipschitz. The maximality of $R$ implies that
its parent $S$ does not satisfy \eqref{f9}, hence $10l(R) = l(S) \geq \tau d_\cF (x_S)$, 
and since 
$$
|d_\cF (x_S) - d_\cF (x_R)| \leq |x_S - x_R| \leq \diam(S) \leq 20 l(R) 
\leq 20\tau d_\cF (x_R) \leq d_\cF (x_R)/5, 
$$
we get that
\begin{equation} \label{ff15}
l(R) \geq 10^{-1} \tau d_\cF(x_S) \geq 20^{-1} \tau d_\cF (x_R) \ \text{ if } R \in \cR.
\end{equation}
We claim that
\begin{equation} \label{f11}
E \sm F \ 
\text{ is the disjoint union of the cubes } R, R \in \cR. 
\end{equation}
The fact that $R \subset E \sm F$ comes from \eqref{f10}. Conversely,
if $x\in E \sm F$, then small cubes that contain $x$ satisfy \eqref{f9}, and
are contained in a cube of $\cR$ (because large cubes $Q$ fail \eqref{f9},
because $d_\cF (x_Q) \leq \dist(x_Q, Q_0) + l(Q_0) \leq \dist(x,Q_0) + 2 l(Q) + l(Q_0)
< \tau^{-1} l(Q)$ for $l(Q)$ large). Finally, the cubes of $\cR$ are 
disjoint by maximality.

We need a little more geometric information on $d_\cF$ and $F$ before we start.

\begin{lem} \label{Lf3}
Set $B_0 = B(x_{Q_0}, M_0 + 10)$. Then
\begin{equation} \label{f16}
\dist(z,E) \leq C \varepsilon d_\cF(z) \ \text{ for } z\in \Sigma \cap B_0,
\end{equation}
\begin{equation} \label{f17}
\dist(z,\Sigma) \leq C \varepsilon d_\cF(z) \ \text{ for } z\in E \cap B_0,
\end{equation}
and 
\begin{equation} \label{f18}
F \subset E \cap \Sigma \cap \ol Q_0.
\end{equation}
\end{lem}

We start with the easy part of \eqref{f18}. If $d_\cF(z) = 0$, then we can find
cubes $Q \in \cF$ such that $\dist(x,Q) + l(Q)$ is arbitrarily small, and since
$Q \subset Q_0 \subset E$, we get that $z\in E \cap \ol Q_0$. 

Next let $z\in B_0$ be given, 
pick $\delta > d_\cF(z)$ close to $d_\cF(z)$,
and choose a first cube $Q_1 \in \cF$ such that $\dist(z,Q_1) + l(Q_1) < \delta$.
Then let $Q$ denote the element of $\cF$ 
that contains $Q_1$ 
and whose generation
$k = k(Q)$ is the smallest possible, but with the constraint that $l(Q) \leq \delta$.
Such a cube $Q$ exists, since $Q_1$ satisfies the constraint.
First assume that $l(Q) < 1$. Then the parent of $Q$ does not satisfy the constraint,
even though it lies in $\cF$, and this forces $l(Q) \geq \delta/10$. 
Obviously $\dist(z,x_Q) \leq M l(Q)$, so Lemma \ref{Lf2} says that
\begin{equation} \label{f19}
\text{$\dist(z,E) \leq C \varepsilon l(Q)$ 
if $z\in \Sigma$, 
and $\dist(z,\Sigma) \leq C \varepsilon l(Q)$ if $z\in E$.}
\end{equation}
If instead $l(Q) = 1$, i.e., $Q = Q_0$, then $\dist(z,x_Q) \leq M l(Q)$
in this case too, because $z \in B_0$, and \eqref{f19} holds as well.

Recall that $l(Q) \leq \delta$ and we can pick any $\delta > d_\cF(z)$;
\eqref{f16} and \eqref{f17} follow. 
Also, in the case when $z\in F$,
we already know that $x\in E$, so \eqref{f19} says that $\dist(z,\Sigma)$
is as small as we want. Hence $z\in \Sigma$.
The lemma follows.
\qed

\ms
We are now ready to define a measure $\sigma$ on $\Sigma$
that approximates $\mu$ reasonably well. Since we have a nice set $F \subset E \cap \Sigma$,
we do not change the measure there, and set
\begin{equation} \label{f21}
\sigma_0 = \mu_{\vert F}.
\end{equation}
Next we consider the set 
\begin{equation} \label{f22}
\cR_0 = \big\{ R \in \cR \, ; \, \dist(R,Q_0) \leq 1\big\},
\end{equation}
and to simplify some of the notation, enumerate $\cR_0$ as a collection
$\{ R_j \}$, $j \in \cJ$. We want to replace each $\mu_j = \mu_{\vert R_j}$,
$j\in \cJ$, by a measure $\sigma_j$ with the same mass. Set
\begin{equation} \label{f23}
l_j = l(R_j) \ \text{ and } \ x_j = x_{R_j}.
\end{equation}
We want to take 
\begin{equation} \label{f24}
\Sigma(j) = \Sigma \cap B(x_j, M_2 l_j),
\end{equation}
where $M_2 = 10^{-1} \tau^{-1} \geq 10$, and 
\begin{equation} \label{f25}
\sigma_j = a_j \H^d_{\vert \Sigma(j)}, \ \text{ with } a_j = \frac{\mu(R_j)}{\H^d(\Sigma(j))},
\end{equation}
but we need to check some things. First observe that
\begin{equation} \label{f26}
l_j \leq 2\tau \inf_{z\in R_j} d_\cF(z) \leq 2\tau 
(1+\dist(R_j,Q_0) \leq 4\tau
\end{equation}
by \eqref{f10} and \eqref{f22}. Pick a first cube $Q' \in \cF$ such that 
\begin{equation} \label{f27}
l(Q') + \dist(x_j,Q') \leq 2d_\cF(x_j) \leq 40 \tau^{-1} l_j. 
\end{equation}
Set $Q_j = Q'$ if $l(Q') \geq 10l_j$, and otherwise let $Q_j$ be the ancestor
of $Q'$ such that $l(Q_j) = 10l_j$; notice that $Q_j\in \cF$
(by heredity and because \eqref{f26} says that $10l_j < 1$), and
\begin{equation} \label{f28}
\dist(x_j,Q_j) \leq 40 \tau^{-1} l_j 
\ \text{ and } \  10 l_j \leq l(Q_j) \leq 40 \tau^{-1} l_j.
\end{equation}
If $M$ is large enough (compared to $\tau^{-1}$), we may apply Lemma \ref{Lf2}, 
with the cube $Q_j$, and to the point $x_j \in E$. We find $\xi_j \in \Sigma$
such that $|\xi_j-x_j| \leq l_j$, and this is good because this implies that
\begin{equation} \label{f29}
\H^d(\Sigma(j)) \geq C^{-1} l_j^d.
\end{equation}
We also need to know that 
\begin{equation} \label{f30}
\text{the $\Sigma(j)$ have bounded overlap and do not meet $F$.} 
\end{equation}
and indeed, if $\xi \in \Sigma(j)$, then 
\begin{equation} \label{f31}
|\xi -x_j| \leq M_2 l_j \leq M_2 \tau d_\cF(x_j) = 10^{-1} d_\cF(x_j)
\end{equation}
by \eqref{f9} and the definition of $M_2$, hence  
\begin{equation} \label{f32}
0 < \frac9{10} d_\cF(x_j) \leq d_\cF(\xi) \leq \frac{11}{10} d_\cF(x_j)
\ \text{ for } \xi \in \Sigma(j).
\end{equation}
Hence $\xi \notin F$, and also the size of $d_\cF(\xi)$ determines
roughly the generation of $j$; 
\eqref{f30} follows at once. We complement the measures $\sigma_j$ by 
\begin{equation} \label{f33}
\sigma_\infty = \H^d_{\vert \Sigma(\infty)},
\ \text{ with } \Sigma(\infty) = \big\{ \xi\in \Sigma \, ; \, \dist(\xi,Q_0) \geq 1/2 \big\}.
\end{equation}
Finally we set
\begin{equation} \label{f34}
\sigma = \sigma_0 + \sigma_\infty + \sum_{j\in \cJ} \sigma_j.
\end{equation}
It will be good to know that 
\begin{equation} \label{f35}
\sigma \ \text{ is an Ahlfors regular measure with support } \Sigma, 
\end{equation}
with the AR constant depending on $n,d,C_0$ only, provided that $M$ is sufficiently large 
and $\eps_0$ is sufficiently small. The verification is twofold. 
First we check that the density $f$ of $\sigma$ with respect
to $\H^d$ is bounded. On $F$, this is because $\mu$ is Ahlfors regular, and the measures
$\sigma_j$ and $\sigma_\infty$ do not charge $F$. Concerning the $\sigma_j$, 
their density $a_j$ is bounded, by \eqref{f29} and because $\mu(R_j) \leq C l_j$,
and then the global density is bounded because of \eqref{f30}.
As for $\sigma_\infty$, its density is $1$. 

Conversely, $f$ is bounded from below on $F$ (by \eqref{f21}) and on $\Sigma(\infty)$.
Now let $\xi \in \Sigma \sm (F \cup \Sigma(\infty))$ be given. 
Thus $0 < d_\cF(\xi) < 3/2$ (because $\dist(\xi,Q_0) \leq 1/2$). Let 
$Q' \in \cF$ be such that $l(Q')+ \dist(\xi,Q') \leq 2 d_\cF(\xi)$. 
Keep $Q=Q'$ if $l(Q') \geq d_\cF(\xi)/100$, and otherwise replace it with an ancestor
$Q$ such that $d_\cF(\xi)/100 \leq l(Q) \leq d_\cF(\xi)/10$. Notice that $Q\in \cF$
because $l(Q) \leq 1$. Thus
\begin{equation} \label{f36}
\dist(\xi,Q) \leq 2 d_\cF(\xi) \ \text{ and } \  d_\cF(\xi)/100 \leq l(Q) \leq 2 d_\cF(\xi).
\end{equation}
We can apply Lemma \ref{Lf2} to $\xi$, and find $x\in E$ such that
\begin{equation} \label{f37}
|x-\xi| \leq C \varepsilon l(Q) \leq C \varepsilon d_\cF(\xi);
\end{equation} 
then $d_\cF(x) \geq d_\cF(\xi)/2$, $x\in E \sm F$, there is a cube $R \in \cR$ 
that contains $x$, $\dist(R,Q_0) < 1$ because $x\in R$ and $\dist(\xi,Q_0) \leq 1/2$,
and hence $R$ is one of the $R_j$. In addition 
$l_j \leq \tau d_\cF(x_j)$ by \eqref{f9}, so
$d_\cF(x) \leq d_\cF(x_j) + l_j \leq 2 d_\cF(x_j)$ (by \eqref{a3}), 
$d_\cF(\xi) \leq 2 d_\cF(x) \leq 4 d_\cF(x_j) \leq 80\tau^{-1} l_j$
by \eqref{ff15}, and \eqref{f37} says that $\xi$ lies well inside $\Sigma(j)$.
The coefficient $a_j$ is also bounded from below, so $f \geq a_j \geq C^{-1}$
near $\xi_j$; \eqref{f35} follows.

\ms
We are now about ready to prove that $\Psi$ in \eqref{f3} satisfies a Carleson measure
condition. By \eqref{ff10}, $\frac{D_{E,\mu}}{D_{\Sigma,\sigma}}$ is bounded 
and bounded from below on $\Omega_\cF$ (because $D_{E,\mu}$ is equivalent 
to the distance to $E$,  and $D_{\Sigma,\sigma}$ to the distance to $\Sigma$), 
so it is enough to prove that
$\Psi_1$ satisfies a Carleson measure condition on $\Omega_\Sigma = \R^n \sm \Sigma$, where
\begin{equation} \label{f38}
\Psi_1(z) := \1_{\Omega_\cF} \dist(z,E)^{-1} \,\Big| D_{E,\mu}(z) - D_{\Sigma,\sigma}(z) \Big|.
\end{equation}
Notice that we chose $\dist(z,E)^{-1}$ because it seems simpler, but $\dist(z,\Sigma)^{-1}$
would have been equivalent. This remark also implies that although we decided to advertise
a Carleson measure condition on $\Omega_\Sigma$, we would 
obtain a Carleson measure condition on $\Omega_E = \R^n \sm E$ just the same way.

Given the definition of $D_{E,\mu}(z)$ and $D_{\Sigma,\sigma}$ and the same equivalences as
above, it will be equivalent, and simpler, to prove that $\Psi_2$ satisfies a Carleson measure 
condition, where
\begin{equation} \label{f39}
\Psi_2(z) := \1_{\Omega_\cF} \dist(z,E)^{\alpha} 
\,\big| D_{E,\mu}(z)^{-\alpha} - D_{\Sigma,\sigma}(z)^{-\alpha} \big|.
\end{equation}
So we gives ourselves $z\in \Omega_\cF$, and we want to estimate
\begin{eqnarray} \label{f40}
\Delta(z) &:=& \dist(z,E)^{-\alpha} \Psi_2(z) 
= \big| D_{E,\mu}(z)^{-\alpha} - D_{\Sigma,\sigma}(z)^{-\alpha} \big|
\nn\\
&=& \Big| \int |z-y|^{-d-\alpha} [d\mu(y)-d\sigma(y)]\Big|.
\end{eqnarray}
Recall that $\mu_j = \mu_{\vert R_j}$; let us also set 
$\mu_0 = \mu_{\vert F}$ and 
\begin{equation} \label{f42}
\mu_\infty = \mu_{\vert E(\infty)}, \ \text{ with } 
E(\infty) = E \sm \Big( F \cup \bigcup_{j \in \cJ} R_j \Big),
\end{equation}
so that $\mu = \sum_j \mu_j$ (a sum that includes $j=0$ and $j=\infty$).
Naturally we write $\Delta(z) \leq \sum_j \Delta_j(z)$, where
\begin{equation} \label{f41}
\Delta_j(z) = \Big| \int |z-y|^{-d-\alpha} [d\mu_j(y)-d\sigma_j(y)]\Big|.
\end{equation}
A priori the sum contains $0$ as well as $\infty$, but 
the term with $j=0$ drops, because $\mu_0 = \sigma_0$.
We now give an estimate on $\Delta_j(z)$ for $j \in \cJ$. Set
\begin{equation} \label{f43}
d_j(z) = \dist(z, R_j \cup \Sigma(j)).
\end{equation}
We claim that
\begin{equation} \label{f44}
d_j(z) \geq  C^{-1} \, l_j, \ \text{ where we may take } 
C = \max(M_2, 
3\tau (M_0+1) M_0)  > 0.
\end{equation}
Me may assume that $d_j(z) \leq  M_2 l_j$, because otherwise \eqref{f44} holds.
Then $|z-x_j| \leq 2M_2 l_j$ (see \eqref{f24} and the definition of $M_2$ below \eqref{f24}).
Hence 
$$
d_\cF(z) \geq d_\cF(x_j) - 2M_2 l_j \geq (1 - 2 M_2 \tau) d_\cF(x_j) 
\geq d_\cF(x_j)/2 \geq (2\tau)^{-1} l_j
$$
by \eqref{f9}, the definition of $M_2$ below \eqref{f24}, and \eqref{f9}.

Let  $Q \in \cF$ be such that $z\in W(Q)$; then $d_\cF(z) \leq l(Q) + \dist(z,Q)
\leq (M_0+1) l(Q)$ by the definition  
\eqref{f2}, and now
$$
l_j \leq 2\tau d_\cF(z) \leq 2\tau (M_0+1) l(Q) \leq 2\tau (M_0+1) M_0 \dist(z,E)
$$
by \eqref{f2}. On the other hand,
$$
d_j(z) = \dist(z, R_j \cup \Sigma(j)) \geq \dist(z, E \cup \Sigma)
\geq \frac23 \dist(z, E)
$$
by \eqref{ff10}; our claim \eqref{f44} follows. 

In the next computations, we no longer record the dependence of our various constants
on $M_0$, $M_2$, or $\tau$. We now use \eqref{f44} to prove that
\begin{equation} \label{f45}
\Delta_j(z) \leq C\mu(R_j) l_j  [l_j + d_j(z)]^{-d-\alpha-1} 
\leq C l_j^{d+1} [l_j + d_j(z)]^{-d-\alpha-1}.
\end{equation}
Set $\delta_j = \diam(\Sigma(j) \cup R_j) \leq C l_j$. 
When  $d_j(z) \leq  2\delta_j$, we just use \eqref{f44} and the fact that the total 
masses of $\sigma_i$ and $\mu_j$ are 
$\mu(R_j) \leq C l_j^d$ 
to get the result.
Otherwise, set $a_0 = |z-x_i|^{-d-\alpha}\mu(R_j)$,
$a_1 = \int |z-y|^{-d-\alpha} d\sigma_j(y)$, and $a_2 = \int |z-y|^{-d-\alpha} d\mu_j(y)$.
Notice that
\begin{equation} \label{f46}
a_1 - a_0 = \int_{\Sigma(j)} \Big[ |z-y|^{-d-\alpha} - |z-x_i|^{-d-\alpha}\Big] d\mu_j(y),
\end{equation}
then observe that for $y\in \Sigma(j)$, 
$$
\Big| |z-y|^{-d-\alpha} - |z-x_i|^{-d-\alpha}\Big| \leq C \delta_j |z-x_i|^{-d-\alpha-1}
$$
(differentiate the integrand along the line segment $[y,x_i]$); this yields
$$
|a_1 - a_0| \leq \delta_j |z-x_i|^{-d-\alpha-1} ||\sigma_j|| 
\leq C l_j \mu(R_j)  d_j(z)^{-d-\alpha-1}.
$$
We have the same estimate on $|a_2 - a_0|$, and \eqref{f45} follows.

\ms
At this point we have enough information to prove the desired Carleson bound
on $\Psi_3$. The most important part will come from the sum over $j\in \cJ$. 
Each $\Delta_j$ gives a bump function with an $L^2$ norm controlled by $\mu(R_j)$,
and sufficiently localized or smooth for the different pieces to be almost
orthogonal. The computations that follow are reminiscent of other computations
done in a similar context, but it seems that we need to be courageous and do them.

We take a Carleson box $B = B(X,R)$ centered on $\Sigma$ 
(but $E$ would give the same), 
and want to prove that 
\begin{equation} \label{f47}
J := \int_{z\in \Omega_\cF \cap B} \Big|\sum_j \Delta_j(z)\Big|^2 
\delta(z)^{d-n+2\alpha} dz
\leq C R^d.
\end{equation}
where we set $\delta(z) = \dist(z,E)$ (but $\dist(z,\Sigma)$ is equivalent on $\Omega_\cF$),
and the extra $2\alpha$ come from the fact that we have to multiply $\sum_j \Delta_j(z)$) 
by $\delta(z)^\alpha$ before checking the Carleson condition; see \eqref{f39}.

We may assume that $R \leq 1$, because $\Omega_\cF \leq B(x_{Q_0},M_0)$
anyway. Also we first concentrate on 
\begin{equation} \label{f47b}
J _1:= \int_{z\in \Omega_\cF \cap B} \Big|\sum_{j\in \cJ(B)} \Delta_j(z)\Big|^2 
\delta(z)^{d-n+2\alpha} dz,
\end{equation}
where $\cJ(B)$ is the collection of $j\in \cJ$ for which 
$R_j \subset CB$.
The value of $C$ will be decided when we deal with the rest of the sum.
 
We write $J_1 \leq  2\sum_i \sum_j J(i,j)$, where it is enough
to sum on the pairs such that $l_i = l(R_i) \leq l_j = l(R_j)$, and
\begin{eqnarray} \label{f48}
J(i,j) &=& \int_{z\in \Omega_\cF \cap B} \Delta_i(z) \Delta_j(z) \delta(z)^{d-n+2\alpha} dz
\nn\\
&\leq& C \int_{z\in \Omega_\cF \cap B} 
\frac{l_j^{d+1}}{[l_j + d_j(z)]^{d+\alpha+1}}
\, \frac{l_i^{d+1}}{[l_i + d_i(z)]^{d+\alpha+1}}
\, \delta(z)^{d-n+2\alpha} dz.
\end{eqnarray}
Recall that $d_j(z) = \dist(z,R_j \cup \Sigma(j))$, but in fact here \eqref{ff10}
says that $d_j(z) \geq \frac12 \dist(z,\Sigma \cup E) \geq \frac12 \delta(z)$.
Then $(d_j(z)+l_j)^{-1} \leq d_j(z)^{-1} \leq 2 \delta(z)^{-1}$.
We may use this to replace a negative power of $(d_j(z)+l_j)$ by the same power of
$\delta(z)$ and simplify some things. And we can do the same thing with $i$.
Set
\begin{equation} \label{f49}
r_i(z) = \frac{l_i^{d+1}}{[l_i + d_i(z)]^{d+1}} \leq 1;
\end{equation}
we first use the fact that
\begin{equation} \label{f50}
\frac{l_i^{d+1}}{[l_i + d_i(z)]^{d+\alpha+1}} 
= r_i(z)  [l_i + d_i(z)]^{-\alpha} \leq 2^\alpha r_i(z) \delta(z)^{-\alpha}
\end{equation}
to get that
\begin{equation} \label{f51}
J(i,j) \leq C \int_{z\in \Omega_\cF \cap B} 
\frac{l_j^{d+1}}{[l_j + d_j(z)]^{d+\alpha+1}} \, r_i(z) \, \delta(z)^{d-n+\alpha} dz.
\end{equation}
Next we divide the integral into annuli $A_{j,k}$ where $d_j(z) \sim 2^k l_j$,
and get contributions $J(i,j,k)$. The smallest annulus should be replaced by a ball, 
but we still get the same estimate, namely
\begin{equation} \label{f52}
J(i,j,k) \leq C 2^{-k(d+\alpha +1)} l_j^{-\alpha}
\int_{z\in \Omega_\cF \cap B \cap A_{j,k}} \, r_i(z) \, \delta(z)^{d-n+\alpha} dz.
\end{equation}
It will be good to know that
\begin{equation} \label{f53}
\int_{B_r} \delta(z)^{d-n+\alpha} dz \leq C r^{d+\alpha}
\end{equation}
when $B_r$ is a ball of radius $r$ centered on $\Sigma$ (this is just easier). 
This estimate is very easy when $\Sigma$ is a $d$-plane; the main point then 
is that the integral in the direction orthogonal to $\Sigma$ converges because of 
the additional exponent $\alpha$. When $\Sigma$, as here, is bilipschitz-equivalent
to a $d$-plane, this is as easy because we can change variables. But this would also
be true with $E$, with just a bit more work, because it is Ahlfors regular and the measure 
of tubes of width $tr$ near $E$ are easy to estimate.

Denote by $I(i,j,k)$ the integral on the right of \eqref{f52} and further cut the domain
of integration into annuli $A_l$ where $d_i(z) \sim 2^l l_i$. We get integrals
\begin{equation} \label{f54}
I(i,j,k,l) \leq \int_{A_{j,k} \cap A_l} r_i(z) \delta(z)^{d-n+\alpha} dz
\leq C 2^{-l(d+1)} \min[(2^l l_i)^{d+\alpha}, (2^kl_j)^{d+\alpha}],
\end{equation}
where the first piece is an estimate of $r_i(z)$ and the second one comes 
from \eqref{f53}, with the two different choices of diameter.
For strategic reasons (we do not want to distinguish between $\alpha \leq 1$ 
and $\alpha > 1$, we choose $\tau \in (0,1]$, smaller than $\alpha$, and replace
$2^{-l(d+1)}$ with $2^{-l(d+\tau)}$ in \eqref{f54}.
No relation with our previous constant $\tau$, though; this one will just be here for the duration of
the computation.

We need to be careful about the region of integration. Set
$d_{i,j} = \dist(R_i, R_j)$, and first asume that $d_{i,j} \leq C 2^k l_j$.
In this case $A_{j,k}$ stays within $C 2^k l_j$ of $R_i$, and we can content
ourselves with $l$ such that $2^l l_i\leq C 2^k l_j$, because for larger ones $A_l$
does not meet $A_{j,k}$. Then use the first option in \eqref{f54} and observe that
\begin{equation} \label{f55}
I(i,j,k) \leq \sum_l I(i,j,k,l) \leq C \sum_l 2^{-l(d+\tau)} (2^l l_i)^{d+\alpha}
= C \sum_l 2^{l(\alpha-\tau)} l_i^{d+\alpha}
\end{equation}
and the largest terms are when $l$ is as large as possible, i.e., when
$2^l l_i\sim 2^k l_j$. 
This yields $I(i,j,k) \leq (2^k l_j/l_i)^{\alpha-\tau} l_i^{d+\alpha}$
and 
\begin{equation} \label{f56}
J(i,j,k) \leq C 2^{-k(d+\alpha +1)} l_j^{-\alpha} I(i,j,k)
\leq C 2^{-k(d+\tau +1)} l_j^{-\tau} l_i^{d+\tau}.
\end{equation}
Now assume that $d_{i,j} \geq C 2^k l_j$.
If $C$ was chosen large enough, all the points of the annulus $A_{j,k}$ lie at distance
roughly $d_{i,j}$ from $R_i$, which means that we just need to sum over the few
$l$ such that $2^l l_i \sim d_{i,j}$. For all these $l$, it is actually better to
use the second option in \eqref{f54}, which yields
$I(i,j,k,l) \leq C 2^{-l(d+\tau)}(2^kl_j)^{d+\alpha} \leq C (d_{i,j}/l_i)^{-d-\tau}
(2^kl_j)^{d+\alpha}$ and then 
\begin{equation} \label{f57}
\begin{aligned}
J(i,j,k) &\leq C 2^{-k(d+\alpha +1)} l_j^{-\alpha} I(i,j,k)
\leq C 2^{-k(d+\alpha +1)} l_j^{-\alpha} (d_{i,j}/l_i)^{-d-\tau} (2^kl_j)^{d+\alpha}
\cr&\leq C 2^{-k} d_{i,j}^{-d-\tau} l_i^{d+\tau} l_j^{d}.
\end{aligned}
\end{equation}
Now we have to sum all these numbers. We start with the first case, and sum first 
over $i\in \cJ$ such that $l_i = 2^{-m}l_j$ for a given $m \geq 0$. 
These cubes lie in a ball of size roughly $2^k l_j$ (because $d_{i,j} \leq C 2^k l_j$), 
so there is roughly $(2^k l_j/(2^{-m}l_j))^d = 2^{(m+k)d}$ of them. 
We get a sum bounded by
$$
C 2^{(m+k)d} 2^{-k(d+\tau+1)} l_j^{-\tau} (2^{-m}l_j)^{d+\tau} 
$$
The exponent for $m$ is $-\tau$, so we may sum over $m$ and get
$C 2^{kd} 2^{-k(d+\tau+1)} l_j^{d}$. Then we sum over $k$
and get $C l_j^{d} \leq C \mu(R_j)$. Then we sum over $j$, recall that the
$R_j$ are disjoint and we only sum over those that are contained in $CB$,
and get less than $C R^d$, as needed.

Now we consider the second case and sum the terms from \eqref{f57}.
Fix $j, k$, and $m\geq 0$, and sum over $i$ such that $l_i = 2^{-m}l_j$.
Further decompose into annuli where $d_{i,j} \sim 2^n 2^k l_j$.
The number of cubes $R_i$ in the annulus is less than $C 2^{(n+k+m)d}$, 
which gives a contribution smaller than
$$
C 2^{(n+k+m)d} 2^{-k} d_{i,j}^{-d-\tau} l_i^{d+\tau} l_j^{d}
= C 2^{(n+k+m)d} 2^{-k} (2^n 2^k l_j)^{-d-\tau} (2^{-m}l_j)^{d+\tau} l_j^{d}
$$
The power for $m$ is $-\tau$, so we sum over $m$ and forget about it; 
we are left with $C 2^{-k} 2^{-\tau(n+k)} l_j^{d}$. We sum over $k$, $n$,
get less than $C l_j^d$, which as before we can sum over $j$ to get at most 
$CR^d$.

We are now finished with $J_1$, but we still need to estimate
\begin{equation} \label{f59}
J _2 = \int_{z\in \Omega_\cF \cap B} \big|\wt\Delta(z) \big|^2 
\delta(z)^{d-n+2\alpha} dz,
\end{equation}
where $\wt\Delta$ is the remaining part of $\Delta$, i.e.,
\begin{equation} \label{f60}
\wt\Delta(z) = \sum_{j \in \cJ_e \cup \{ \infty \}} \Delta_j(z),
\end{equation}
where $\cJ_e$ is the set of indices $j\in \cJ$ such that $R_j$ is not 
contained in $CB$. 

We start with $\cJ_e$. Let $j\in \cJ$ be given.
One possibility is that $l(R_j) \geq R$; then \eqref{f44} implies that 
$\dist(z,R_j\cup \Sigma(j)) \geq C^{-1}R$. Otherwise, and if $C$ is chosen large
enough, the fact that $R_j$ is not contained in $CB$ implies that
$R_j \cup \Sigma(j)$ does not meet $2B$, and hence $\dist(z,R_j\cup \Sigma(j)) \geq R$.

Now set $\mu_e = \sum_{j \in\cJ_e} [\mu_j-\sigma_j]$, and observe that
\begin{equation} \label{f61}
\Delta_e(z) := \sum_{j \in \cJ_e} \Delta_j(z) \int |z-y|^{-d-\alpha} d\mu_e(y)
\end{equation}
is bounded by $C R^{-\alpha}$, so 
\begin{equation} \label{f62}
\int_{z\in \Omega_\cF \cap B} \big|\Delta_{e}(z) \big|^2 \delta(z)^{d-n+2\alpha} dz 
\leq C R^{-2\alpha} \int_{z\in B} \delta(z)^{d-n+2\alpha} dz 
\leq C R^d,
\end{equation}
as needed. We are left with $\Delta_\infty$. We need to control
\begin{equation} \label{f63}
J _3 = \int_{z\in \Omega_\cF \cap B} \big|\Delta_\infty(z)\big|^2 
\delta(z)^{d-n+2\alpha} dz.
\end{equation}
Observe that 
\begin{equation} \label{f64}
\dist(z,\Sigma(\infty)\cup E(\infty)) \geq C^{-1}
\end{equation}
because either $z\in W(Q)$ for a cube $Q$ such that
$\diam(W(Q)) \leq 10^{-1}$, and then \eqref{f64} holds because
$Q$ is centered on $Q_0$ and $\Sigma(\infty)\cup E(\infty)$ is far from
$Q_0$ (see the definitions and \eqref{f33} and \eqref{f42}), or else
\eqref{f64} comes directly from \eqref{f2} and \eqref{ff10}.

Once we know \eqref{f64}, we see that $\Delta_\infty(z) \geq C^{-1}$
and \eqref{f63} follows as for $\Delta_e$.

This concludes our proof of Proposition \ref{Pf1}.
\qed

\section{Construction of the extrapolation saw-tooth region and related results for the harmonic measure based on Sections~\ref{S2}--\ref{Sf}, \label{sSum}}
\subsection{Collecting the results of Sections~\ref{S2}--\ref{Sf}}
In this section we summarize the results of Sections~\ref{S2}--\ref{Sf} in application to the particular stopping time region which we will use for the extrapolation procedure in the forthcoming discussion. At this point let us start by recalling Tolsa's $\alpha$-numbers \cite{To}, which give a good control on sums of squares of local Wasserstein 
distances to flat measures, for every Ahlfors-regular measure on a uniformly rectifiable set.

First, recall the Tolsa's $\alpha$-coefficients defined in \eqref{2.5bis}. Fix some $M>0$ which will eventually be chosen sufficiently large. Given a $d$-dimensional  Ahlfors regular set $E$ and an Ahlfors regular measure $\mu$ on $E$, equipped with the usual dyadic grid $\dd=\dd(E)$ (see Section~\ref{intro}), we denote 
$$  \alpha (Q):= \alpha (x_Q, M\, l(Q)).$$

Recall the notion of uniform rectifiability from Definition~\ref{def-ur}. 
In Theorem~1.2, \cite{To}, Tolsa proves the following result 
(a modification of which corresponding to balls rather than dyadic cubes we already used in Section~\ref{SDist}). 
\begin{theorem}\label{tUR-alpha}
For every uniformly rectifiable set $E$,  every Ahlfors regular measure $\mu$ on $E$, and every
dyadic cube $R\in \dd(E)$, 
\begin{equation} \label{a6.6}
\sum_{Q \in \dd(R)}   \alpha (Q)^2 \mu(Q) \leq C \mu(R)^d,
\end{equation}
where 
the constant $C$ depends on $M$, $n$, $d$, $C_0$, and the UR constants of $\mu$. 
\end{theorem}
Strictly speaking, \cite{To} defines $\alpha(Q)$ slightly differently, indexed by the dyadic cubes of $\RR^n$ rather than (pseudo-)dyadic cubes of $E$ and using $M=3$, but possibly adjusting the values of $M$ (say, directly in the proof), the two statements are equivalent. Indeed, every $B(x_Q, M\, l(Q))$ for $Q\in \dd(E)$ is contained in some $B(x_{Q'}, 2M\, l(Q'))$ for $Q'\in \dd(\RR^n)$ with $l(Q')\approx l(Q)$ and the number of different $Q$'s corresponding to the same $Q'$ is uniformly bounded.

Now let us return to the construction of a stopping time region. 
\begin{definition}\label{def-st} Let $E$ be a $d$-dimensional  Ahlfors regular set and let $\mu$ be an Ahlfors regular measure on $E$, equipped with the usual dyadic grid $\dd=\dd(E)$. Fix some $\eps_0>0$, $\delta_0>0$, and  $Q_0\in \dd$. Then the stopping time region $\cF=\cF_{\eps_0, \delta_0}(Q_0)$ is constructed as follows (cf. the procedure in Section~\ref{intro}). We start from the top cube $Q_0$, and decide to remove a cube
$Q \in \bD(Q_0)$, as well as all its descendants, as soon as 
\begin{equation}\label{eq-st}
\mbox{either $\alpha (Q) > \varepsilon_0$
or } J_\alpha(Q) := \sum_{k(Q_0) \leq k \leq k(Q)}  \alpha (R_k(Q))^2 \geq \delta_0.
\end{equation}
Here, as before, for $k(Q_0) \leq k \leq k(Q)$, we denote by $R_k(Q)$ the cube of $\bD_k$ that contains $Q$.
The remaining collection of cubes will be referred to as $\cF=\cF_{\eps_0, \delta_0}(Q_0)$. For a fixed $M_0>0$ a saw tooth region based on $\cF$ is, as before, $\Omega_\cF$ defined in \eqref{f1}--\eqref{f2}.

It has become customary to also remove the siblings of any cube $Q$ that we remove as above, because this gives a little more regularity to the decomposition of $\bD$ into stopping time regions like $\cF$, and this costs essentially nothing. Here we do not need this because of our
specific description, but it would not hurt either.

It is convenient to write $\dd(Q_0)\setminus \cF$ as $\cup_{Q_j\in \kF} \dd(Q_j)$ 
where $\kF=\{Q_j\}_j$ is a disjoint collection 
of cubes $Q_j\in \dd(Q_0)$ 
maximal under our stopping time procedure. That is, the collection  $\kF=\{Q_j\}_j$ consists exactly of maximal cubes 
$Q \subset Q_0$ satisfying 
\eqref{eq-st} and the entire collection of all removed cubes (that is, cubes of $\kF$ and its descendents) is then $\cup_{Q_j\in \kF} \dd(Q_j)$. 
\end{definition}

\begin{remark}\label{r8.6} At this point and throughout Section~\ref{sSum} we assume that $\kF\neq \{Q_0\}$, for otherwise $\cF =\emptyset$ and there is no $\Sigma$ to be constructed. 

We remark also that we do not require $Q_0\in \dd_0$ (i.e., $k(Q_0)=0$ hereafter). This was often an assumption in previous chapters but the corresponding results rescale easily -- we will mention this in due time. 
\end{remark}

We shall now collect results from Sections~\ref{S2}--\ref{Sf}, to arrive
at the following. 

\begin{thm}\label{t8.8}
Let $E$ be a $d$-dimensional  Ahlfors regular set and $\mu$ be an Ahlfors regular measure on $E$. 
Then for $M_0>1$ large enough depending on $n, d, C_0$, $M>1$ large enough depending on 
$n, d, C_0, M_0$, and $\eps_0, \delta_0> 0$ small enough depending on $n, d, C_0, M_0, M$, 
for any $Q_0\in \dd(E)$ and the associated $\cF=\cF_{\eps_0, \delta_0}(Q_0)$ built as in Definition~\ref{def-st}, there exists a $d$-Ahlfors regular set $\Sigma \subset \RR^n$ and a $d$-Ahlfors 
regular measure $\sigma$ supported on $\Sigma$ with the following properties: 
\begin{enumerate}
\item for any $\alpha>0$ the harmonic measure associated to 
$L_\sigma=-\div D_{\sigma}^{-(n-d-1)} \nabla$,  with 
\begin{equation} \label{eq8.9}
D_{\sigma}(X) = \Big\{ \int_\Sigma |X-y|^{-d-\alpha} d\sigma(y) \Big\}^{-1/\alpha},
\end{equation}
is $A_\infty$ with respect to $\sigma$;
\item for any $\alpha>0$ the harmonic measure associated to 
$\wh L=-\div \wh D^{-(n-d-1)} \nabla$,  with 
$\wh D$ defined in \eqref{f4}--\eqref{f5},
is $A_\infty$ with respect to $\sigma$.
\end{enumerate}
If $M_0, M$ are chosen large enough and $\eps_0, \delta_0$ are chosen small enough as postulated above, the 
$A_\infty$ constants 
of $\omega_{L_\sigma}$ and $\omega_{\wh L}$ and the Ahlfors regularity constant of $\sigma$ depend only on $n, d, C_0, \alpha$.
\end{thm}

\bp At this point the proof is a collection of results from Sections~\ref{S2}--\ref{Sf}. 

 First of all, since the assumptions and results are scale invariant, we assume without loss of generality that $Q_0\in \dd_0$. 
 
 Going further, our stopping time region falls under the scope of  Sections~\ref{S2}--\ref{Sf}. Indeed, the list of restrictions on $\cF$ in Sections~\ref{S2}--\ref{Sf} is exhausted by two properties, \eqref{a9} and \eqref{a12}. At the same time, as pointed out in Section~\ref{intro}, the numbers $\alpha (Q)$ control the properties \eqref{a9} and \eqref{a12}. To be precise, if $J_\alpha(Q) \leq \delta_0$ and we use a suitable fixed multiple of $M$ in place of $M$ in the definition of $J_\alpha(Q)$, then $J(Q)\leq \delta_1$, with $\delta_1$ being a fixed multiple of $\delta_0$. This is due to \eqref{a6.5} and the fact that $\alpha$ numbers control (bilateral) $\beta_1$ numbers proved in Lemma~3.2 in \cite{To}. Furthermore, bilateral $\beta_1$ numbers control powers of bilateral $\beta_\infty$ numbers (see \cite{DS}, p. 27) and hence, slightly adjusting the choice of $M$ as above, similar considerations assure that the condition $\alpha (Q)\leq \varepsilon_0$ implies \eqref{a9} with $\eps_1$ being a fixed multiple of a power of $\eps_0$. Here, the power depends on the dimension only and by a ``fixed multiple" we mean multiplication by a numerical constant which is allowed to depend on dimension only.
 
Thus, we can follow the construction of the closed set $\Sigma\in \RR^n$  through Sections~\ref{S2}--\ref{Sf}. According to Lemma~\ref{Lb1}, if $\eps_0$ (and hence, $\eps_1$) is small enough, there exists a closed set $\Sigma\in \RR^n$ and  $g: \RR^n \to \RR^n $, mapping  $\RR^d$ to $\Sigma$, bilipschitz, with Lipschitz constants depending on $n, d, C_0, \delta_0.$ Hidden in this statement is the dependence on choices of $M_0$ and $M$ as well (as the choice of $\eps_0$ ultimately depends on them) but all this is harmless. We remark that 
since $g$ is bilipschitz, $\Sigma$ is $d$-Ahlfors regular, with
constants depending on the same parameters as the biLipschitz constant for $g$.
Furthermore, by Proposition~\ref{Pd1} the matrix $J=Dg\, Q$ has the a form \eqref{d1}, where $A^1$ is bounded and invertible
by comment before the proof of Lemma~\ref{pRd3}.

Recalling that the $R_{k(t)}(x)$ are isometries, it follows that 
\begin{equation} \label{a5.2}
\Ab(x,t) = |\det(J(x,t))| (J(x,t)^{-1})^T J(x,t)^{-1}
\end{equation}
has the same structure as $J$ (with the same control of the Carleson and $L^\infty$ norms of its components), due to 
Lemma ~\ref{pRd3}
That is, 
\begin{equation} \label{d1-bis}
\Ab(x,t)= \left( \begin{array}{cc}
\Ab^1(x,t) & \CC^2(x,t) \\ \CC^3(x,t) & 
I_{n-d}+\CC^4(x,t)  
\end{array} \right), 
\end{equation}
where $\CC^2$, $\CC^3$, and $\CC^4$ 
are bounded (with an $L^\infty$ constant which goes to 0 as  $\eps_0\to 0$) and satisfy Carleson measure conditions (with a constant which goes to 0 as $\delta_0\to 0$). Also, if $\delta_0$, $\eps_0$ are small enough (depending on $n, d, C_0, M_0, M$ only), then $\Ab^1$ and $\Ab$ are bounded and invertible, with uniform bounds.

At this point we introduce the
$d$-dimensional Ahlfors regular measure $\sigma$ on $\Sigma$,
with a uniform control on the AR constant and good approximation properties,
which has been constructed in Section~\ref{Sf} (see \eqref{f35}), and use it
to define the two matrices
\begin{equation} \label{}
\A_\sigma(x, t):=\left(\frac{|t|}{D_{\sigma}(g(x, t))}\right)^{n-d-1} \Ab(x,t)
\, \text{ and } \, \wh\A:=\left(\frac{|t|}{\wh D (g(x, t)}\right)^{n-d-1} \Ab(x,t),
\end{equation}
defined for $(x, t)\in \RR^n$. These define degenerate elliptic operators
on $\R^n\sm \R^d$, and we want to show that they
satisfy the conditions imposed 
in Theorem~\ref{Itsf1}. 

First we observe that $\frac{|t|}{D_{\sigma}(g(x, t))}$ and $\frac{|t|}{\wh D (g(x, t)}$ are both bounded from above and below on $\RR^n\setminus\RR^d$ and hence, $\A$ and $\wh A$ are bounded and elliptic.
Next, let us write 
\begin{multline} \label{d1-bis}
\A_{\sigma}(x, t)=\left(\frac{|t|}{D_{\sigma} (g(x, t)}\right)^{n-d-1} \Ab(x,t)\\[4pt]= \left(\frac{|t|}{D_{\sigma} (g(x, t)}\right)^{n-d-1}\left( \begin{array}{cc}
\Ab^1(x,t) & 0 \\ 0& 
I_{n-d} 
\end{array} \right)+ \left(\frac{|t|}{D_{\sigma} (g(x, t)}\right)^{n-d-1}\left( \begin{array}{cc}
0 & \CC^2(x,t) \\ \CC^3(x,t) & 
\CC^4(x,t)  
\end{array} \right). 
\end{multline}
The second term satisfies the same Carleson measure conditions as the original $\CC^j$'s, $j=2,3,4$, since the multiplicative factor is bounded from above and below. The 
multiple of $\Ab_1$ is harmless and anyway we did not impose any specific
conditions on $\A_1$. 
It remains to analyze $\left(\frac{|t|}{D_\sigma (g(x, t)}\right)^{n-d-1} I_{n-d}$. 
Turning to this task, we recall that according to Proposition~\ref{Pe1}
\begin{equation} \label{e28-bis}
\Big|\frac{D_\sigma(g(y,t))}{|t|} - C_\alpha \lambda_\sigma(f(y),|t|)^{-1/\alpha} \Big| 
\end{equation}
satisfies a Carleson condition on $\RR^n\setminus \RR^d$. 
Since both  $\frac{D_\sigma(g(y,t))}{|t|}$ and $\lambda_\sigma(f(y),|t|)$ are bounded from above and below, applying the 
fundamental theorem of calculus
to $s\mapsto s^{-(n-d-1)}$, we deduce that 
\begin{equation} \label{e28-bis2}
\Big|\left(\frac{|t|}{D_\sigma(g(y,t))}\right)^{n-d-1} - C'_\alpha \lambda_\sigma(f(y),|t|)^{\frac{n-d-1}{\alpha}} \Big| 
\end{equation}
satisfies a Carleson condition on $\RR^n\setminus \RR^d$. Thus, a multiple of $I_{n-d}$ by the expression in \eqref{e28-bis2} can be absorbed into $\C^4$ and we are left with 
$$ C'_\alpha \lambda_\sigma(f(y),|t|)^{\frac{n-d-1}{\alpha}}=:b.$$
Since $\lambda$ is bounded from above and below, the condition \eqref{a2.32} is verified. On the other hand, \eqref{a2.33} follows immediately from \eqref{2.4}, and this finishes the argument for $L_\sigma$.

At this point we can apply Theorem~\ref{Itsf1} to the operator $L_\sigma$ associated to
$\A_{\sigma}$, and get that its harmonic measure is $A_\infty$ with respect to
the Lebesgue measure on $\R^d$. But this operator is conjugated to $L_\sigma$
by the bilipschitz mapping $g$ (see for instance Lemma 6.17 in \cite{DFM3}), 
so the $A_\infty$ result for $L_\sigma$ follows.

For $\wh L$, we proceed similarly with $\wh A$, write 
$\frac{|t|}{\wh D} =\frac{D_\sigma}{\wh D}\,\frac{|t|}{D_\sigma}$ and use an argument analogous to \eqref{d1-bis}--\eqref{e28-bis} along with Proposition~\ref{Pf1} to conclude.
\ep

\subsection{Further geometric constructions: replacement sets, sawtooth domains, and projections}

To set up the
extrapolation procedure, we rest on the strategy pioneered by \cite{HM}. 
However, a large portion of our work happens on replacement sets rather than saw-tooth domains, and the resulting geometric set up is necessarily different. 

Let us start recalling (and adapting to our scenario) the definitions of the dyadic saw-tooth domains from \cite{HM} and \cite{MZ}. These are morally similar to the domains defined in \eqref{f1}, but a more precise geometric structure will be helpful below. 
Since $\Omega=\RR^n\setminus E$ is an open set, it has a Whitney decomposition -- see Theorem~1 on p. 167 of \cite{S}. 
We will perform the same argument as in  \cite{S}, 
using powers of 10 instead of powers of 2,
and then subpartition emerging cubes further into subcubes of sidelength $10^3$ times smaller. As a result, we get a collection of closed ``Whitney'' boxes in $\Omega$, denoted by 
$\W = \W(\Omega)$, which form a covering of $\Omega$ with pairwise non-overlapping 
interiors and satisfy
\begin{equation}\label{Wbox}
	9\cdot 10^2\,\diam \,I \leq \dist(10^2 I, \po) \leq \dist(I,\po) \leq 21\cdot 10^3\, \diam \,I, \quad \forall \, I\in\W,
\end{equation}
and 
\begin{equation}\label{ngWbox}
	\frac{1}{20} \,\diam \,I_1 \leq \diam\, I_2 \leq 20\,\diam\, I_1
\end{equation}
whenever $I_1$ and $I_2$ in $\W$ touch. Let $X_I$ denote the center of $I$ and $\ell(I) $ the side length of $I$; 
then $\diam\, I \sim \ell(I)$. 

Let $\mathbb{D}$ be a collection of dyadic cubes for the Ahlfors regular set $E$, as in \eqref{a2}--\eqref{a4}. Pick two parameters $\eta\ll 1$ and $K\gg 1$, and for any cube $Q\in\dd$  define
\begin{equation}\label{def:WQ0}
	\W_Q^0 : = \{I\in\W: \eta^{\frac{1}{4}} \ell(Q) \leq \ell(I) \leq K^{\frac{1}{2}} \ell(Q), \dist(I,Q) \leq K^{\frac{1}{2}} \ell(Q)\}.
\end{equation} 
Recall from Definition~\ref{cs} that $A_Q$ denotes a corkscrew point for the surface ball $\Delta(x_Q, C^{-1}r_k(Q))$, with the constant $C$ from  \eqref{a3}. We can guarantee that we can choose $A_Q=X_I$ for some $I\in\W_Q^0$ provided that we choose $\eta$ small enough and $K$ large enough. 

We will further augment our collection to include pertinent Harnack chains. 
To this end, recall the following definition. 

\begin{definition}  \label{def1.hc} We say that an open domain $\Omega$ satisfies the Harnack Chain condition if there is a uniform constant $C$ such that for every $\rho >0,\, \Lambda\geq 1$, 
and every pair of points $X,X' \in \Omega$ with $\delta(X),\,\delta(X') \geq\rho$ 
and $|X-X'|<\Lambda \,\rho$, there is a chain of open balls
$B_1,\dots, B_N \subset \Omega$, $N\leq C(\Lambda)$,
with $X\in B_1$, $X'\in B_N$, $B_k\cap B_{k+1}\neq \emptyset$
and $C^{-1}\diam (B_k) \leq \dist (B_k,\partial\Omega)\leq C\diam (B_k)$.  
Here $C(\Lambda)$ does not depend on $\Omega$, $\rho$, $x$, or $y$.
The chain of balls is called a ``Harnack Chain''.
\end{definition}

The boundary of the domain in the definition 
above is not presumed to exhibit any particular dimension, but we recall that for $\Omega=\RR^n\setminus E$  
for some Ahlfors regular set $E$ of dimension $d < n-1$,
the Harnack chain condition, and even something stronger, holds. 
\begin{lemma}[Lemma 2.1 of \cite{DFM1}]\label{lm:Hcc}
	Let $E$ be a $d$-Ahlfors regular set in $\RR^n$ with $d<n-1$ and $\Omega=\RR^n\setminus E$. Then there exists a constant $c\in(0,1)$, that depends only on $d, n, C_0$, such that for $\Lambda \geq 1$ and $X_1, X_2 \in\Omega$ such that $\delta(X_i) \geq s$ and $|X_1 -X_2| \leq \Lambda s$, we can find two points $Y_i\in B(X_i, s/2) $ such that 
$\dist([Y_1, Y_2],E) 
\geq c\Lambda^{-d/(n-1-d)} s$. That is, there is a thick tube in $\Omega$ that connects the balls $B(X_i, s/2)$.
\end{lemma}
This a stronger property because it ensures 
that two points are connected by a thick tube rather than just a chain, but we did verify that it formally implies the Harnack chain condition from Definition~\ref{def1.hc}, with the constants depending on the ADR constants of $E$ and the dimension only - see Remark~2.2 in \cite{MZ}.
We review some of this for the convenience of the reader.

\begin{remark}\label{rmk:Hcc}
Note that  in the situation above, 
\begin{equation}\label{eq:Hclength}
	|Y_1 - Y_2| \leq |Y_1 - X_1| + |X_1 - X_2| + |X_2 + Y_2| < 2\Lambda s. 
\end{equation} 
Let $\tau = c\Lambda^{-d/(n-1-d)} s$ and $Z_1 = Y_1$. For $2\leq j \leq N$ let $Z_j$ be consecutive points on the line segment $[Y_1, Y_2]$ such that $|Z_j - Z_{j-1}| = \tau/3$. Then
\[ (N-1) \frac{\tau}{3} \leq |Y_1 - Y_2| < N \frac{\tau}{3}. \]
Combined with \eqref{eq:Hclength} we get 
a bound for the length of the chain, namely 
\begin{equation}\label{eq:Hccount}
	N \sim \frac{|Y_1- Y_2|}{\tau/3} \lesssim \Lambda^{\frac{n-1}{n-1-d}}.
\end{equation}
Let $B_0 = B(X_1, s/2), B_j = B(Z_j, \tau/4)$ for $1\leq j\leq N$ and 
$B_{N+1} = B(X_2, s/2) $. Clearly $B_j \cap B_{j+1} \neq \emptyset$ for all $0\leq j\leq N$. Moreover 
$\dist(B_0, E), \dist(B_{N+1}, E) \geq s/2$ 
and for $1\leq j\leq N$,
\begin{equation}\label{eq:Hclb}
\dist(B_j,E) 
\geq \frac{3}{4}\tau = \frac{3}{4} c\Lambda^{-\frac{d}{n-1-d}} s,
\end{equation}
and
\begin{equation}\label{eq:Hcub}
	\dist(B_j,E) 
\leq \min\{\delta(X_1), \delta(X_2)\} + \frac{s}{2} + |Y_1 - Y_2| < \min\{\delta(X_1), \delta(X_2)\} + 3\Lambda s.
\end{equation}
\end{remark}

For each $I\in\W_Q^0$, by Lemma \ref{lm:Hcc} and the discussions after that, 
there is a Harnack chain connecting 
its center
$X_I$ to the corkscrew point  $A_Q$ ; 
we call it $\mathcal{H}_I$. By the definition of $\W_Q^0$ we may construct this Harnack chain so that it consists of a bounded number of balls (depending on the values of $\eta, K$), and stays a distance at least $c\eta^{\frac{n-1}{4(n-1-d)}}\ell(Q)$ away from $\po$ (see \eqref{eq:Hclb}). We let $\W_Q$ denote the set of all $J\in\W$ which meet at least one of the Harnack chains 
$\mathcal{H}_I$, with $I\in \W_Q^0$, i.e.
\begin{equation}
	\W_Q := \{J\in\W: \text{there exists } I\in\W_Q^0 \text{ for which } \mathcal{H}_I \cap J \neq \emptyset\}.
\end{equation}
Clearly $\W_Q^0 \subset \W_Q$. Besides, it follows from the construction of the augmented collections $\W_Q$ and the properties of the Harnack chains 
(in particular \eqref{eq:Hclb} and \eqref{eq:Hcub}) that there are uniform constants 
$c$ and $C$ such that
\begin{equation}\label{def:WQ}
	c\eta^{\frac{n-1}{4(n-1-d)} } \ell(Q) \leq \ell(I) \leq CK^{\frac{1}{2}} \ell(Q), 
	\quad \dist(I,Q) \leq CK^{\frac{1}{2}} \ell(Q) 
\end{equation}
for any $I\in \W_Q$. In particular once $\eta, K$ 
are fixed,
for any $Q\in \dd$ the cardinality of $\W_Q$ is uniformly bounded, 
by an integer
which we denote by $N_0$.

Next we choose a small parameter $\theta\in (0,1)$ so that for any $I\in\W$, the concentric 
dilation $I^* = (1+\theta)I$ still satisfies the Whitney property
\begin{equation}\label{eq:Wboxdl}
	\diam\, I \sim \diam\, I^* \sim \dist(I^*,\po) \sim \dist(I,\po).
\end{equation}
Moreover by taking $\theta$ small enough we can guarantee that 
$\dist(I^*, J^*) \sim \dist(I,J)$ for every $I, J\in\W$, that 
$I^*$ meets $J^*$ if and only if $\partial I$ meets $\partial J$, 
 and that $\frac{1}{2} J\cap I^*=\emptyset$ for any distinct $I, J\in\W$. In what follows we will need to work with further dilations $I^{**} = (1+2\theta) I$ or $I^{***} = (1+4\theta)I$ etc. (We may need to take $\theta$ even smaller to make sure the above properties also hold for $I^{**}, I^{***}$ etc.)
Given an arbitrary $Q\in\dd$, we define associated Whitney regions $U_Q$, $U_Q^*$ by 
\begin{equation}
	U_Q: = \bigcup_{I\in\W_Q} I^*, \quad U_Q^*: = \bigcup_{I\in\W_Q} I^{**}.
\end{equation}

Let $\dd_Q = \{Q'\in\dd: Q' \subset Q\}$. For any $Q\in\dd$ and any family $\mathcal{F}= \{Q_j\}$ of disjoint cubes in $ \dd_Q\setminus\{Q\}$, we define the local and global discretized 
sawtooth regions relative to $\kF$ by
\begin{equation}
	\dd_{\kF,Q} : = \dd_Q \setminus \bigcup_{Q_j\in\kF} \dd_{Q_j}, \quad \dd_{\kF} : = \dd \setminus \bigcup_{Q_j\in\kF} \dd_{Q_j}.
\end{equation} 
We also define the local sawtooth domain relative to $\kF$ by
\begin{equation}\label{def:OFQ} 
	\Omega_{\kF,Q}:= \Int \bigg( \bigcup_{Q'\in\dd_{\kF,Q}} U_{Q'}\bigg), \quad \Omega_{\kF,Q}^*:= \Int \bigg( \bigcup_{Q'\in\dd_{\kF,Q}} U_{Q'}^*\bigg).
\end{equation}
For convenience we set
\begin{equation}\label{def:WFQ}
	\W_{\kF,Q} := \bigcup_{Q'\in\dd_{\kF,Q}} \W_{Q'},
\end{equation}
so that in particular, we may write
\begin{equation}
	\Omega_{\kF,Q} = \Int \bigg(\bigcup_{I\in\W_{\kF,Q}} I^* \bigg), \quad 
	\Omega_{\kF,Q}^* = \Int \bigg(\bigcup_{I\in\W_{\kF,Q}} I^{**} \bigg).
\end{equation}
We will need further fattened sawtooth domain $\Omega_{\kF,Q}^{**}$ etc. whose definitions follow the same lines as above. We remark that by \eqref{def:WQ}, there is a constant $C_3$ depending on $K, \theta$ such that
\begin{equation}\label{stinball}
	\Omega_{\mathcal{F},Q} \subset B(x_Q, C_3 \ell(Q)) \cap \Omega
\end{equation}
for any $Q\in\dd$ and collection of maximal cubes $\mathcal{F}$, where $x_Q$ is the ``center'' of $Q$ as in \eqref{a3}.

The global versions $\Omega_\kF, \Omega^*_\kF, \W_{\kF}$ are defined analogously using $\dd_\kF$ in place of $\dd_{\kF, Q}$.

The sawtooth domains thus defined, of course, have boundaries with portions of different dimension: 
parts of their boundary are given by the intersection with the original $d$-dimensional set $E$, 
and other parts are composed of the $(n-1)$-dimensional faces of Whitney cubes in $\Omega$. 
Yet, they are amenable to the analysis in \cite{DFM-mixed}. 
In particular, $\Omega_\kF$ itself satisfies corkscrew and Harnack chain conditions and $\po_\kF$ can be equipped with a doubling measure $\mu_*$ defined as follows. 
For each Borel set $E\subset\partial\Omega_{\kF}$, let
\begin{equation}\label{eq.nu}
\mu_*(E)=\H^d|_{\Gamma}(E\cap\Gamma)+\int_{E\backslash\Gamma}\dist(X,\Gamma)^{d+1-n}\,d\H^{n-1}|_{\partial\Omega_{\kF}\backslash\Gamma}(X).
\end{equation}
It has been demonstrated in \cite{Bruno} that $\mu_*$ is doubling, and moreover, the domain $\Omega_\kF$ equipped with the measure $m(A)=\dint_A \dist (x, E)^{-n+d+1}\, dX$ and the boundary measure $\mu_*$, satisfies the conditions (H1)--(H6) from \cite{DFM-mixed} (see Theorem~4.2 in \cite{Bruno}). In particular, the harmonic measure corresponding to the operator $L$ is well-defined on $\po_\kF$, and we will denote the latter by $\omega_*$. It satisfies the usual properties of an elliptic measure: doubling, comparison principle. change-of-pole inequalities - we send the reader to \cite{DFM-mixed} for details. We point out that the dyadic cubes in \cite{Bruno} have been built in powers of 2, while ours are built in powers of 10, but otherwise the construction is identical and we will freely use the results from \cite{Bruno}.

As long as $\eta$ and $\theta$ are chosen small enough and $K$ large enough 
depending on the dimension and the AR 
constants of $E$ only, all the properties listed 
above are satisfied with constants depending on the dimension and the AR 
constants of $E$ only, uniformly for all $\kF$ (in particular, because for a $d$-dimensional $E$, $d<n-1$, its own corkscrew and Harnack chain properties are satisfied with the constants depending on the dimension and the AR 
constants of $E$ only). We might adjust the choice of $\eta, \theta$ and $K$ as we go along but we will always make sure that it depends on the dimension and the AR 
constants of $E$ only.

At this point we want to compare the two domains $\Omega_{\kF}$ and $\Omega_{\cF}$
that were constructed above in \eqref{def:OFQ} and \eqref{f1}. 
Retain the notation from Definition~\ref{def-st}, where we are given
a top cube $Q_0$ and a stopping time region $\cF \subset \bD(Q_0)$.
Then we consider the now very specific class $\kF=\{Q_j\}_j$ of 
maximal cubes $Q \subset Q_0$ that are not contained in $\cF$
(we may call them the stopped cubes). Then $\bD_{\kF}$ consists of all the cubes of $\cF$ 
(those that were used to construct $\Omega_\cF$ in \eqref{f1}), and also the cubes $Q$ 
that are not contained in $Q_0$. In other words, 
\begin{equation}\label{eqtwoom1} 
\bD_{\kF, Q_0} = \cF, \quad \mbox{and, respectively,} \quad 
\bD_\kF=(\dd \setminus \bD(Q_0))\cup \cF.
\end{equation}

Ultimately, we will be choosing $\eta, \theta$ and $K$ first and then $M_0$ large enough, depending on $\eta, \theta$ and $K$ so that $U^*(Q)\subset W(Q)$ and hence
\begin{equation}\label{dygeomst}
\Omega_{\kF, Q_0}\subset \Omega^*_{\kF, Q_0}
 \subset \Omega_\cF. \end{equation}

More generally, the geometric statements in the remainder of this subsection implicitly assume that we are allowed to adjust our choices, while keeping their order intact, that is, if $\eta, \theta$ are small enough and $K$ is large enough depending on $d,n,$ and AR constants of $E$; $M_0$ is large enough depending on all these parameters; $M$ is large enough depending on $d,n,$ and AR constants of $E$ and our choice of $\eta, \theta, K, M_0$; and $\eps_0, \delta_0$ are small enough depending on all of the above, then the statements are valid. Since all of these ultimately depend on $d,n,$ and AR constants of $E$, 
we shall suppress this in many 
statements.
We will need to define a ``projection" of cubes $Q_j\in \kF$ on $\Sigma$ and a projection 
of 
$Q_j$ on $\po_\kF$. 
To this end, recall the collection $\cR$ defined in \eqref{f7}--\eqref{f11}. 
We recall that by \eqref{f11} $E\setminus F = \bigcup_{R_k\in \cR} R_k$ and hence, by  \eqref{f18}, 
\begin{equation}\label{eq9.15-0}
E\setminus \bigcup_{R_k\in \cR} R_k 
= F
\subset E\cap \Sigma \cap\overline{Q_0},
\end{equation} 
that is,  $E$ coincides with $\Sigma$ on the complement of $\cup_{R_k\in \cR} R_k$. 

Furthermore, we claim 
that every $R_k\in \cR$ is contained in some $Q_j\in \kF$.
Indeed, if $R_k$ is not contained in any $Q_j\in \kF$ then $R_k\in \dd(Q_0)\setminus \bigcup_{Q_j\in \kF} \dd_{Q_j}$. Then $d_\cF (x_{R_k})\leq l(R_k)$ by definition. However, by \eqref{f10} we have $d_\cF (x_{R_k})\geq (2\tau)^{-1} l(R_k)\geq 50\, l(R_k),$ which is a contradiction.

Having this and \eqref{eq9.15-0} in mind, we can write any $Q_j\in \kF$ as
\begin{equation}\label{eq9.19}
Q_j = \bigg(\bigcup_{R_k\subseteq Q_j} R_k\bigg) \bigcup 
\bigg(Q_j\setminus \bigcup_{R_k\subseteq Q_j} R_k\bigg)\quad \mbox{with}\quad Q_j\setminus \bigcup_{R_k\subseteq Q_j} R_k\subset E\cap \Sigma, 
\end{equation}
and then let 
\begin{equation}\label{eq9.20}
\pi(Q_j): = \bigg(\bigcup_{R_k\subseteq Q_j} \Sigma(k)\bigg) \bigcup \bigg(Q_j\setminus \bigcup_{R_k\subseteq Q_j} R_k\bigg)\subset \Sigma,
\end{equation}
denote a ``projection" of $Q_j$ on $\Sigma$. Here, $\Sigma(k)$ is defined in \eqref{f24}. The definition of $\pi(Q)$ depends on the choice of the numerical constant $\tau$ in the definition of $R_k$ and hence, of $\Sigma(k)$, but this is, as usual, harmless, as long as $\tau$ is small enough for our considerations from previous chapters to apply.

In order to define the projection of $Q_j$ on $\po_\kF$ we recall the following result from \cite{Bruno}.

\begin{proposition}\label{prop.cubeinface} Fix a disjoint family of cubes $\kF\subset\dd$. Then for each $Q_j\in\kF$, there is an $(n-1)-$dimensional cube $P_j\subset\partial\Omega_{\kF}$, which is contained in a face of $I^*$ some $I\in\W$, and which satisfies
	\begin{equation}\label{eq.cubeinface}
	\ell(P_j)\approx\dist(P_j,Q_j)\approx\dist(P_j,\partial\Omega)\approx\ell(I)\approx\ell(Q_j),
	\end{equation}
with uniform constants. 
\end{proposition}
Then we let 
\begin{equation}\label{eq9.20-bis}
\pi_*(Q_j): = P_j
\end{equation}
denote our ``projection" of $Q_j\in \kF$ on $\po_\kF$. 
We first point out the following. 

\begin{lemma}\label{geom-pi-1} Retain the definitions above. Then
\begin{enumerate}
\item
for any $Q_j\in\kF$ the set $\pi(Q_j)$ is contained in some surface ball $\Delta(\wh x_{Q_j}, \wh r_{Q_j})$ where $\wh x_{Q_j} \in \Sigma,$ $\wh r_{Q_j}=C\,l(Q_j)$ for some constant $C$ depending on our choice of $\tau$ only;
\item for any $Q_j\in\kF$ the set $\pi(Q_j)$ contains some surface ball $\Delta(\wh x'_{Q_j}, \wh r'_{Q_j})$, where $\wh x'_{Q_j} \in \pi(Q_j),$ $\wh r'_{Q_j}=c\,l(Q_j)$ for some numerical constant $c$ depending on the AR constant of $\mu$, $d$, and on the choice of $\tau$; 
\item the 
$\pi(Q_j)$, $Q_j\in\kF$, have bounded overlap  and  for any $Q_j\in\kF$ we have $\sigma(\pi(Q_j))\approx \mu(Q_j)$, with all implicit constants depending on the dimension, AR constant of $\mu$, and our choice of $\tau$ only.
\end{enumerate}
Similar statements are valid for $\pi_*$ in place of $\pi$, $\mu_*$ in place of $\sigma$, and $\po_\kF$ in place of $\Sigma$. The relevant constants then depend on our choice of $\eta, \theta, K$ and hence, ultimately, on the dimension and AR constant of $E$ only.
\end{lemma}
\bp We start from the statement that
\begin{equation} \label{eq9.21}
\text{the 
$\pi(Q_j)$ have bounded overlap and $\sigma(\pi(Q_j))\approx \mu(Q_j)$.} 
\end{equation}
The fact that the 
$\pi(Q_j)$ have bounded overlap is a direct consequence of \eqref{f30} and of 
the fact that the 
$Q_j\in \kF$ are disjoint. The equivalence of sizes will follow from statements (1) and (2) of the Lemma. 
Next, \begin{equation} \label{eq9.22}
\text{every $\pi(Q_j)$ is contained in some surface ball $\Delta(\wh x_{Q_j}, \wh r_{Q_j})$,} 
\end{equation}
where $\wh x_{Q_j} \in \Sigma,$ $\wh r_{Q_j}=C\,l(Q_j)$ for some constant $C$ depending on our choice of $\tau$ only. This follows from the  
observation that for every $R_k\subset Q_j$ we have $l(R_k)\leq l(Q_j)$ and from definition \eqref{f24} upon recalling that $M_2=10^{-1} \tau^{-1}$ and $\tau \in (0, 10^{-2})$. 
Thirdly, 
\begin{equation} \label{eq9.22-1}
\text{every $\pi(Q_j)$ contains some surface ball $\Delta(\wh x'_{Q_j}, \wh r'_{Q_j})$,} 
\end{equation}
where $\wh x'_{Q_j} \in \pi(Q_j)$ and 
$\wh r'_{Q_j}=c\,l(Q_j)$ for some numerical constant $c$ depending on the AR 
constant of $\mu$, $d$, and on the choice of $\tau$. First, we show that 
\begin{equation}\label{eq9.22-0}
\mbox{every decomposition in \eqref{eq9.19} contains at least one $R_k$ with 
$l(R_k)\approx l(Q_j)$. }
\end{equation}
Indeed, $Q_j$ contains a surface ball $E\cap B(x_Q, C^{-1}l(Q_j))$ by \eqref{a3}. Here, $C$ is the constant denoted by $C$  in \eqref{a3}. Using now the second inclusion in \eqref{a3} we deduce that there is a dyadic cube $Q'$ contained in $E\cap B(x_Q, C^{-1}l(Q_j)/2)$ with the sidelength $l(Q')=c_1(C, d) l(Q_j)$ for a sufficiently small $c_1(C, d)$ depending on $C$ and $d$ only. For any point $z\in Q'$
\begin{multline*}
d_{\cF}(z)= \inf_{Q \in \cF} \dist(z,Q) + l(Q)\\[4pt]\geq 
\min \Big\{\inf_{Q \in \cF: \,Q_j\subset Q} (\dist(z,Q) + l(Q)) , 
\inf_{Q \in \cF:\, Q_j\cap Q=\emptyset} (\dist(z,Q) + l(Q))\Big\}\\[4pt]
\geq 
\min \left\{l(Q_j),\,C^{-1}l(Q_j)/2\right\}= C^{-1}l(Q_j)/2.
\end{multline*} 
Now, if $C^{-1}/2\geq 1/\tau$ then $C^{-1}l(Q_j)/2\geq C^{-1}l(Q')/2\geq l(Q')/\tau$, and so $Q'$ is a subset of some $R_k$ by definition \eqref{f9}. And recalling that any 
$R_k$ intersecting $Q_j$ has to be contained in $Q_j$, we have $l(R_k)\approx l(Q_j)$ with the implicit constant depending on $C$ and $d$ only. If, on the other hand, $C^{-1}/2\leq 1/\tau$ then we can keep subdividing $Q'$ (the number of times only depending on $\tau$, $C^{-1}$, and $d$), until we reach a cube $Q''\subset Q'$ such that $C^{-1}l(Q_j)/2\geq l(Q'')/\tau$. Then, similarly to above, there must be an $R_k$ containing $Q''$ with $l(Q'')\leq l(R_k)\leq l(Q_j)$ and therefore $l(R_k)$ depending on $\tau$, $C^{-1}$, and $d$ only. 

Finally, once we know that there exists an $R_k\subset Q_j$ with  $l(R_k)\approx l(Q_j)$, it follows from the definition of $\Sigma(k)$ in \eqref{f24} that (possibly slightly enlarging $M_2$) $\Sigma(k)$ contains a surface ball of radius $\approx l(R_k)$, and this surface ball will be used as $\Delta(\wh x'_{Q_j}, \wh r'_{Q_j})$. This finishes the proof of \eqref{eq9.22-1}.

Proving analogous  statements for the projection on $\po_\kF$ is actually much easier: (1) and (2) follow directly from the fact that $P_j$ is a cube, (3) is due to the fact that $P_j$'s have a bounded overlap and the definition of $\mu_*$ (the reader can consult \cite{Bruno} for more general upper and lower estimates on $\mu_*$, but in fact, in this case it is perhaps easier to see the desired bound directly from definitions). 
\ep

Next, for any $Q\in \dd(Q_0)\setminus \cup_{Q_j\in \kF} \dd(Q_j)$ we let 
\begin{multline}\label{eq9.23}
\pi(Q): = \Big(\bigcup_{Q_j\in \kF, \, Q_j\subset Q} \pi(Q_j)\Big) 
\bigcup \Big(Q\setminus \bigcup_{Q_j\in \kF, \, Q_j\subset Q} Q_j\Big)\\[4pt]
= \Big(\bigcup_{R_k\in \cR:\,R_k\subset Q} \Sigma(k)\Big) \bigcup \Big(Q\setminus \bigcup_{R_k\in \cR:\,R_k\subset Q} R_k\Big) \subset \Sigma, 
\end{multline}
and 
\begin{equation}\label{eq9.23-bis}
\pi_*(Q): = \Big(\bigcup_{Q_j\in \kF, \, Q_j\subset Q} \pi_*(Q_j)\Big) \bigcup \Big(Q\setminus \bigcup_{Q_j\in \kF, \, Q_j\subset Q} Q_j\Big) \subset \po_\kF.
\end{equation}
\begin{lemma}\label{geom-pi-2} Retain the definitions above. Then
\begin{enumerate}
\item
for any $Q\in \dd(Q_0)\setminus \cup_{Q_j\in \kF} \dd(Q_j)$ the set $\pi(Q)$ is contained in some surface ball $\Delta(\wh x_{Q}, \wh r_{Q})$ where $\wh x_{Q} \in \Sigma,$ $\wh r_{Q}=C\,l(Q)$ for some constant $C$ depending on our choice of $\tau$ only; 
\item for any $Q\in \dd(Q_0)\setminus \cup_{Q_j\in \kF} \dd(Q_j)$ the ball $\Delta(\wh x_{Q}, 2\,\wh r_{Q})$ above is contained in $\Delta(\wh x_{Q_0}, 4\,\wh r_{Q_0})$; for any $Q_j\in \kF$ the ball $\Delta(\wh x_{Q_j}, 2 \,\wh r_{Q_j})$ from Lemma~\ref{geom-pi-1} is contained in $\Delta(\wh x_{Q_0}, 4\,\wh r_{Q_0})$;
\item for any $Q\in \dd(Q_0)\setminus \cup_{Q_j\in \kF} \dd(Q_j)$ the corkscrew point of $Q$ with respect to $E$, $A_{Q}$, is a corkscrew point for $B(\wh x_{Q}, 4\,\wh r_{Q})\cap \Sigma$ on $\Sigma$ with implicit constants depending on the AR constant of $\mu$, on our resulting choice of $M_0$, and on the choice of $\tau$; 
\item for any $Q\in \dd(Q_0)\setminus \cup_{Q_j\in \kF} \dd(Q_j)$ the set $\pi(Q)$ contains some surface ball $\Delta(\wh x'_{Q}, \wh r'_{Q})$ where $\wh x'_{Q} \in \pi(Q),$ $\wh r'_{Q}=c\,l(Q)$ for some numerical constant $c$ depending on the AR constant of $\mu$, $d$, and on the choice of $\tau$; 
\item for any $Q\in \dd(Q_0)\setminus \cup_{Q_j\in \kF} \dd(Q_j)$ we have $\sigma(\pi(Q))\approx \mu(Q)$, with all implicit constants depending on the dimension, AR constant of $\mu$, $d$, and our choice of $\tau$ only.
\end{enumerate}
The analogues of statements (1)--(3) and (5) are valid for $\pi_*$ in place of $\pi$, $\mu_*$ in place of $\sigma$, and $\po_\kF$ in place of $\Sigma$. 
The relevant constants then depend on our choice of $\eta, \theta, K$ and hence, ultimately, 
on the dimension and AR constant 
of $E$ only. The statement (4) for $\pi_*$ should be substituted by the following property.

For $Q_j\in\kF$, let $B(x^*_j,r^*_j)$ be a ball, 
centered on
$\pi_*(Q_j)=P_j$, 
and such that $r^*_j\approx l(Q_j)$ and 
$$\overline{\bigcup_{Q\in \dd_{Q_j}} U_Q} \subset B(x^*_j,r^*_j). $$ Then for each $Q\in \dd(Q_0)\setminus \cup_{Q_j\in \kF} \dd(Q_j)$, there is a surface ball $\Delta'(x^*_{Q}, r^*_{Q})$ where $x^*_{Q} \in \po_{\kF},$ $r^*_{Q}\approx l(Q)$, and such that 
	\begin{equation}\label{eq.starball}
	\Delta'(x^*_{Q}, r^*_{Q})\subset\Big(Q\cap\partial\Omega_{\kF}\Big)\cup\Big(\cup_{Q_j\in\kF:Q_j\subset Q}B(x^*_j,r^*_j)\cap \po_\kF\Big),
	\end{equation}
	with$\dist(Q,\Delta'(x^*_{Q}, r^*_{Q}))\lesssim\ell(Q)$.
\end{lemma}
\bp 
The fact that 
\begin{equation} \label{eq9.24}
\text{every $\pi(Q)$ is contained in some surface ball $\Delta(\wh x_{Q}, \wh r_{Q})$,} 
\end{equation}
follows from the definition of $\pi(Q)$:  we recall that 
the surface ball
$\Delta(x_Q, l(Q))$ contains $Q$ by \eqref{a3} and then observe that all $\Sigma(k)$ such that  $R_k\subset Q$ are contained $\Delta(x_Q, (1+M_2)l(Q))$. Taking $\wh r_{Q}=2 (1+M_2)l(Q))$ and any $\wh x_{Q}\in \pi(Q)$ we then have \eqref{eq9.24}. 

We note that the construction above and a virtually identical construction in the proof of (1) in Lemma~\ref{geom-pi-1} ensure 
that (2) in the statement of the Lemma is satisfied. 
Moreover, for any $Q\in \dd(Q_0)\setminus \cup_{Q_j\in \kF} \dd(Q_j)$, 
$A^{Q}\in B(x_{Q}, C^{-1}l(Q))\subset B(\wh x_{Q}, 4\,\wh r_{Q})$ and $\dist(A^{Q}, E)\geq \tau_0C^{-1} l(Q)$ by definition. It follows that having chosen $M_0$ large enough depending on the AR constant only, $A^{Q} \in \Omega_\cF$, and hence, by \eqref{ff10}, 
$$\dist(A^{Q}, \Sigma)\geq (1-C\,M_0\eps) \tau_0C^{-1} l(Q).$$
Assuming, as usual, that $\eps$ is small depending on $M_0$ and other parameters depending on the AR constant, we have $\dist(A^{Q}, \Sigma)\geq \tau_0(2C)^{-1} l(Q).$ It follows that $A^{Q}$ is a corkscrew point for $B(\wh x_{Q}, 4\,\wh r_{Q})\cap \Sigma$ with implicit constants depending on the AR constant of $\mu$ and our resulting choices of $M_0$ and smallness of $\eps$, as well as on $\tau$.
 
As a result, $\wh x_{Q} \in \Sigma$ and 
$\wh r_{Q}=C\,l(Q)$ for some constant $C$ depending on our choice of $\tau$ only. 

The opposite inclusion, while less straightforward, is also true: for any $Q\in \dd(Q_0)\setminus \cup_{Q_j\in \kF} \dd(Q_j)$
\begin{equation} \label{eq9.22-4}
\text{$\pi(Q)$ contains some surface ball $\Delta(\wh x'_{Q}, \wh r'_{Q})$,} 
\end{equation}
where $\wh x'_{Q} \in \pi(Q),$ $\wh r'_{Q}=c\,l(Q)$ for some numerical constant $c$ depending on the AR constant of $\mu$, $d$, and on the choice of $\tau$. 
Let us prove this.
Fix some $Q\in \dd(Q_0)\setminus \cup_{Q_j\in \kF} \dd(Q_j)$. Recall \eqref{f11} which, in particular, means that 
$$Q\cap F=Q\cap (Q\setminus [(E\setminus F)\cap Q])= Q\cap \bigg( 
Q\setminus \bigcup_{R_k\in \cR:\,R_k\subset Q} R_k\bigg)\subset \pi(Q).$$

Now, by \eqref{a3} there exists a constant $C_1$ 
such that 
$E\cap B(x_Q, C_1^{-1} l(Q))\subset Q$. 
Fix $\kappa= C_1^{-1}/4\leq  1/4$. If $B(x_Q, \kappa l(Q))\cap \Sigma \subset Q\cap F\subset \pi(Q)$ then we can assign $\Delta(\wh x'_{Q}, \wh r'_{Q}):=B(x_Q, \kappa l(Q))\cap\Sigma$ and finish the argument.

If, on the other hand, there exists a $\wh x_Q \in B(x_Q, \kappa l(Q))\cap\Sigma$ such that $\wh x_Q\in \Sigma\setminus F$, we choose $\Delta(\wh x'_{Q}, \wh r'_{Q})=\wh \Delta:= B(\wh x_Q, \kappa l(Q))\cap\Sigma$. Let us show that $B(\wh x_Q, \kappa l(Q))\cap\Sigma\subset \pi(Q).$ 

Assume, on the contrary, that there exists $\xi\in \wh \Delta$ such that $\xi\not\in \pi(Q)$. First of all, $\wh \Delta\subset B(x_Q, 2\kappa l(Q))$ (since $|\wh x_Q-x_Q|\leq \kappa l(Q)$) and hence, $\wh \Delta\cap F\subset B(x_Q, C_1^{-1}l(Q))\cap F\subset Q\cap F\subset \pi(Q)$ for our choice of $\kappa$. Therefore, $\xi\in (\wh B\cap\Sigma)\setminus F$ and the discussion near \eqref{f36}--\eqref{f37} applies. In particular, there exists a $k_0$ such that $\xi\in \Sigma(k_0)$ and there exists $x\in R_{k_0}$ such that $|x-\xi|\leq C\eps d_{\cF}(\xi)$. However, 
$$d_{\cF}(\xi)  =\inf_{Q'\in \cF} (l(Q')+\dist(\xi, Q'))\leq l(Q)+|\xi-x_Q|\leq (1+2\kappa)l(Q).$$
Therefore, 
$$|x-x_Q|\leq |x-\xi|+|\xi-x_Q|\leq C\eps (1+2\kappa)l(Q)+2\kappa l(Q)<2C\eps l(Q)+2\kappa l(Q). $$
Assuming that $C\eps \leq C_1^{-1}/4$ (which is always safe because $C_1$ depends only on the AR constant of $\mu$) and recalling that $\kappa= C_1^{-1}/4$ (which is our choice), we conclude that 
$x\in E\cap B(x_Q, C_1^{-1} l(Q)) \subset Q.$ Therefore, $R_{k_0}$ meets $Q$. Then $R_{k_0}$ is necessarily a subcube of $Q$ and therefore, $\xi\in \Sigma(R_{k_0})\subset \pi(Q)$, which is a contradiction. We have finished the proof of \eqref{eq9.22-4}.

The last statement of the Lemma follows from (1), (4), and the 
Ahlfors regularity of $\Sigma$.

Let us turn to the analogous statements for the projection on the saw-tooth. (1) remains true since $\pi_*(Q_j)$ is at a distance proportional to $\ell(Q_j)\leq \ell(Q)$ from $Q$. (2) is proved exactly as it is for $ \pi$. (3) is due to the fact that $Q\notin \kF$ and hence, 
$A_Q$ belongs
to some Whitney region associated to $Q$ whose interior (and even its enlargement) is a subset of $\Omega_\kF$. 

Therefore, $A_Q$ is also a corkscrew point for a corresponding surface ball on $\po_\kF$ with the corkscrew constants depending on $\eta$ and $\theta$ in the definition of $\po_\kF$. 
As for (5), 
\begin{multline}
\mu_*(\pi_*(Q)) 
= \mu\bigg(Q\setminus \bigcup_{Q_j\in \kF, \, Q_j\subset Q} Q_j\bigg)
+ \mu_*\bigg(\bigcup_{Q_j\in \kF, \, Q_j\subset Q} \pi_*(Q_j)\bigg)
\\ \approx \mu\bigg(Q\setminus \bigcup_{Q_j\in \kF, \, Q_j\subset Q} Q_j\bigg)
+ \sum_{Q_j\in \kF, \, Q_j\subset Q}\mu_*\big(\pi_*(Q_j)\big)
\\ \approx
\mu\bigg(Q\setminus \bigcup_{Q_j\in \kF, \, Q_j\subset Q} Q_j\bigg)+ \sum_{Q_j\in \kF, \, Q_j\subset Q}\mu( Q_j) = \mu (Q)
\end{multline}
where the first  equality is due to the fact that the corresponding regions are disjoint, the second one uses the finite overlap property of $P_j$, the third one follows from (3) of Lemma~\ref{geom-pi-1}. Finally, \eqref{eq.starball} is proved in \cite{Bruno}. 
\ep

\subsection{Projections and the harmonic measure(s)}

With this setup, we are ready for the following results. The proofs of the two Lemmas below follow closely those of Lemma B.2 and Lemma B.6 in \cite{HM} where analogous results have been established for projections of the harmonic measure on a saw-tooth region starting for an $(n-1)$ dimensional set $E$. Our geometric set-up is, however, different not only because of mixed dimension but also because we need to tie in $\wh\omega$, the  
harmonic measure on $\Sigma$ rather than $E$. 

\begin{lem}\label{l9.25-0} 
Let $E$ be a $d$-dimensional Ahlfors regular set in $\RR^n$ and $\mu$  an Ahlfors regular measure on $E$ and let $\kF=\{Q_j\}_j$ be 
the collection of disjoint cubes in $\dd(Q_0)$ associated to our 
stopping time region $\cF$ as above.
Then the following two statements are valid.

\begin{enumerate} 
\item Assume that $\wh \omega$ is a doubling measure on $\Sigma$. 
Then the projection of $\wh \omega$ on $\F$ within $Q_0\in \dd(E)$ defined as 
\begin{equation} \label{eq9.26}
\P_{\kF} \,\wh\omega(A):= \wh\omega 
\bigg(A\setminus \big(\bigcup_{Q_j\in \kF} Q_j \big)\bigg) +
\sum_{Q_j\in \kF} \frac{\mu(A\cap Q_j)}{\mu(Q_j)}\, \wh\omega(\pi(Q_j)), 
\quad A\subset Q_0,
\end{equation}
is dyadically doubling on $Q_0$. In particular, the conclusion of this lemma  
holds for the elliptic measure $\wh\omega=\wh\omega^{A_{Q_0}}$,  
associated to the operator $\wh L$ on $\RR^n\setminus \Sigma$. 

\item  Assume that $\omega_*$ is a doubling measure on $\Omega_\kF$. 
Then the projection of $\omega_*$ on $\F$ within $Q_0\in \dd(E)$ defined by 
\begin{equation} \label{eq9.26-bis} 
\P_{\kF} \,\omega_*(A)
:= \omega_*\bigg(A\setminus \big(\bigcup_{Q_j\in \kF} Q_j \big)\bigg) 
+\sum_{Q_j\in \kF} \frac{\mu(A\cap Q_j)}{\mu(Q_j)}\,  \omega_*(\pi_*(Q_j)), \quad 
A\subset Q_0,
\end{equation}
is dyadically doubling on $Q_0$. In particular, the conclusion of this lemma 
holds for the elliptic measure 
$\omega_*=\omega_*^{A_{Q_0}}$ associated to the operator $L$ 
on $\Omega_\kF$.
\end{enumerate}
\end{lem}

\bp 
Let $Q\in \dd(Q_0)$ and $Q'\subset Q$ be a dyadic child of $Q$. 
There are three possible cases: $Q$ (and hence $Q'$) is contained in one of the 
$Q_j \in \kF$; 
$Q'$ coincides with one of the $Q_k\subset \kF$ and hence, 
$Q$ is not contained in any 
$Q_j\in \kF$;
$Q'$ (and hence, $Q$) is not contained in any $Q_j\in \kF$.

If $Q$ (and hence $Q'$) is contained in one of the 
$Q_j \in \kF$, then
$$\P_{\kF} \,\wh\omega(Q)
=\frac{\mu(Q\cap Q_j)}{\mu(Q_j)}\, \wh\omega(\pi(Q_j)), 
\quad \P_{\kF} \,\wh\omega(Q')=\frac{\mu(Q'\cap Q_j)}{\mu(Q_j)}\, \wh\omega(\pi(Q_j)), 
$$
and the desired result follows from the doubling property of $\mu$.

The second case is when $Q'$ coincides with one of the 
$Q_k\in \kF$  
and hence, $Q$ is not contained in any of $Q_j\subset \kF$. 
Notice that  $Q$ contains any cube $Q_j\in \kF$ such that  $Q\cap Q_j\neq\emptyset$. 
Since in addition the
$\pi(Q_j)$ have finite overlap and are disjoint from 
$Q\setminus \big(\bigcup_{Q_j\in \kF} Q_j \big)$, 
\begin{multline}\label{eq9.26-1}
\P_{\kF} \wh\omega (Q)
= \wh\omega\bigg(Q\setminus \big(\bigcup_{Q_j\in \kF} Q_j\big)\bigg) +
\sum_{Q_j\in \kF:\,Q\cap Q_j\neq\emptyset} \frac{\mu(Q\cap Q_j)}{\mu(Q_j)}\, \wh\omega(\pi(Q_j))\\[4pt]
=\wh\omega\bigg(Q\setminus \big(\bigcup_{Q_j\in \kF} Q_j\big)\bigg) +
\sum_{Q_j\in \kF:\,Q\cap Q_j\neq\emptyset} \frac{\mu(Q_j)}{\mu(Q_j)}\, \wh\omega(\pi(Q_j))\\[4pt]
= \wh\omega\bigg(Q\setminus \big(\bigcup_{Q_j\in \kF} Q_j\big)\bigg) +
\sum_{Q_j\in \kF:\,Q\cap Q_j\neq\emptyset} \wh\omega(\pi(Q_j)) \lesssim \wh\omega (\pi(Q)).
\end{multline}
We claim that there exists $\wh x_{Q'}\in \pi(Q')$ and $c', C'>0$ such that
\begin{equation}\label{eq9.26-2}
\Delta(\wh x_{Q'}, cl(Q')) \subset \pi(Q'), \quad \pi(Q)\subset  \Delta(\wh x_{Q'}, Cl(Q')).
\end{equation}
The first statement follows from \eqref{eq9.22-1}. The second one follows from \eqref{eq9.24} harmlessly enlarging $C$ so that the ball in \eqref{eq9.24} is contained in $\Delta(\wh x_{Q'}, Cl(Q'))$ for any $\wh x_{Q'}\in \pi(Q')$. Now, recalling that $\wh \omega$ is doubling (here, we have to use (2) and (3) of Lemma~\ref{geom-pi-2} to ensure that the pole is properly placed), we have
\begin{equation}\label{eq9.26-2} \P_{\kF} \wh\omega (Q) \lesssim \wh\omega (\pi(Q)) \leq \wh\omega (\Delta(\wh x_{Q'}, Cl(Q'))) \lesssim \wh\omega (\Delta(\wh x_{Q'}, cl(Q'))) \lesssim \wh\omega (\pi(Q')).
\end{equation}
Since we are in the case
when $Q'$ coincides with one of the 
$Q_k\in \kF$, 
we have, in particular, 
$\P_{\kF} \,\wh\omega(Q') = \wh\omega(\pi(Q_k))=\wh\omega(\pi(Q'))$, so that the right-hand side of \eqref{eq9.26-2} is equal to $\P_{\kF} \,\wh\omega(Q')$. 

Finally, we consider the third case when $Q'$ (and hence, $Q$) is not contained in any $Q_k\in \kF$. 
The very same argument as above, using \eqref{eq9.22-4} in place of \eqref{eq9.22-1}, yields 
\begin{equation}\label{eq9.26-1-bis}\P_{\kF} \wh\omega (Q)\lesssim \wh\omega (\pi(Q')). 
\end{equation}
However, much as in \eqref{eq9.26-1}, 
\begin{equation}\label{eq9.26-3}\P_{\kF} \wh\omega (Q')
= \wh\omega\bigg(Q'\setminus (\bigcup_{Q_j\in \kF} Q_j)\bigg) +
\sum_{Q_j\in \kF:\,Q'\cap Q_j\neq\emptyset} \wh\omega(\pi(Q_j))\gtrsim \wh\omega (\pi(Q')),
\end{equation}
as desired.

Passing to the projection of $\omega_*$, we remark that the treatment of the first 
two cases is literally the same. The only significant difference is the argument for the 
third case when $Q'$ (and hence, $Q$) is not contained in any 
$Q_k\in \kF$.
Much as for $\wh \omega$ we have 
$$ \P_{\kF}\omega_* (Q) \lesssim \omega_* (\pi_*(Q)) \lesssim \omega_* (\Delta (x^*_Q, r_Q^*)),$$
where the first inequality follows from the argument analogous to \eqref{eq9.26-1} and the second one is due to Lemma~\ref{geom-pi-2}, a version of (1) for $\po_\kF$, and the surface ball $\Delta (x^*_Q, r_Q^*)$, $x_Q^*\in \po_\kF$, $r_Q^*\approx l(Q)$, is the one containing $\pi_*(Q)$. Using the doubling property of $\omega_*$ we have
$$\omega_* (\Delta (x^*_Q, r_Q^*))\lesssim \omega_*(\Delta'(x_{Q'}^*, r_{Q'}^*)), $$
where $\Delta'(x_{Q'}^*, r_{Q'}^*)$ is the surface ball from \eqref{eq.starball} for the cube $Q'$. Now, 
by \eqref{eq.starball} and the doubling property of $\omega_*$,
\begin{multline}\label{eqnew1}
\omega_*(\Delta'(x_{Q'}^*, r_{Q'}^*)) 
\leq \omega_*\left(Q'\cap\partial\Omega_{\kF}\right)
+ \sum_{Q_j\in\kF:Q_j\subset Q'} \omega_*\big(B(x^*_j,r^*_j)\cap \po_\kF\big)
\\ \lesssim \omega_*\big(Q'\cap\partial\Omega_{\kF}\big)
+ \sum_{Q_j\in\kF:Q_j\subset Q'} \omega_*\big(\pi_*(Q_j)\big) 
\lesssim \omega_*(\pi_*(Q')).
\end{multline}
Then, analogously to \eqref{eq9.26-3}, we can finish the argument. 

We remark that the elliptic measure of $L$ on $\po_\kF$ fits the hypothesis and, in particular, is doubling by \cite{Bruno, DFM-mixed}.

\ep

\begin{lem}\label{l9.25} 
Let $E$ be a $d$-dimensional Ahlfors regular set in $\RR^n$ and $\mu$ be an Ahlfors regular measure on $E$. Under the conditions of Theorem~\ref{t8.8}, the projection of the harmonic measure of $\wh L$ on $\F$ within $Q_0\in \dd(E)$ defined by \eqref{eq9.26} with  
$\wh\omega=\wh\omega^{A_{Q_0}}$
is $A^\infty_{\dd}(Q_0)$ with respect to $\mu$.
\end{lem}
\bp As per Remark~\ref{rAinfty}, we aim to show \eqref{AiTh-gen-2} 
for some $\delta, \epsilon\in (0,1)$.  

For brevity, we will write $\wh\omega=\wh\omega^{A_{Q_0}}$ throughout the proof. We fix $0<\eta<1/2$ and $A\subset Q\in \dd(Q_0)$ with $\mu(A)\geq (1-\eta) \mu(Q)$. 

If $Q\subseteq Q_j$ for some $Q_j\subset \kF$, then 
$$ \frac{\P_{\kF} \wh \omega (A)}{\P_{\kF} \wh \omega (Q)}=\frac{\mu(A\cap Q_j)}{\mu(Q_j)}\, \wh\omega^{A_{Q_0}}(\pi(Q_j))\left(\frac{\mu(Q\cap Q_j)}{\mu(Q_j)}\, \wh\omega^{A_{Q_0}}(\pi(Q_j))\right)^{-1}=\frac{\mu (A)}{\mu(Q)} \geq (1-\eta),$$
as desired.

If $Q$ is not contained in any cube of $\kF$, then it belongs to $\dd(Q_0)\setminus\cup_{Q_j\in \kF} \dd(Q_j)$. It might or might not intersect with the cubes of $\kF$ and we let 
$$\kF_G=\{Q_j\in \kF: \,\mu(A\cap Q_j)\geq (1-2\eta) \, \mu(Q_j)\},$$
and 
$$E_0:=Q\setminus \bigcup_{Q_j\in \kF} Q_j, \quad G:= \bigcup_{Q_j\in \kF_G} Q_j, \quad B:=\bigcup_{Q_j\in \kF\setminus \kF_G} Q_j.$$
Then simply dropping the cubes of $\kF\setminus \kF_G$ in the sum,

\begin{multline}
\label{eq9.27}
\cP_{\kF} \wh\omega (A) \geq \omega (A\cap E_0) + \sum_{Q_j\in \kF_G} \frac{\mu(A\cap Q_j)}{\mu(Q_j)}\,\wh\omega (\pi(Q_j)) \\[4pt]
\geq \omega(A\cap E_0) + (1-2\eta)\sum_{Q_j\in \kF_G} \wh\omega (\pi(Q_j))\geq (1-2\eta)\,\wh\omega ((A\cap E_0) \cup \pi(G)),
\end{multline}
where we set 
$$\pi(G):=\bigcup_{Q_j\in \F_G} \pi(Q_j),$$
and used the fact that $\pi(Q_j)$ and $E_0$ are disjoint. 
Recall now \eqref{eq9.23}--\eqref{eq9.24}. Using the property that $\wh\omega$ 
is $A_\infty$ with respect 
to $\sigma$ by Theorem~\ref{t8.8}, we have 
\begin{equation}\label{eq9.28}
\frac{\wh\omega ((A\cap E_0) \cup \pi(G))}{\wh\omega(\Delta(\wh x_Q, \wh r_Q))} \gtrsim \left(\frac{\sigma((A\cap E_0) \cup \pi(G))}{\sigma(\Delta(\wh x_Q, \wh r_Q))}\right)^\theta.
\end{equation}
Here again we have to use (2) and (3) of Lemma~\ref{geom-pi-2} to ensure that the pole is properly placed.

Now, from \eqref{eq9.21} and the fact that the 
$Q_j$ are disjoint we conclude that $\sigma(\pi(G))\approx \mu(G)$ and since furthermore  $A\cap E_0$ and $\pi(G)$ as well as $A\cap E_0$ and $G$ are disjoint, 
\begin{equation}\label{eq9.29}
 \sigma((A\cap E_0) \cup \pi(G)) \approx \mu((A\cap E_0) \cup  G).  
 \end{equation}
Also,
$$
\mu(A\cap B) \leq (1-2\eta) \sum_{Q_j\in \kF\setminus \kF_G}\mu(Q_j)
\leq  (1-2\eta)\mu(Q),
$$
and hence, 
$$(1-\eta)\mu(Q) \leq \mu(A) \leq \mu((A\cap E_0)\cup G)+\mu (A\cap B)\leq \mu((A\cap E_0)\cup G)+(1-2\eta)\mu(Q),$$
so that 
\begin{equation}\label{eq9.30}
 \mu((A\cap E_0)\cup G) \geq \eta\, \mu(Q).
\end{equation}
Since by Ahlfors regularity of $\sigma$ 
$$ \sigma(B(\wh x_Q, \wh r_Q)) \approx \wh r_Q^d\approx \mu(Q),$$
invoking \eqref{eq9.29} and \eqref{eq9.30} we 
see that the right-hand side of \eqref{eq9.28} is
bounded from below by $\eta^\theta$, modulo a multiplicative constant. 

Coming back to \eqref{eq9.27}, this yields
$$ \cP_{\kF} \wh\omega (A)\gtrsim (1-2\eta) \eta^\theta\, \wh\omega(\Delta(\wh x_Q, \wh r_Q)).$$ 
Finally, we recall that the ball $\Delta(\wh x_Q, \wh r_Q)$ contains $\pi(Q)$ by \eqref{eq9.23}--\eqref{eq9.24} and the sets $\pi(Q_j)$ have have finite overlap and are disjoint from $Q\cap E_0$. Therefore,  
$$ \wh\omega(\Delta(\wh x_Q, \wh r_Q)) \geq \wh\omega(\pi(Q)) \gtrsim \wh\omega(Q\cap E_0)+\sum_{Q_j\in \F} \wh\omega(\pi(Q_j)) =  \P_{\kF} \wh\omega(Q),$$
as desired. 
\ep

\begin{theorem}\label{l9.39}
Let $E$ be a $d$-dimensional  Ahlfors regular set and $\mu$ be an Ahlfors regular 
measure on $E$. Let $\omega$ be the harmonic measure associated to the operator 
$L=-\div D_{\mu}^{-(n-d-1)} \nabla$ in $\RR^n\setminus E$,  with 
\begin{equation} \label{eq9.40}
D_{\mu}(X) = \Big\{ \int_E |X-y|^{-d-\alpha} d\mu(y) \Big\}^{-1/\alpha}, \quad \alpha>0.
\end{equation}
Then let $Q_0\in \dd(E)$ be given, construct $\cF=\cF_{\eps_0, \delta_0}(Q_0)$
and the complementary collection $\kF$ as in Definition \ref{def-st}, assume that
$\cF \neq \emptyset$, and define the projection of $\omega=\omega^{A_{Q_0}}$ 
on $\F$ within $Q_0$ by
\begin{equation} \label{eq9.41}
\P_{\kF} \,\omega(A):= \omega\bigg(A\setminus (\bigcup_{Q_j\in \kF} Q_j)\bigg) +
\sum_{Q_j\in \kF} \frac{\mu(A\cap Q_j)}{\mu(Q_j)}\, \omega(Q_j), \quad A\subset Q_0.
\end{equation}
If our various constants are chosen correctly (see below this statement), 
$\P_\kF \omega$ lies in the class $A^\infty_{\dd}(Q_0)$ with respect to $\mu$,
with $A_\infty$ constants that depend only on $n, d, C_0$, and $\alpha$.
\end{theorem}

To be precise, if 
$\eta$ (see above \eqref{def:WQ0})and $\theta$ (near \eqref{eq:Wboxdl}) 
are small enough and $K$ (above \eqref{def:WQ0}) is large enough,
depending on $n, d,C_0$; 
$M_0>1$ is large enough depending on $n, d,C_0$ and $\eta, \theta, K$; 
$M>1$ large enough depending on $n, d, C_0, M_0$, 
and $\eps_0, \delta_0> 0$ small enough depending on $n, d, C_0, M_0, M$, then $\P_\kF \omega$ lies in the class $A^\infty_{\dd}(Q_0)$ with respect to $\mu$,
with $A_\infty$ constants that depend only on $n, d, C_0$, and $\alpha$.

As usual, the uniform rectifiability is not needed at this stage. 

\smallskip

\bp 
The plan of the proof is to show that $\P_{\kF} \,\omega_*^{A_{Q_0}}$ 
(where $\omega_*$ is the elliptic measure associated to the operator $L$ on $\Omega_\kF$)
is $A^\infty_{\dd}(Q_0)$ with respect to $\P_{\kF} \,\wh \omega^{A_{Q_0}}$ 
(associated to the operator $\wh L$ on the  domain $\R^n \sm \Sigma$) 
and $\P_{\kF} \,\omega_*^{A_{Q_0}}$ is $A^\infty_{\dd}(Q_0)$ with respect to 
$\P_{\kF} \,\omega^{A_{Q_0}}$, so that $\P_{\kF} \,\omega^{A_{Q_0}}$ is $A^\infty_{\dd}(Q_0)$ with respect to $\P_{\kF} \,\wh \omega^{A_{Q_0}}$ (in the notation of Lemma~\ref{l9.25}) 
at which point we can use Lemma~\ref{l9.25} to achieve the desired result. 
This second part is easier, and is closely related to the Main Lemma for sawtooth projections in \cite{DJK},
its version in \cite{HM}, and similar results. However, working with a ``replacement boundary" $\Sigma$ and the associated harmonic measure is new, and the first part, requiring both a change of the domain and a change of the operator, is more intricate.

We remark that formally speaking, we only defined $A^\infty$ and $A^\infty_\dd$ properties with respect to the Ahlfors regular measure on the boundary, but the same 
Definitions~\eqref{Ainfty-gen}, \eqref{Ainfty-bis} apply to any doubling measure $\mu$ together with the equivalent reformulations, in particular, in Remark~\ref{rAinfty}.

We start with the proof of the 
simpler statement that $\P_{\kF} \,\omega_*^{A_{Q_0}}$ is $A^\infty_{\dd}(Q_0)$ 
with respect to $\P_{\kF} \,\omega^{A_{Q_0}}$.

As usual, we will simplify the notation by writing $\omega=\omega^{A_{Q_0}}$, $\omega_*=\omega_*^{A_{Q_0}}$, and $\wh\omega= \wh \omega^{A_{Q_0}}$ throughout the proof. Recall from Remark~\ref{rAinfty} that 
$A^\infty_\dd$ is an equivalence relationship, and let us concentrate on showing that 
there is a constant 
$C>0$  such that for every $Q\in \dd(Q_0)$ and every Borel set $A\subset Q$, 
\begin{equation} \label{eq9.42}
\frac{\P_\kF \omega_*(A)}{\P_\kF \omega_*(Q)}\leq C\, \frac{\P_\kF  \omega(A)}{\P_\kF  \omega(Q)}.
\end{equation}

The simplest case is when $Q\subseteq Q_j$ for some $Q_j\in \kF$. Then, by definition, 
$$
\frac{\P_\kF \omega(A)}{\P_\kF  \omega(Q)}
=\frac{\P_\kF \omega_*(A)}{\P_\kF  \omega_*(Q)} = \frac{\mu(A\cap Q_j)}{\mu(Q\cap Q_j)} .
$$

Let us assume now that $Q$ is not contained in any $Q_j\in \kF$. In this case, similarly to \eqref{eq9.26-3} and \eqref{eqnew1}
\begin{equation}\label{eq9.43}
\P_{\kF} \omega_* (Q)\gtrsim \omega_* (\pi(Q))\geq \omega_*(\Delta'(x_{Q}^*, r_{Q}^*)),
\end{equation}
where $\Delta'(x_{Q}^*, r_{Q}^*)$ is the surface ball from \eqref{eq.starball}.

 Therefore, with the notation $E_0=Q_0 \setminus \big(\bigcup_{Q_j\in \kF} Q_j \big)$,
\begin{equation}\label{eq9.44}
\frac{\P_\kF  \omega_*(A)}{\P_\kF \omega_*(Q)}\lesssim  \frac{  \omega_*(A\cap E_0)}{\omega_* (\Delta'(x_{Q}^*, r_{Q}^*))}+ 
\sum_{Q_j\in \kF: \, Q_j\subset Q} \frac{\mu(A\cap Q_j)}{\mu(Q_j)}\, \frac{\omega_*(\pi_*(Q_j))}{\omega_* (\Delta'(x_{Q}^*, r_{Q}^*))}.
\end{equation}
Using the change of pole formula for $\omega_*$ in \cite{DFM-mixed} and, if necessary, 
Harnack's inequality to slightly adjust the corkscrew point, we can write the above as 
\begin{equation}\label{eq9.45}
\frac{\P_\kF  \omega_*(A)}{\P_\kF  \omega_*(Q)}\lesssim \omega_*^{A_Q}(A\cap E_0) + 
\sum_{Q_j\in \kF: \, Q_j\subset Q} \frac{\mu(A\cap Q_j)}{\mu(Q_j)}\,  \omega_*^{A_Q}(\pi_*(Q_j)).
\end{equation}
By the maximum principle and the fact that $\omega^{X}(Q_j)\approx 1$ for $X\in \pi_*(Q_j)$, we have

$$\omega_*^{A_Q}(A\cap E_0)\lesssim  \omega^{A_Q}(A\cap E_0)\quad\mbox{and} \quad \omega_*^{A_Q}(\pi_*(Q_j))\lesssim \omega^{A_Q}(Q_j).$$
 Then 
\begin{equation}\label{eq9.45-bis}
\frac{\P_\kF  \omega_*(A)}{\P_\kF  \omega_*(Q)}\lesssim \omega^{A_Q}(A\cap E_0) + 
\sum_{Q_j\in \kF: \, Q_j\subset Q} \frac{\mu(A\cap Q_j)}{\mu(Q_j)}\,  \omega^{A_Q}( Q_j) \lesssim \frac{\P_\kF  \omega(A)}{\P_\kF  \omega(Q)},
\end{equation}
as desired.

Now let us turn to the proof that $\P_{\kF} \,\omega_*^{A_{Q_0}}$ 
(where $\omega_*$ is the elliptic measure associated to the operator $L$
on $\Omega_\kF$)
is $A^\infty_{\dd}(Q_0)$ with respect to $\P_{\kF} \,\wh \omega^{A_{Q_0}}$ 
(associated to the operator $\wh L$ on the  domain $\R^n \sm \Sigma$). The challenge is to change the operator and the domain on which the harmonic measure is evaluated simultaneously. 

To this end, we recall that by Remark~\ref{rAinfty} (which as we mentioned, applies to general doubling measures) in order to show that a doubling measure $\omega$ is $A^\infty_\dd(Q_0)$ with respect to another doubling measure $\nu$, it is sufficient to show that there exists $0<\eps, \delta<1$ such that for every $Q\in \dd(Q_0)$ and every Borel set $F\subset Q$, 
\begin{equation} \label{AiTh-gen-2-bis}
\frac{\omega(F)}{\omega(Q)} < \delta \Rightarrow \frac{\nu(F)}{\nu(Q)} < \varepsilon.
\end{equation}

We claim that it is moreover sufficient to show an even weaker property, where we
only check something like \eqref{AiTh-gen-2-bis} in the middle of $Q$ (but in a uniform way).
First we claim that if we choose the constant $M'$ large enough,
then for each $Q\in \dd(Q_0)$, we can find a dyadic subcube $Q' \subset Q$ such that 
\begin{equation}\label{eq9.77.1}
\frac {1}{10 M'} \ell(Q) \leq \ell(Q') \leq \frac{1}{M'} \ell(Q), 
\quad \mbox{and}\quad  \dist (Q', E \sm Q)  
\geq \frac {M'}{2} \ell(Q').
\end{equation}
Indeed, \eqref{a3} gives a ``center'' $x_Q$ for $Q$ such that 
$\dist(x_Q), E \sm Q) \geq C^{-1} \ell(Q)$, and if $M'$ is large enough, any
cube $Q'$ such that $\frac {1}{10 M'} \ell(Q) \leq \ell(Q') \leq \frac{1}{M'} \ell(Q)$
(this covers a generation of cubes) and that contains $x_Q$ will satisfy the
second part of \eqref{eq9.77.1} automatically. 

Now pick any $M'$ as above, and for each $Q\in \dd(Q_0)$ a cube 
$Q' = Q'(Q) \subset Q$ such that \eqref{eq9.77.1} holds; 
we claim that in order to show that a doubling measure $\omega$ is $A^\infty_\dd(Q_0)$ with respect to another doubling measure $\nu$ it is enough to prove that there exists 
$0<\eps', \delta'<1$ such that for every $Q\in \dd(Q_0)$ and every Borel set $F\subset Q'$, 
\begin{equation} \label{AiTh-gen-2-bis2}
\frac{\omega(F)}{\omega(Q')} < \delta' \Rightarrow \frac{\nu(F)}{\nu(Q')} < \varepsilon'.
\end{equation}

Now assume that  
\eqref{AiTh-gen-2-bis2} holds for every $F \subset Q' \subset Q$  
as stated,  
and that for some $F\subset Q$ we have $\frac{\nu(F)}{\nu(Q)} > \varepsilon$,  with $\eps \in (0,1)$  
to be defined shortly. Obviously  
$$
\nu(Q\setminus Q')+\nu(F\cap Q') = \nu(F) > \eps \nu(Q).
$$

Since $\nu$ is doubling, there exists a small constant $c_{\nu}$ 
depending on the doubling constants of $\nu$, the dimension, and $M'$, 
such that  $\nu(Q)\leq \frac{1}{c_\nu}\, \nu (Q')$.     Hence, 
$$ 
\frac{\nu(F\cap Q')}{\nu(Q')} > \frac{\varepsilon\nu(Q') - \nu(Q\setminus Q')}{\nu(Q')}
= 1 - (1-\eps) \frac {\nu(Q)}{\nu(Q')} 
\geq 1 - (1-\eps) \frac{1}{c_\nu}.   
$$
Choosing $0<\eps<1$ so that $c_\nu(1-\eps')=1-\eps$  
(which is always possible, for any given $0<\eps'<1$), we get $\frac{\nu(F\cap Q')}{\nu(Q')}>\eps'$. Then, by \eqref{AiTh-gen-2-bis2}, 
$\frac{\omega(F\cap Q')}{\omega(Q')} \geq \delta' $. Therefore, 
$$\frac{\omega(F)}{\omega(Q)}\geq c_\omega \frac{\omega(F\cap Q')}{\omega(Q')} \geq c_\omega \delta' =:\delta, $$
where $c_\omega$ is a constant depending on the doubling constants of $\omega$, 
the dimension, and $M'$ only. We arrive at \eqref{AiTh-gen-2-bis}.

Hence, in application to our situation, it is enough to show that there exists $0<\eps, \delta<1$ 
such that for every $Q\in \dd(Q_0)$ and every Borel set $F\subset Q'$, 
\begin{equation} \label{AiTh-gen-2-bis2}
\frac{\P_{\kF} \,\omega_*^{A_{Q_0}}(F)}{\P_{\kF} \,\omega_*^{A_{Q_0}}(Q')} > \varepsilon \Rightarrow
\frac{\P_{\kF} \,\wh \omega^{A_{Q_0}}(F)}{\P_{\kF} \,\wh \omega^{A_{Q_0}}(Q')} > \delta 
\end{equation}
where the large constant $M'$ will  be chosen below, and $Q'$
is a descendent of $Q$ chosen as above.

If $Q'$ is a subcube of some $Q_j\in \kF$, there is nothing to prove, for 
$$
\frac{\P_\kF \wh \omega^{A_{Q_0}}(F)}{\P_\kF  \wh \omega^{A_{Q_0}}(Q')} 
=\frac{\P_\kF \omega_*^{A_{Q_0}}(F)}{\P_\kF  \omega_*^{A_{Q_0}}(Q')} 
= \frac{\mu(F\cap Q_j)}{\mu(Q\cap Q_j)} .
$$
Therefore, we concentrate on the case when $Q'$ and hence, $Q$, 
belongs to $\dd_\kF$. Since 
$$
\P_\kF  \wh \omega^{A_{Q_0}}(Q') \lesssim  \wh \omega^{A_{Q_0}}(\pi(Q')) 
$$
(see \eqref{eq9.26-1}, \eqref{eq9.26-1-bis}), the usual change of pole considerations yield 
\begin{multline*} 
\frac{\P_{\kF} \,\wh \omega^{A_{Q_0}}(F)}{\P_{\kF} \,\wh \omega^{A_{Q_0}}(Q')} 
\gtrsim  \wh\omega^{A_{Q'}} 
\Big(F\setminus \big(\bigcup_{Q_j\in \kF} Q_j \big)\Big) +
\sum_{Q_j\in \kF} \frac{\mu(F\cap Q_j)}{\mu(Q_j)}\, \wh\omega^{A_{Q'}}(\pi(Q_j)) \\ 
\gtrsim  \wh\omega^{A_{Q'}} 
\Big(F\setminus \big(\bigcup_{Q_j\in \kF} Q_j \big)\Big) +
\sum_{Q_j\in \kF:\, \frac{\mu(F\cap Q_j)}{\mu(Q_j)}>\eta} \frac{\mu(F\cap Q_j)}{\mu(Q_j)}\, \wh\omega^{A_{Q'}}(\pi(Q_j)) \\
\gtrsim  \wh\omega^{A_{Q'}} 
\Big(F\setminus \big(\bigcup_{Q_j\in \kF} Q_j \big)\Big) + \eta
\sum_{Q_j\in \kF:\, \frac{\mu(F\cap Q_j)}{\mu(Q_j)}>\eta} \wh\omega^{A_{Q'}}(\pi(Q_j))\\
\gtrsim \eta \, \wh\omega^{A_{Q'}} \Big(\bigl(F\setminus \big(\bigcup_{Q_j\in \kF} Q_j \big)\bigr) \bigcup 
\big(\bigcup_{Q_j \in \kF_G} \pi(Q_j) \big)\Big).
\end{multline*}
Here $0<\eta<1$ is a parameter to be chosen below, and $\kF_G$ is the 
collection of cubes in $\kF$ such that $\frac{\mu(F\cap Q_j)}{\mu(Q_j)}>\eta$. 
Let $Z = \bigl(F\setminus \big(\bigcup_{Q_j\in \kF} Q_j \big)\bigr) \bigcup \big(\bigcup_{Q_j \in \kF_G} \pi(Q_j) \big)$
be the set from the last line.
By the maximum principle and the fact that 
$L=\wh L$ in $\Omega_\kF \cap \Omega_\cF$, 
$\wh\omega^{A_{Q'}}(Z)$ is 
larger than the solution to $Lu=0$ 
on $\Omega_\kF \cap \Omega_\cF$ 
with the data given by 
the restriction of $\1_Z$ to 
the boundary of $\Omega_\kF \cap \Omega_\cF$, evaluated at $A_{Q'}$. 

We claim that for some constants $C, \alpha>0$ depending on the usual geometric parameters only, 
the latter is larger than 
$\omega_*^{A_{Q'}}(Z) - C\,(M')^{-d+1-\alpha}$.
In order to prove this, by the maximum principle we only need to show that 
\begin{equation}\label{eq9.80}
\omega_*^{X} \left(\pi(Q') \right) \leq C (M')^{-d+1-\alpha}  
\quad \mbox{for } X\in \partial (\Omega_\kF \cap \Omega _\cF) \setminus \partial\Omega_{\kF} =
\Omega_\kF \cap \partial\Omega _\cF, 
\end{equation}
because $Z \subset \pi(Q')$ since $F \subset Q'$. 

Let us postpone for now the proof of \eqref{eq9.80} and try to finish the argument. At this point, collecting all of the above, we have 
\begin{equation}\label{eq9.79}
\frac{\P_{\kF} \,\wh \omega^{A_{Q_0}}(F)}{\P_{\kF} \,\wh \omega^{A_{Q_0}}(Q')} 
\gtrsim \eta \, 
\omega_*^{A_{Q'}}(Z) 
-\eta\,(M')^{-d+1-\alpha}. 
\end{equation}
However, 
\begin{multline}\label{eq9.81}
\omega_*^{A_{Q'}} (Z)
\gtrsim \omega_*^{A_{Q'}}\bigl(F\setminus \big(\bigcup_{Q_j\in \kF} Q_j \big)\bigr) +\sum_ {Q_j \in \kF_G} \omega_*^{A_{Q'}} (\pi(Q_j))\\
\geq  \omega_*^{A_{Q'}}\bigl(F\setminus \big(\bigcup_{Q_j\in \kF} Q_j \big)\bigr) +\sum_ {Q_j \in \kF_G} \frac{\mu(F\cap Q_j)}{\mu(Q_j)}\,\omega_*^{A_{Q'}} (\pi(Q_j))\\
\geq \omega_*^{A_{Q'}}\bigl(F\setminus \big(\bigcup_{Q_j\in \kF} Q_j \big)\bigr) +\sum_ {Q_j \in \kF} \frac{\mu(F\cap Q_j)}{\mu(Q_j)}\,\omega_*^{A_{Q'}} (\pi(Q_j)) -\eta \sum_ {Q_j \in \kF\setminus \kF_G} \omega_*^{A_{Q'}} (\pi(Q_j))\\
\geq \P_{\kF}\omega_*^{A_{Q'}}(F) -\eta.
\end{multline}
where the first inequality is due to the finite overlap property of $\pi(Q_j)$, 
the second one is simply due to the fact the density does not exceed 1, while  
in the third one we have added back and subtracted the cubes where $\frac{\mu(F\cap Q_j)}{\mu(Q_j)}\leq \eta$, and  
the fourth one uses once again the finite overlap property and the fact that 
$$\sum_ {Q_j \in \kF\setminus \kF_G} \omega_*^{A_{Q'}} (\pi(Q_j)) \leq  \omega_*^{A_{Q'}} (\pi(Q'))\leq 1.$$

Then \eqref{eq9.79} and \eqref{eq9.81} give 
\begin{equation} \label{9a83}
\frac{\P_{\kF} \,\wh \omega^{A_{Q_0}}(F)}{\P_{\kF} \,\wh \omega^{A_{Q_0}}(Q')} 
\gtrsim \eta \left( \P_{\kF}\omega_*^{A_{Q'}}(F) -\eta-(M')^{-d+1-\alpha}\right).
\end{equation}

We started with the assumption that 
$\frac{\P_{\kF} \,\omega_*^{A_{Q_0}}(F)}{\P_{\kF} \,\omega_*^{A_{Q_0}}(Q')} > \varepsilon$
as in \eqref{AiTh-gen-2-bis2}. 
After the change the poles  \eqref{eq9.45} 
yields $\P_{\kF}\omega_*^{A_{Q'}}(F) \geq C^{-1} \varepsilon$,   where of course 
$C$ does not depend on $\eta$. Choosing $\eta$ and $M'$ so that $(M')^{-d+1-\alpha}=C^{-1}\eps/4$, $\eta=C^{-1}\eps/4$, we get $\eqref{AiTh-gen-2-bis2}$ with $\delta\approx \eps^2$. Note that  $\eps \in (0,1)$ does not have to be small. In fact,  to be more accurate,  there are some a priori constraints on the choice on $M'$ (depending on the allowable geometric parameters only) so we will first choose $M'$ and then $\eps$ and $\eta$ accordingly (recall that all of the above are in our disposal).

To finish the proof, it remains to show \eqref{eq9.80}. This is essentially  the comparison principle once we observe that 
\begin{equation}\label{eq9.85}\dist(Q', \Omega_\kF\cap \po_\cF) \geq C M' \ell(Q'),
\end{equation}
where $C$ depends on $\eta, \theta$, and $K$, and the 
dimension and AR constants of $E$ only. This estimate was the entire reason for introducing $M'$ into the argument. To this end, we observe that by \eqref{dygeomst}
$$
\Omega_\kF\cap \po_\cF \subset \cup_{\wt Q\in \dd \setminus \dd(Q_0)} U_{\wt Q}^*. 
$$

We split into two cases. Assume first that $\wt Q\in \dd \setminus \dd(Q_0)$ is such that $\ell(\wt Q)\geq c_0\,\ell(Q_0)$, with a small constant $c_0$ to be determined below. Then 
\begin{multline}\label{eq9.82}
\dist ( U_{\wt Q}^*, Q') \geq \dist ( U_{\wt Q}^*, E) \geq  C(\eta, \theta) \ell(\wt Q) \\
\geq c_0\,C(\eta, \theta)\,\ell(Q_0) \geq c_0\,C(\eta, \theta)\,M'\,\ell(Q'),
\end{multline}
where we used \eqref{def:WQ} and \eqref{eq:Wboxdl} for the third inequality above and \eqref{eq9.77.1} along with the fact that $\ell(Q_0)\geq \ell(Q)$ for the last one. 

If, on the other hand, $\wt Q\in \dd \setminus \dd(Q_0)$ is such that 
$\ell(\wt Q)\leq c_0\,\ell(Q_0)$, with 
$c_0<1$,  then  
$\wt Q \cap Q_0=\emptyset$ and  
\begin{eqnarray} \label{eq9.83}
\dist ( U_{\wt Q}^*, Q') &\geq& \dist (\wt Q, Q')-\dist (U_{\wt Q}^*, \wt Q) 
\geq \dist(Q', E\sm Q)  
-C(K, \theta)\,\ell(\wt Q) \nn\\
&\geq& \frac{M'}{2}\ell(Q') -C(K, \theta)\,\ell(\wt Q) \geq \frac{M'}{2}\ell(Q') -C(K, \theta)\,c_0\,\ell(Q_0)
\\ &\geq& \frac{M'}{2}\ell(Q') -10\,C(K, \theta)\,c_0\, M' \ell(Q'),
\nn
\end{eqnarray}
 where we used \eqref{def:WQ} for the third and the fifth inequality and \eqref{eq9.77.1} 
for the fourth one. Now choosing $c_0$ so that $10\,C(K, \theta)\,c_0=1/4$, 
and combining \eqref{eq9.82}  
and \eqref{eq9.83},  
we arrive at \eqref{eq9.85}, with a small constant $C$ depending on $\eta, \theta, K$, the  
dimension, and the AR constant  
of $E$ only.

With \eqref{eq9.85} at hand, by Lemma~15.28 of \cite{DFM-mixed} 
$$ \omega_*^X(\pi(Q')) \leq C \frac{m(B_{Q'} \cap \Omega_\kF)}{l(Q')^2} g_*(X,A^{Q'}),  
 $$
 where $g_*$ is the Green function of $L$ on $\Omega_\kF$, $B_{Q'}$ is a ball of radius $Cl(Q')$  centered in $Q'$, and the measure $m$ is given by 
 $m(A)=\dint_A \dist (X, E)^{-n+d+1}\, dX$.  
 Next, according to Lemma~14.83 in \cite{DFM-mixed}
$$g_*(X,A^{Q'}) \lesssim \ell(Q')^\alpha \frac {|X-A_{Q'}|^{2-\alpha}}{m(B(X, |X-A_{Q'}|)} \lesssim  \ell(Q')^\alpha \frac {(M'\ell(Q'))^{2-\alpha}}{m(B(A_{Q'}, CM'\ell(Q'))\cap \Omega_\kF)}.
 $$ 
Combining the two estimates above,  
$$ \omega_*^X(\pi(Q')) \leq C  \frac {(M')^{2-\alpha}}{m(B(A_{Q'}, CM'\ell(Q'))\cap \Omega_\kF)} \,m(B_{Q'} \cap \Omega_\kF) =C (M')^{-d+1-\alpha},  
 $$ 
 as desired.
\ep

\section{Extrapolation}\label{Extrap}
Let us start with the following definitions of dyadic Carleson measures. 
\begin{definition}\label{cm-dyad}
Let $E$ be a $d$-dimensional  Ahlfors regular set, $\mu$ be an Ahlfors regular measure on $E$, and $\dd(E)$ be our usual collection of dyadic cubes on $E$ (associated to $\mu$). Let 
$\{\ka(Q)\}_{Q\in \dd(E)}$ be a sequence of non-negative numbers indexed by $Q\in \dd(E)$. For any subcollection $\dd'\subset \dd(E)$, $Q_0\in \dd(E)$ we let
$$\km(\ka,\dd'):=\sum_{Q\in \dd'} \ka(Q)^2\,\mu(Q),$$
and 
$$\|\km(\ka)\|_{CM_\dd}:=\sup_{Q\in \dd(E)} \frac{\km(\ka,\dd(Q))}{\mu(Q)}, \quad  \|\km(\ka)\|_{CM_{\dd(Q_0)}}:=\sup_{Q\in \dd(Q_0)} \frac{\km(\ka,\dd(Q))}{\mu(Q)},$$ 
and if the latter two quantities are finite, we say that $\km(\ka)\in CM_\dd$ or $\km(\ka) \in CM_{\dd (Q_0)}$, respectively. 
\end{definition}

These definitions pertain to the measure (rather than to the sequence) and could seem different in homogeneity from their continuous analogue in geometric saw-tooth in Definition~\ref{cm}. To reconcile these differences we say, slightly abusing the notation, that the sequence $\{\ka(Q)\}_{Q\in \dd(E)} \in CM_\dd$ if the corresponding $\km(\ka)\in CM_\dd$ and similarly $\{\ka(Q)\}_{Q\in \dd(Q_0)} \in CM_{\dd(Q_0)}$ if $\km(\ka)\in CM_{\dd(Q_0)}$.

Furthermore, for any family of pairwise disjoint dyadic cubes $\kF=\bigcup_j Q_j$ we define the restriction of $\km$ on the sawtooth by 
$$\km_\kF(\ka,\dd')=\sum_{Q\in \dd'\setminus\bigcup_{Q_j\in \kF}\dd(Q_j)} \ka(Q)^2\,\mu(Q), $$
and for any $Q\in \dd(E)$
$$\dd_Q^{short}:=\dd(Q)\setminus\{Q\}. $$

As the reader may correctly guess, we aim to use the forthcoming Lemmas for Tolsa's $\alpha$-numbers in place of $\ka(Q)$. At this point, however, we keep the statements in full generality and note that throughout Section~\ref{Extrap} $\{\ka(Q)\}_{Q\in \dd(E)}$ denotes any sequence with non-negative entries.

\begin{lemma}\label{l10.2} Let $E$ be a $d$-dimensional  Ahlfors regular set, $\mu$ be an Ahlfors regular measure on $E$. Fix some $Q_0\in \dd(E)$ and some sequence of non-negative numbers $\{\ka_Q\}_{Q\in \dd(Q_0)}$ such that the corresponding $\km(\ka)$ satisfies 
\begin{equation}\label{eq10.5} \km(\dd(Q_0))\leq (a_0+b_0)\, \mu(Q_0),\quad \mbox{for some }a_0\geq 0, \,b_0>0. 
\end{equation}
Fix some $K\geq 1$ 
and construct a (maximal) family $\kF$ of pairwise disjoint cubes obtained by subdividing $Q_0$  and stopping when 
\begin{equation}\label{eq-st-gen}
\mbox{either 
$\ka(Q)^2 > 2 b_0$
or } J_\ka(Q) = \sum_{k(Q_0) \leq k \leq k(Q)}  \ka(R_k(Q))^2 \geq 2 K b_0, 
\end{equation} 
(at which point we assign $Q\in \kF$). Then
\begin{equation}\label{eq10.3}
\|\km_\kF\|_{CM_{\dd(Q_0)}} \leq 4 K b_0,
\end{equation}
and 
\begin{equation}\label{eq10.4}
\mu (B) \leq \frac{a_0+b_0}{a_0+2b_0}\,\mu(Q_0),
\end{equation}
where $B$ is the union of cubes $Q_j \in \kF$ such that $\km (\ka,\dd_{Q_j}^{short})>a_0\mu (Q_j).$
\end{lemma}

The Lemma is analogous to Lemma~7.2 in \cite{HM}. However, we have to slightly change the stopping time region being constructed and carefully track emerging constants as ultimately only for very special stopping time regions will we be able to use the results of Sections~\ref{intro}-\ref{sSum}.

\bp If $\ka(Q_0)^2\geq 2b_0$ then the result is trivial for the following reason. 
We stop immediately with $Q_0$, so $\kF = \{ Q_0 \}$,
The left-hand side of \eqref{eq10.3} is simply equal to 0, and 
\begin{multline*}\km (\ka,\dd_{Q_0}^{short})= \km (\ka,\dd_{Q_0})-\ka(Q_0)^2 \mu(Q_0)\leq (a_0+b_0)\mu(Q_0)-2b_0\mu(Q_0)\\[4pt]
=(a_0-b_0)\mu(Q_0)< 
a_0\mu (Q_0), 
\end{multline*}
so that $B=\emptyset$.

Therefore, it is safe to assume from now on that $\ka(Q_0)^2< 2b_0$ 
and so $\kF\neq \{Q_0\}.$ As usual, we write $\kF=\{Q_j\}_j$,
where the $Q_j$ are thus disjoint (by maximality) cubes $Q_j\in \dd(Q_0)$.
Then 
\begin{multline}\label{eq10.7}
\sum_{Q_j\in \kF} (\km (\ka,\dd_{Q_j}^{short}) + J_\ka (Q_j) \,\mu(Q_j)) \\[4pt]
= \sum_{Q_j\in \kF} \km (\ka,\dd_{Q_j}) - \sum_{Q_j\in \kF} \ka(Q_j)^2\,\mu(Q_j) + \sum_{Q_j\in \kF} \sum_{k(Q_0) \leq k \leq k(Q_j)}  \ka(R_k(Q_j))^2 \mu(Q_j)
\\[4pt]
= \sum_{Q_j\in \kF} \km (\ka,\dd_{Q_j}) + \sum_{Q_j\in \kF} 
\sum_{R\in \dd(Q_0):\, Q_j \subset R, Q_j \neq R}   \ka(R)^2 \mu(Q_j)
\\[4pt]
= \sum_{Q_j\in \kF} \km (\ka,\dd_{Q_j}) + \sum_{Q\in \dd_{Q_0}\setminus \bigcup_{Q_j\in\kF} \dd_{Q_j}} \ka(Q)^2 \mu(Q)
\sum_{Q_j\in\kF:\, Q_j\subset Q, Q_j \neq R}   
\frac{\mu(Q_j)}{\mu(Q)}
\\[4pt]
\leq  \sum_{Q_j\in \kF} \km (\ka,\dd_{Q_j}) + \sum_{Q\in \dd_{Q_0}\setminus \bigcup_{Q_j\in\kF} \dd_{Q_j}} \ka(Q)^2 \mu(Q) 
\leq \km(\ka, \dd_{Q_0})
\end{multline}
where we used the fact that $Q_j$ are disjoint
for the last line. 
Let $\kF_B:=\{Q_j\in \kF:\, \km (\ka,\dd_{Q_j}^{short})>a_0\mu (Q_j)\}.$  Now, 
\begin{multline}\label{eq10.8}
(a_0+2b_0)\,\mu(B)=(a_0+2b_0)\sum_{Q_j\in \kF_B} \mu(Q_j) =a_0\sum_{Q_j\in \kF_B} \mu(Q_j)+2b_0\sum_{Q_j\in \kF_B} \mu(Q_j)\\[4pt]
< \sum_{Q_j\in \kF_B} \km(\ka, \dd_{Q_j}^{short}) + \sum_{Q_j\in \kF_B} J_\ka(Q_j)\,\mu(Q_j),
\end{multline}
where we used the definition of $\kF_B$ and the fact that $Q_j$ is a stopping cube, hence, either $J_\ka(Q_j)\geq 2K b_0\geq 2b_0$ or 
$J_\ka(Q_j)\geq \ka(Q_j)^2 \geq 2 b_0$, 
so that in any case $J_\ka(Q_j)\geq  2 b_0$. Next, using \eqref{eq10.7} and \eqref{eq10.5}, the last expression in \eqref{eq10.8} is bounded from above by 
\begin{equation}\label{eq10.9}
\sum_{Q_j\in \kF_B} \left(\km(\ka, \dd_{Q_j}^{short}) + J_\ka(Q_j)\,\mu(Q_j)\right) \leq \km(a, \dd_{Q_0}) \leq (a_0+b_0)\,\mu(Q_0),
\end{equation}
and \eqref{eq10.8}--\eqref{eq10.9} yields \eqref{eq10.4}.

Turning to \eqref{eq10.3}, we observe that the latter amounts to showing that 
for every $Q\in \dd(Q_0)$,  
$$ \sum_{Q'\in \dd(Q)\setminus \bigcup_{Q_j\in\kF} \dd_{Q_j}} \ka(Q')^2 \,\frac{\mu(Q')}{\mu(Q)} \leq 4Kb_0.$$ 
Having fixed any such $Q$, it is convenient to introduce,
for any large integer $N$, the collection $\kF_N$ of maximal cubes (by inclusion) of
$$\kF\cup\{Q'\in \dd(Q_0):\,l(Q')\leq 10^{-N-1}\}$$
and
the corresponding smaller family of cubes
$$
{\mathcal H}_Q = 
\dd(Q)\setminus \bigcup_{Q_k\in\kF_N} \dd_{Q_k}
=\Big\{Q'\in \dd(Q)\setminus \bigcup_{Q_j\in\kF} \dd_{Q_j}:\,l(Q')\geq 10^{-N}\Big\}.$$
Clearly, it is sufficient to prove that for every $Q\in \dd(Q_0)$
\begin{equation}\label{eq10.10}
 \sum_{Q'\in {\mathcal H}_Q}   \ka(Q')^2 \,\mu(Q') \leq 4Kb_0\, \mu(Q)
\end{equation}
uniformly in $N$. The main difference between $\kF$ and $\kF_N$ is that the (disjoint by construction) cubes of $\kF_N$ cover any $Q\in \dd(Q_0)$, which implies that 
\begin{multline}\label{eq10.11}
\sum_{Q'\in {\mathcal H}_Q}
\ka(Q')^2 \,\mu(Q')
= \sum_{Q'\in {\mathcal H}_Q}
\sum_{Q_k\in\kF_N:\,Q_k\in \dd(Q')} \frac{\mu(Q_k)}{\mu(Q')}\,\ka(Q')^2 \,\mu(Q')\\[4pt]
= \sum_{Q_k\in\kF_N:\,Q_k\in \dd(Q)}\mu(Q_k)
\sum_{Q'\in \dd(Q): \,Q_k \subset Q', Q_k \neq Q'} \ka(Q')^2 . 
\end{multline}
We separately consider the elements of $\kF_N$ which belong to $\kF$ and those that don't. First, 
\begin{multline}\label{eq10.12}
\sum_{Q_k\in\kF_N\setminus \kF:\,Q_k\in \dd(Q)}\mu(Q_k) 
\sum_{Q'\in \dd(Q): \,Q_k \subset Q', Q_k \neq Q'} 
\ka(Q')^2 
\\[4pt]
\leq \sum_{Q_k\in\kF_N\setminus \kF:\,Q_k\in \dd(Q)}\mu(Q_k) J_\ka (Q_k)
<2Kb_0 \sum_{Q_k\in\kF_N\setminus \kF:\,Q_k\in \dd(Q)}\mu(Q_k)\leq 2Kb_0 \,\mu(Q).
\end{multline}
Here in the next-to-the-last inequality we used the fact that by the stopping time construction $J_\ka (Q_k)<2Kb_0$ when $Q_k\in\kF_N\setminus \kF$, 
for otherwise it would belong to $\F$. 
Next consider $Q_k\in\kF_N \cap \kF$, denote by $\wt Q_k$ denotes the parent of $Q_k$,
observe that $\sum_{Q'\in \dd(Q): \,\,Q_k \subset Q', Q_k \neq Q'} \ka(Q')^2 
\leq J_\ka (\wt Q_k) < 2Kb_0$ because otherwise $\wt Q_k \in \F$ and this would contradict
the maximality of $Q$. We are also using the fact that $Q_k \neq Q_0$ here.
Now
\begin{multline}\label{eq10.12}
\sum_{Q_k\in\kF_N\cap \kF:\,Q_k\in \dd(Q)} \mu(Q_k) 
\sum_{Q'\in \dd(Q): \,Q_k \subset Q', Q_k \neq Q'} 
\ka(Q')^2 
\\[4pt]
\leq \sum_{Q_k\in\kF_N\setminus \kF:\,Q_k\in \dd(Q)}\mu(Q_k) J_\ka (\wt Q_k)
<2Kb_0 \sum_{Q_k\in\kF_N\setminus \kF:\,Q_k\in \dd(Q)}\mu(Q_k)\leq 2Kb_0 \,\mu(Q),
\end{multline}

This finishes the proof of \eqref{eq10.3} and Lemma \ref{l10.2} follows. \ep

\ms
The next result is the main extrapolation step, 
analogous to Lemma~8.5 from \cite{HM}. We have to state it differently, however, because once again we can only afford to work with very special stopping time regions. 

\begin{lemma}\label{l10.14} 
Let $E$ be a $d$-dimensional  Ahlfors regular set and $\mu$ be an Ahlfors regular measure on $E$. Fix some $Q_0\in \dd(E)$ and a dyadically doubling Borel measure $\omega$ on $Q_0$. Assume that there is some sequence of non-negative numbers $\{\ka_Q\}_{Q\in \dd(Q_0)}$ such that the corresponding $\km(\ka)$ satisfies 
\begin{equation}\label{eq10.15}
\|\km(\ka)\|_{CM_{\dd(Q_0)}} \leq L_0 
\end{equation}
for some $L_0<\infty$. Furthermore, assume that there exists $b_0\geq 0$ 
such that for  
some $K\geq 1$, and any
$a_0\in [0, L_0]$ the stopping time region from the statement of Lemma~\ref{l10.2},  
built according  
to \eqref{eq-st-gen}, 
 satisfies the property that the projection of $\omega$ on $\kF$ within $Q_0$, defined by \eqref{eq9.41}, is $A^\infty_\dd (Q_0)$ with respect to $\mu$. Finally, assume that whenever $\km(\ka, \dd(Q))=0$ we have that $\omega$ is $A^\infty_\dd (Q)$ with respect to $\mu$.

Then $\omega$ is $A^\infty_\dd (Q_0)$ with respect to $\mu$.
\end{lemma}

\bp The proof can be carried out closely following that of Lemma~8.5 from \cite{HM}. 
There the authors have a seemingly stronger hypothesis that there exists $\gamma\geq 0$ such that for every $Q\in \dd(Q_0)$ and family of pairwise disjoint dyadic cubes $\kF=\{Q_j\}_j\subset \dd(Q)$ such that $\|\km_\kF\|_{CM_Q}\leq \gamma$ the projection of $\omega$ on $\kF$ within $Q_0$, defined by \eqref{eq9.41}, is $A^\infty_\dd (Q_0)$ with respect to $\mu$. The actual proof, however, relies only on the stopping time regions built in their analogue of our Lemma~\ref{l10.2}. 

In a few words, the proof proceeds by induction argument with a continuous parameter, with the main hypothesis
\begin{multline}\label{eq10.16}
\mbox{$\exists\, \eta_a\in (0, 1), \, \exists\, C_a<\infty: \, \forall\, Q\in \dd(Q_0)$ with $\km(\ka,\dd_Q)\leq a\mu(Q)$ we have }\\
A\subset Q, \quad \frac{\mu(F)}{\mu(Q)}\geq 1-\eta_a\quad\Longrightarrow \quad \frac{\omega(F)}{\omega(Q)}\geq \frac{1}{C_a}, 
\end{multline}
referred to as $H(a)$, $a\geq 0$. 

The induction proceeds in two steps. Step I is that $H(0)$ holds. This is a straightforward consequence of one of our assumptions as $\km(\ka, \dd_Q)=0$ implies that $\omega$ is $A^\infty_\dd (Q)$ with respect to $\mu$. 

Turning to the induction step, one aims to show that $H(a)$ implies $H(a+b_0)$, for all $a\in [0,L_0]$ so that the conclusion of the theorem 
can be reached in $k$ steps where $k$ is such that $kb_0\geq L_0$. To this end, we fix $0\leq a:=a_0\leq L_0$ and $Q\in \dd(Q_0)$ such that $\km(\ka,\dd_Q)\leq (a_0+b_0)\mu(Q)$. In the notation of \cite{HM} one would take $\gamma=4Kb_0$. Then, using the results of Lemma~\ref{l10.2}, the proof of the induction  step follows the lines of the argument for Lemma~8.5 from \cite{HM}, and we omit the details.
\ep

\section{Conclusion}\label{SConcl} 

At this point we finally collect the results of Sections~\ref{intro}--\ref{Extrap} towards the proof of the Main Theorem. 
\begin{theorem}\label{t11.1} Let $E$ be a $d$-dimensional  uniformly rectifiable set in $\RR^n$, $d\leq n-2$,  and $\mu$ be a uniformly rectifiable measure on $E$. Let $\omega$ be the harmonic measure associated to the operator $L=-\div D_{\mu}^{-(n-d-1)} \nabla$ in $\RR^n\setminus E$,  with 
\begin{equation} \label{eq9.40}
D_{\mu}(X) = \Big\{ \int_E |X-y|^{-d-\alpha} d\mu(y) \Big\}^{-1/\alpha}, \quad \alpha>0.
\end{equation}
Then $\omega$ is $A^\infty$ with respect to $\mu$ in the sense of Definition~\ref{Ainfty-bis}.
\end{theorem} 
\bp Our first task is to show that for every fixed $Q_0\in \dd(E)$ the harmonic measure $\omega=\omega^{A_{Q_0}}$ is $A^\infty_\dd(Q_0)$ with respect to $\mu$. Once this is established, we recover that by Harnack inequality we also have $\omega=\omega^{A_{Q_0}}$ is $A^\infty_\dd(Q'_0)$ with respect to $\mu$ for any $Q'_0$ with $l(Q'_0)=l(Q_0)$ and $\dist(Q_0, Q'_0) \leq C l(Q_0)$. And then, using the doubling property of $\omega$ and Harnack inequality once again, we can show that $\omega$ is $A^\infty$ with respect to $\mu$ in the sense of Definition~\ref{Ainfty-bis} (not only dyadically).

Now, as before, we let $\{\alpha(Q)\}_Q\in \dd(E)$ stand for the Tolsa $\alpha$ numbers and recall that by Theorem~\ref{tUR-alpha} and our assumptions there exists $L_0<\infty$ such that 
$$\|\km(\alpha)\|_{CM_{\dd(Q_0)}}\leq L_0,$$
for any fixed $Q_0\in \dd(E)$. We need to verify the rest of the assumptions of Lemma~\ref{l10.14}. 

To this end, recall the statement of Theorem~\ref{l9.39}. Take, as in Theorem~\ref{l9.39}, $M_0>1$ large enough depending on $n, d,C_0$, $M>1$ large enough depending on $n, d, C_0, M_0$, and $\eps_0, \delta_0>0$ small enough depending on $n, d, C_0, M_0, M$, for any $Q_0\in \dd(E)$ so that the conclusion of the Theorem is verified for  $\cF=\cF_{\eps_0, \delta_0}(Q_0)$ and the  
complementary collection $\kF$ built in Definition~\ref{def-st}. 

We can safely assume that $\delta_0>\eps_0^2$, because we decided to chose
$\varepsilon_0$ last. We choose $b_0=\eps_0^2/2$
and $K=\delta_0/\eps_0^2$ so that $2Kb_0=\delta_0$. 
The stopping time region from Lemma~\ref{l10.2} is then the same as the stopping time region from Theorem~\ref{l9.39} (provided that $\kF\neq \{Q_0\}$) and hence, the desired property that the projection of $\omega$ on $\kF$ within $Q_0$, defined by \eqref{eq9.41}, is $A^\infty_\dd (Q_0)$ with respect to $\mu$, is verified. 

If it happens that $\kF=\{Q_0\}$, then by definition \eqref{eq9.41} $\cP_\kF\omega (A) =\frac{\mu(A)}{\mu(Q)}\,\omega^{A_{Q_0}}(Q_0)\approx \mu(A)$, so that the hypothesis of Lemma~\ref{l10.14} ($\cP_\kF\omega$ is $A^\infty_\dd(Q_0)$ with respect to $\mu$) is trivially valid.

Finally, it remains to consider the case $\km(\ka, \dd(Q))=0$ which we will reformulate as  $\km(\ka, \dd(Q_0))=0$. Then by definition $\alpha(Q)=0$ for all $Q\in \dd(Q_0)$ and hence, $J_\alpha(Q)=0$ for all $Q\in \dd(Q_0)$. Therefore, for any  $\eps_0,\delta_0>0$ we have $\kF=\emptyset$ so that $\cP_\kF \omega=\omega$ and Theorem~\ref{l9.39} then gives the desired result.
\ep

\Addresses

\end{document}